\documentclass[10pt]{article}
  \usepackage[T1]{fontenc}
  \usepackage{lmodern}

\usepackage{amsmath,amssymb,mathtools}
  \usepackage{amsthm}
  \usepackage[margin=1in]{geometry}
  \usepackage{graphicx}
\usepackage{microtype}
\usepackage{needspace}
\usepackage{bm}
\usepackage{xcolor}
\usepackage{cases}
\usepackage{booktabs}
\usepackage{multirow}
  \usepackage[numbers]{natbib}
\usepackage{hyperref}
  \hypersetup{
    hidelinks,
    pdftitle={Auxiliary Gradient-Flow Solvers for Generalized Newtonian Models},
    pdfauthor={Stefano Zampini, Daniele Boffi, Gurt Dovletov, Peter Markowich}
  }

  \theoremstyle{plain}
\newtheorem{theorem}{Theorem}[section]
\newtheorem{lemma}[theorem]{Lemma}
\newtheorem{remark}[theorem]{Remark}
\newtheorem{assumption}[theorem]{Assumption}
\newtheorem{definition}[theorem]{Definition}

\def\mybegproof{\noindent{\bf Proof: }}
\def\myendproof{\par\nobreak\rightline{\vrule height5pt width5pt depth0pt}\medskip}

\newcommand{\dx}{\,\mathrm{d}x}
\newcommand{\ds}{\,\mathrm{d}s}
\newcommand{\R}{\mathbb{R}}
\newcommand{\Rinf}{(-\infty,+\infty]}
\newcommand{\grad}{\nabla}
\newcommand{\diver}{\nabla\cdot}
\newcommand{\eps}{\varepsilon}
\DeclareSymbolFont{CMeletters}{OML}{cmm}{m}{it}
\DeclareMathSymbol{\epsilon}{\mathord}{CMeletters}{"0F}

\newcommand{\fro}{\!:\!}
\newcommand{\dd}[2]{\frac{d#1}{d#2}}

\newcommand{\bu}{{\bm{u}}}
\newcommand{\ba}{\bm{a}}
\newcommand{\bv}{\bm{v}}

\newcommand{\btau}{\bm{\tau}}

\newcommand{\Ltwo}{{L^2(\Omega)}}
\newcommand{\Linf}{{L^\infty(\Omega)}}
\newcommand{\Lphi}{{L^\phi(\Omega)}}

\newcommand{\LPhip}{{L_0^\Phi(\Omega)}}

\newcommand{\Wphiz}{{W_0^{1,\phi}(\Omega)}}

\newcommand{\vWphiz}{{W_0^{1,\phi}(\Omega;\mathbb{R}^d)}}
\newcommand{\LPhi}{{L^\Phi(\Omega)}}

\newcommand{\Lpp}{{L^{p'}(\Omega)}}

\newcommand{\Honez}{{H_0^1(\Omega)}}

\newcommand{\vHonez}{{H_0^1(\Omega;\mathbb{R}^d)}}

\newcommand{\geni}[3]{{{#1}^{#2}_{#3}}}
\newcommand{\smphi}{\geni{S}{-}{\phi}}
\newcommand{\spphi}{\geni{S}{+}{\phi}}
\newcommand{\simphi}{\geni{R}{-}{\phi}}
\newcommand{\sipphi}{\geni{R}{+}{\phi}}

\newcommand{\pos}[1]{\left(#1\right)_+}

\newcommand{\Ee}{{\mathcal{E}_\epsilon}}
\newcommand{\Eh}{{\mathcal{E}_h}}
\newcommand{\wEe}{\widetilde{\Ee}}

\newcommand{\DEh}{{\mathcal{DE}_h}}
\newcommand{\Jc}{\mathcal{J}_c}
\newcommand{\Ju}{\mathcal{J}_u}
\newcommand{\Jbu}{\mathcal{J}_\bu}
\newcommand{\wJc}{\widetilde{\mathcal{J}}_c}
\newcommand{\wJu}{\widetilde{\mathcal{J}}_u}
\newcommand{\wJbu}{\widetilde{\mathcal{J}}_\bu}
\newcommand{\J}{\mathcal{J}}

\newcommand{\Ce}{{L_\epsilon^\Phi(\Omega)}}
\newcommand{\Le}{{L_\epsilon^2(\Omega)}}
\newcommand{\dom}[1]{{\text{dom}(#1)}}

\DeclareMathOperator*{\argmin}{argmin}

\DeclareMathOperator*{\essinf}{ess\,inf}
\DeclareMathOperator*{\esssup}{ess\,sup}

\newcommand{\st}{:}
\newcommand{\kac}{{Ka{\v{c}}anov }}

\renewcommand{\P}[1]{\mathcal{P}_{#1}}
\newcommand{\Pk}{\P{k}}
\newcommand{\Pkc}{\P{2k-2}}
\newcommand{\Pkm}{\P{k-1}}
\newcommand{\Q}[1]{\mathcal{Q}_{#1}}
\newcommand{\Qk}{\Q{k}}

\newcommand{\Uh}{U_h}

\newcommand{\bUh}{\bm{U}_h}

\newcommand{\bUhdiv}{\bm{U}_{h,\mathrm{div}}}

\newcommand{\Ch}{C_h}

\newcommand{\auc}{\zeta}
\newcommand{\pauc}{\pi}
\newcommand{\bauc}{\boldsymbol{\zeta}}
\newcommand{\evar}{\delta}
\usepackage{authblk}
  
  \newcommand{\authororcid}[2]{\href{https://orcid.org/#2}{#1}}

\begin{document}
  \title{Auxiliary Gradient-Flow Solvers for Generalized Newtonian Models}
  \date{}
\author[1,*]{\authororcid{Stefano Zampini}{0000-0002-0435-0433}}
\author[1,2,3]{\authororcid{Daniele Boffi}{0000-0003-0717-2940}}
\author[1]{\authororcid{Gurt Dovletov}{0009-0002-2607-6655}}
\author[4]{Peter Markowich}
\affil[1]{Computer, Electrical and Mathematical Sciences and Engineering\\King Abdullah University of Science and Technology}
\affil[2]{Dipartimento di Matematica ``F. Casorati'', University of Pavia}
\affil[3]{IMATI ``E. Magenes'', CNR, Pavia}
\affil[4]{Faculty of Mathematics, University of Vienna, Vienna}
\affil[*]{Corresponding author: \href{mailto:stefano.zampini@gmail.com}{stefano.zampini@gmail.com}}

\newcommand{\paperabstract}{%
We introduce an auxiliary gradient-flow framework for variational problems with generalized Newtonian structure governed by an N-function. The key idea is to replace the nonlinear constitutive dependence on the gradient, or symmetric gradient, by an auxiliary scalar variable representing its squared magnitude. This shifts the nonlinearity from the state equation to the auxiliary variable, yielding a sequence of uniformly elliptic weighted linear problems.
At the continuous level, we construct an auxiliary energy on a metric space adapted to the growth of the underlying N-function. In this topology, we prove lower semicontinuity, geodesic $\lambda$-convexity, and exponential convergence of the associated minimizing-movement scheme. At the finite element level, we derive a metric gradient flow through an explicit Riesz map, prove global well-posedness of the resulting semi-discrete ODE, and establish convergence to the finite element solution of the Euler--Lagrange equations of the generalized Newtonian energy. For the $p$-Laplacian and $p$-Stokes models, this gives a rigorous convergence result for $4/3\le p\le 4$, $p\ne2$, with asymptotic rate estimates beyond this range.
We also propose practical time discretizations, including an operator-splitting scheme that gives the \kac iteration as a special case, and an adaptive pseudo-transient method that can be implemented using scalable linear solvers. Numerical experiments for power-law, Carreau--Yasuda, regularized Bingham, and optimal-design models demonstrate robustness and performance competitive with Newton's method, with essentially mesh-independent iteration counts in the power-law and Bingham tests.
}
\newcommand{\paperkeywords}{Gradient flows; Generalized Newtonian models; Finite element methods; \kac iteration; Pseudo-transient continuation; p-Laplacian; p-Stokes; Carreau--Yasuda model.}
  \maketitle
  \begin{abstract}
  \paperabstract
  \end{abstract}

\section{Introduction}

This work concerns the numerical solution of minimization problems of the type
\begin{equation}\label{eq:energy_general}
\argmin_{u\in W}\int_{\Omega} \phi(|Du|)\dx - F(u),
\end{equation}
where $\Omega\subset\R^d$, $d=2,3$, is an open, bounded domain,
the integrand is governed by an N-function $\phi$ \cite{musielak1983orlicz}, and $D$ denotes either the gradient of a scalar field
or the symmetric gradient $\eps(\bu) = \frac{1}{2}(\grad\bu + \grad\bu^T)$ of an incompressible vector field in an Orlicz--Sobolev functional space $W$. Such energies arise in nonlinear diffusion,
optimal design, and in the modeling of generalized Newtonian fluids. The associated Euler--Lagrange equations
are nonlinear and may be singular or degenerate, which presents challenges for both analysis and numerical
approximation. %

The finite element approximation of problems of this type has a long history, beginning with classical works on the $p$-Laplacian and on generalized Newtonian models such as Carreau and power-law (also known as Ostwald--de Waele fluids)
\cite{glowinski1975approximation,wei1989finite,barrett1993finite,barrett1993finite_2,barrett1994quasi}.
FEM optimality was shown in \cite{ebmeyer2005quasi} in terms of a quasi-norm; adaptive finite element methods for the $p$-Laplacian were analyzed in
\cite{DieningKreuzer2008,belenki2012optimality,carstensen2003posteriori}, while the finite element approximation of $p$-Stokes systems was studied in \cite{belenki2012finite}. High-order FEM approximation properties for the $p$-Laplacian have been investigated in \cite{ainsworth}.
Discretizations of energies with non-standard growth have also been investigated, for instance in \cite{breitDiening}.
Structurally related developments include discontinuous Galerkin methods for quasilinear elliptic problems
\cite{Houston1,Houston2}, and, more recently, Hybrid High-Order methods
\cite{botti2021hybrid,di2017hybrid,carstensen2021unstabilized}
and Virtual Element methods
\cite{antonietti2024virtual,antonietti2026virtual}.

The nonlinear algebraic systems arising from these discretizations are commonly solved either by Newton-type methods or by fixed-point linearizations. Among the latter, the \kac method is one of the most widely used approaches. It originated in nonlinear magnetostatics \cite{kavcur1968convergence} and was later developed as a general linearization technique for nonlinear PDEs and variational inequalities \cite{fuvcik1974kavcanov,han1997kavcanov, nevcas1983solution}, see also  
\cite{zeidler2013nonlinear}.
For energies of the type \eqref{eq:energy_general}, a relaxed version of the \kac method has been studied for the $p$-Laplacian \cite{diening2020relaxed,balci2023relaxed}, shear-thinning Carreau-type fluids \cite{heid2022convergence}, and regularized Bingham models \cite{heid2022adaptive}. For a damped version of the algorithm, see \cite{heid2022modified}, and for recent work on the relaxed version, see \cite{diening2025guaranteed}. Finally, closely related to the \kac method is the adaptive iterative linearized Galerkin method \cite{heid2020adaptive}.

The starting point of the present work is the differential-algebraic  system introduced in \cite{haskovec2026} for modeling the evolution of biological transportation networks
\[
\dd{s}{t} = |\grad u|^2 - s^{\gamma-1},
\quad
-\diver\bigl((s+r)\grad u\bigr) = f,\quad \gamma > 1,
\]
where $r>0$ is an arbitrarily small regularization parameter ensuring the well-posedness of the linear elliptic problem. The steady states of this system solve the Euler--Lagrange equations of the regularized $p$-Laplacian
\[
\int_\Omega (|\grad u|^{p-2} + r)\grad u \cdot \grad v \dx
=
\int_\Omega f v \dx,
\quad
p = \frac{2\gamma}{\gamma-1}.
\]
Rewriting $\gamma$ in terms of $p$ and introducing the change of variables
$c=s^{2/(p-2)}$, we obtain
\[
\dd{c}{t}
=
\frac{2}{p-2} c^{\frac{4-p}{2}}\bigl(|\grad u|^2-c\bigr),
\quad
-\diver \left( \left(c^{\frac{p-2}{2}}+r\right) \grad u\right) = f .
\]
If the goal is to recover only the steady state, the prefactor in the ordinary differential equation (ODE) is not essential, since it changes the transient dynamics but not the equilibria. One is thus tempted to deliberately drop it and use the simpler evolution
\begin{equation}\label{eq:ode_intro}
\dd{c}{t} = |\grad u[c]|^2 - c,
\end{equation}
where $u[c]$ denotes the solution of the preceding weighted linear elliptic problem for the given $c$.
This seemingly minor simplification was the starting point of the present paper: it produces an auxiliary evolution that is no longer tied to the specific algebraic structure of the power-law model (which is only encoded through the elliptic constraint) and can therefore be extended to more general constitutive laws. Note that, for $r=0$, a forward Euler discretization of \eqref{eq:ode_intro} with time step equal to one gives the \kac iteration for the $p$-Laplacian \cite[Equation~1.4]{diening2020relaxed},
\[
c_{n+1} = |\grad u_n|^2, \quad \int_\Omega c_{n+1}^{(p-2)/2}\,\grad u_{n+1}\cdot \grad v \dx
=
\int_\Omega f v \dx .
\]

Bartels, Diening, and Nochetto \cite{bartels2018unconditional} studied gradient flows that evolve $u$ and proved unconditional energy stability for a class of Orlicz energies with nonincreasing viscosity, together with fully discrete error estimates for singular power-law flows; we instead evolve $c$. The idea of isolating the nonlinearity in an auxiliary variable is also behind the recent work \cite{dediego2026transform} for the singular regime $1<p<2$, where the flux is introduced as an independent unknown and condensed out of each linearization rather than evolved. Here the auxiliary variable is a scalar field representing, at equilibrium, the squared magnitude $|Du|^2$, and the aim is not to derive an equivalent minimization problem for \eqref{eq:energy_general}, but a variationally consistent framework for designing and analyzing approximation schemes based on a sequence of uniformly elliptic linear state problems of the type
\begin{equation}\label{eq:constraint_gen}
\int_{\Omega}\mu(c(x))\,D u\fro D v\dx
= F(v)
,\quad\mu(r):=\frac{\phi'(\sqrt{r})}{\sqrt{r}}.
\end{equation}
Note that, in the non-Newtonian fluids literature (see e.g. \cite[Section 4.1]{bird1986dynamics}), the viscosity is written as a function $\widehat\mu(\dot\gamma)$, where the shear rate $\dot\gamma$ is the magnitude of the rate-of-strain tensor, i.e. $\dot\gamma = \sqrt{2}|\eps(\bu)|$.
Here the auxiliary variable tracks, modulo the $\sqrt{2}$ scaling, the squared modulus $|\eps(\bu)|^2$.

We first construct the auxiliary energy on uniformly truncated functions, $\epsilon\le c(x)\le\epsilon^{-1}$, for fixed $\epsilon\in(0,1)$. The viscosity derivative defines a metric under which the auxiliary space is a complete geodesic space, isometric to $\Ltwo$. Under suitable assumptions, lower semicontinuity and geodesical $\lambda$-convexity yield exponential convergence of minimizing movements to a minimizer in the truncated set, with a rate independent of $\epsilon$ \cite{ambrosio2005gradient,degiorgi}. For power laws, the result covers $4/3\le p\le4$, excluding the linear case $p=2$. We do not pass to the limit $\epsilon\to0$ at the continuous level: the role of this theory is to identify the metric and the convexity structure that we then exploit to prove convergence of the finite element flow to the nonlinear Galerkin solution.

Instead of discretizing the method of minimizing movements directly, we derive a finite-dimensional evolution using the theory of metric gradient flows \cite{liero2013gradient} (see also \cite[Section 4]{mielke2023introduction}). The discrete metric induces an explicit Riesz map, whose inverse converts the first variation of the auxiliary energy into its metric gradient. This identification yields the evolution \eqref{eq:ode_intro} without an additional variational solve, with the rheology entering through the state constraint \eqref{eq:constraint_gen}. Under the assumptions of Theorem~\ref{thm:ode_bound}, we prove global existence and uniqueness of the solution on the positive cone of the finite element space for the auxiliary variable. Additional structural assumptions give convergence to the nonlinear Galerkin solution, in the same power-law range as the continuous theory. Outside this range, the framework suggests asymptotic rates near equilibrium consistent with the experiments.

We consider time discretizations that include the \kac iteration. Inspired by pseudo-transient continuation \cite{kelley1998convergence,coffey2003pseudotransient}, we propose a time-adaptive iteration compatible with scalable linear solvers. We do not analyze convergence of these fully discrete schemes. Experiments with power-law and Carreau--Yasuda models \cite{yasuda1979investigation,yasuda1981shear} support the predicted rates within the theoretical range and the asymptotic estimates outside it. We observe essentially mesh-independent iteration counts for the Pseudo method in the power-law and Bingham tests, whereas the counts increase with mesh refinement in our tests of the optimal-design models from \cite{bartels2008convergent,carstensen2021unstabilized}. The methods are competitive with Newton's method and succeed in several cases where it fails, including the $p$-Laplacian with $p=100$.

Section~\ref{sec:preliminaries} introduces the functional setting and metric. Section~\ref{sec:auxiliary_energy} develops the auxiliary energy and continuous flow. Section~\ref{sec:impl} treats finite element discretization and time stepping. Section~\ref{sec:experiments} reports numerical experiments, and Section~\ref{sec:conclusions} discusses future directions. Given the length of the analysis, all proofs are collected in Appendix~\ref{sec:appendix} to make the main text easier to follow.

\section{Preliminaries}\label{sec:preliminaries}

\subsection{Orlicz--Sobolev spaces}

\begin{definition}[N-function]
We say that $\phi:[0,+\infty)\to[0,+\infty)$ is an N-function if $\phi$ is convex, continuous,
$\phi(0)=0$, $\phi(r)>0$ for $r>0$, $\lim_{r \to 0} \frac{\phi(r)}{r}=0$, and
$\lim_{r \to +\infty}\frac{\phi(r)}{r}=+\infty$.
\end{definition}
The right derivative $\phi'_+$ is non-decreasing and right-continuous, satisfies $\phi'_+(0) = 0$ and $\phi'_+(r) > 0$ for $r>0$, and obeys $\lim_{r \to +\infty}\phi_+'(r)=+\infty$ \cite{musielak1983orlicz}. An equivalent characterization of an N-function is the integral representation
\begin{equation}\label{eq:n_func_def}
\phi (r) = \int^r_0 \phi'_+(s)\ds,
\end{equation}
with an integrand satisfying these properties. For simplicity, we assume throughout that $\phi$ has enough regularity for the following derivatives to be well defined:
\begin{assumption}\label{ass:phi_regular}
$\phi \in C^1([0,+\infty))\cap C^3((0,+\infty))$.
\end{assumption}
\noindent We also assume the following growth condition:
\begin{assumption}\label{ass:phi_assum_reflex}
The Simonenko indices  of $\phi$ \cite{simonenko1963interpolation} are bounded, i.e.
\[
\simphi:=\inf_{r> 0} \frac{r\phi'(r)}{\phi(r)} > 1, \quad \sipphi =   \sup_{r> 0} \frac{r\phi'(r)}{\phi(r)} < +\infty.
\]
\end{assumption}
\noindent This assumption ensures that $\phi$ and its complementary function %
$\phi^\ast(\eta):=
\sup_{r\ge 0}\big(\eta r-\phi(r)\big)$
satisfy a $\Delta_2$ condition \cite{phi_reflex}, i.e. 
there exists $C\ge 1$ such that $\phi(2r)\le C\phi(r)$ for all $r\ge 0$, and similarly for $\phi^*$.
Under Assumption \ref{ass:phi_assum_reflex}, the natural spaces to study energies of the form \eqref{eq:energy_general} are the Orlicz--Sobolev spaces \cite{diening2011lebesgue}. Associated with $\phi$ is the Orlicz space
\[
\Lphi:=\left\{v:\Omega\to\R\text{ measurable}:\int_\Omega\phi(|v|/\lambda)\dx<+\infty\text{ for some }\lambda>0\right\},
\]
which is a separable and reflexive Banach space equipped with the Luxemburg norm
\[
\|v\|_{\Lphi}
:=\inf\Big\{\lambda>0:\int_{\Omega}\phi(|v|/\lambda)\dx\le 1\Big\},
\]
and, under the sole $\Delta_2$ condition, we have the equivalence (see e.g. \cite{donaldson1971orlicz})
$$v\in \Lphi \iff \int_\Omega \phi(|v|) \dx < +\infty.$$
The Orlicz--Sobolev space $W^{1,\phi}(\Omega)$ consists of all $v\in \Lphi$
such that $\grad v\in L^\phi(\Omega;\R^d)$,
and the following continuous embeddings hold for bounded domains \cite[Lemmas 3.7.7 and 6.16]{harjulehto2019generalized}
\begin{align}
\phi(r)/r^2 \text{ non-decreasing}&\Rightarrow \Lphi  \hookrightarrow \Ltwo,\quad W^{1,\phi}(\Omega) \hookrightarrow H^1(\Omega),\label{eq:inclusions_convex}\\
\phi(r)/r^2 \text{ non-increasing}&\Rightarrow \Ltwo  \hookrightarrow \Lphi,\quad  H^1(\Omega)  \hookrightarrow W^{1,\phi}(\Omega).\label{eq:inclusions_concave}
\end{align}

Instead of working with the energy functional \eqref{eq:energy_general}, we consider a reparameterization in terms of $|D u|^2$, and introduce the function
\begin{equation}\label{eq:Phi}
\Phi(r):=\phi(\sqrt{r}),
\end{equation}
so that $\phi(|D u|)=\Phi(|D u|^2)$; from \eqref{eq:n_func_def} we have the following characterization 
\begin{equation}\label{eq:Phi_mu}
\Phi(r) := \frac{1}{2}\int^r_0\mu(s)\ds,\quad
\mu(r):=2\Phi'(r)=\frac{\phi'(\sqrt{r})}{\sqrt{r}}.
\end{equation}
The function $\Phi$ is continuous and strictly increasing, with $\Phi(0)=0$ and $\Phi(r)>0$ for $r>0$. Moreover, $\mu(r)>0$ for $r>0$. However, since the mapping $r\mapsto \sqrt{r}$ is concave, $\Phi$ is generally not an N-function. A direct calculation yields
\begin{equation}\label{eq:Phipp}
\Phi''(r) = \frac{1}{4 r\sqrt{r}}\left(\phi''(\sqrt{r})\sqrt{r} - \phi'(\sqrt{r})\right) = \frac{1}{4 t^3}\left(\phi''(t)t - \phi'(t)\right) =\frac{1}{4 t} \dd{}{t}\left(\frac{\phi'(t)}{t}\right), \quad t=\sqrt{r}.
\end{equation}
This leads to the next assumption, which we require to hold throughout this work.
\begin{assumption}\label{ass:Phi_concace_or_convex_strictly}
$\Phi$ is either strictly convex (i.e. $\mu'(r) >0$) or strictly concave ($\mu'(r) <0$) on $(0,+\infty)$.
\end{assumption}
Under this assumption, $\phi'$ is a strictly increasing function on $(0,+\infty)$ and the complementary function $\phi^\ast$ admits the representation 
\begin{equation}\label{eq:phi_complement}
\phi^\ast(r)=\int^r_0 (\phi')^{-1}(s)\ds.
\end{equation}
From \eqref{eq:Phipp} we thus have that $\Phi$ is strictly convex (resp. concave) if and only if $\phi'(r)/r$ is increasing (resp. decreasing),
and we derive the convenient representation
\begin{equation}\label{eq:mup_exact_final}
\mu'(r)=\frac{1}{2r}\left(\frac{t\phi''(t)}{\phi'(t)}-1\right)\mu(r), \quad t = \sqrt{r}.
\end{equation}
The monotone form of de L'H\^opital's rule gives the following lemma, which ensures that one of the embeddings \eqref{eq:inclusions_convex}--\eqref{eq:inclusions_concave} holds. %
\begin{lemma}\label{lem:from_Phi_to_inclusions}
If $\frac{\phi'(r)}{r}$ is non-decreasing (resp. non-increasing) then $\phi(r)/r^2$ is non-decreasing (resp. non-increasing). Therefore
\[
\begin{aligned}
\Phi \text{ convex on }(0,+\infty) &\Rightarrow \Lphi  \hookrightarrow \Ltwo,\quad W^{1,\phi}(\Omega) \hookrightarrow H^1(\Omega),\\
\Phi \text{ concave on }(0,+\infty) &\Rightarrow \Ltwo  \hookrightarrow \Lphi,\quad  H^1(\Omega)  \hookrightarrow W^{1,\phi}(\Omega).
\end{aligned}
\]
\end{lemma}

For our analysis, we will need two further assumptions. The first is the so-called \emph{uniform convexity} assumption introduced in \cite{diening2008fractional}.
\begin{assumption}\label{ass:sindices}
Let $\smphi$ and $\spphi$ be defined as
\[
\smphi := \inf_{r> 0} \frac{r\phi''(r)}{\phi'(r)} ,\quad  \spphi :=\sup_{r> 0} \frac{r\phi''(r)}{\phi'(r)},
\]
and assume that
\[
0 < \smphi \le \spphi < +\infty.
\]
\end{assumption}
\Needspace{6\baselineskip}
The next lemma gives two-point growth bounds for $\phi'$ and $\mu$ (see Appendix Lemma~\ref{lem:phiprime_bounds_proof}).

\begin{lemma}\label{lem:phiprime_bounds}
Let $0<s\le t$ and $\smphi, \spphi$ be given in Assumption \ref{ass:sindices}. Then
\begin{enumerate}
\item[(i)] $\phi'$ satisfies the two-point growth:
$\displaystyle \left(\frac{s}{t}\right)^{\spphi}\le \frac{\phi'(s)}{\phi'(t)}\le \left(\frac{s}{t}\right)^{\smphi},$
\item[(ii)] $\mu$ satisfies the two-point growth:
$\displaystyle \left(\frac{s}{t}\right)^{\frac{\spphi-1}{2}}\le \frac{\mu(s)}{\mu(t)}\le \left(\frac{s}{t}\right)^{\frac{\smphi-1}{2}}.$
\end{enumerate}
\end{lemma}
Lemma~\ref{lem:phiprime_bounds} also yields the following bounds (see Appendix Lemma~\ref{lem:integ_E_c_proof}).
\begin{lemma}\label{lem:integ_E_c}
Let $\smphi, \spphi$ be given as in Assumption \ref{ass:sindices}. Then, for all $r\ge0$,
\[
\frac{1}{\spphi+1}\,r\phi'(r) \le \phi(r) \le \frac{1}{\smphi+1}\,r\phi'(r),\quad \frac{1}{\spphi+1}\,r\mu(r) \le \Phi(r) \le \frac{1}{\smphi+1}\,r\mu(r),
\]
where $r\phi'(r)$ and $r\mu(r)$ at $r=0$ are understood through their continuous extensions $\lim_{s\to 0^+} s\phi'(s)=0$ and $\lim_{s\to 0^+} s\mu(s)=0$.
\end{lemma}
The preceding bounds imply the chain of inequalities
\[
1 < 1 + \smphi \le \simphi \le \sipphi \le 1 + \spphi < +\infty.
\]

\subsection{Metric space}
Associated with $\Phi$ is the function space $\LPhi$, defined, with some abuse of notation, as
\begin{equation}\label{eq:LPhi_def}
\LPhi := \left\{ c : \Omega \rightarrow \R \text{ measurable} : \sqrt{|c|} \in \Lphi\right\}.
\end{equation}
Note that $\LPhi$ is a reflexive and separable Orlicz space only in the convex $\Phi$ case.
When $\Phi$ is concave, the Luxemburg norm is not well defined, and $\LPhi$ is not an Orlicz space. Nevertheless, %
the modular is always well defined, and we have the equivalences
\begin{equation}\label{eq:c_modular}
c\in \LPhi
\Longleftrightarrow
\int_\Omega \phi \left(\sqrt{|c|}\right)\dx<+\infty
\Longleftrightarrow
\int_\Omega \Phi (|c|)\dx<+\infty.
\end{equation}
To prove that $\LPhi$ is a function space also in the concave case, we only need to show it is closed under addition, since the scaling by constants clearly yields functions in $\LPhi$. Since $\Phi$ is non-decreasing, $\Phi(0)=0$, and it is subadditive on $[0,+\infty)$, we have
$\Phi(|r_1+r_2|)\le \Phi(|r_1|+ |r_2|)\le \Phi(|r_1|)+\Phi(|r_2|)$, %
for all $r_1,r_2$,
which, when integrated, gives $c_1+c_2\in \LPhi$ for all $c_1,c_2\in\LPhi$. 

In order to define a topology on $\LPhi$ which is valid in both the convex and concave cases, we introduce a further assumption.
\begin{assumption}\label{ass:sindices_strict}
Let $\smphi$ and $\spphi$, defined in Assumption \ref{ass:sindices}, be such that one of the following two conditions holds
\[
0 < \smphi \le \spphi < 1,
\quad \text{or} \quad
1 < \smphi \le \spphi < +\infty.
\]
\end{assumption}
We then introduce the odd function
\begin{equation}\label{eq:psi_def}
\psi(r) := \operatorname{sgn}(r)\int^{|r|}_0 \sqrt{\frac{\omega}{2}\mu'(s)}\ds,
\end{equation}
where
\begin{equation}\label{eq:omega_def}
\omega := \operatorname{sgn}(\mu'(r)),\quad r>0,
\end{equation}
and note that $\psi$ is a continuous, strictly increasing function under Assumption \ref{ass:Phi_concace_or_convex_strictly}. The following lemma proves that under Assumption \ref{ass:sindices_strict} it is a homeomorphism  and thus its inverse $\psi^{-1}(r)$ is also a continuous, strictly increasing function; this allows us to identify functions of $\LPhi$ with functions of $\Ltwo$ via the bijective map
\begin{equation}\label{eq:psi_map_def}
\Psi : \LPhi \rightarrow \Ltwo, \quad \Psi(c)(x) := \psi(c(x)),
\end{equation}
and obtain the characterization $c \in \LPhi \iff \Psi(c) \in \Ltwo$. %
For the proof, see Appendix Lemma~\ref{lem:psi_bounds_mu_proof}.
\begin{lemma}\label{lem:psi_bounds_mu}
Let $\psi$ be defined as in \eqref{eq:psi_def} and let Assumption \ref{ass:sindices_strict} hold. Then:
\begin{enumerate}
\item[(i)] There exist constants $0<C_1< C_2<+\infty$ depending only on $\smphi, \spphi$ such that
\[
C_1\mu(r)r \le \psi(r)^2 \le C_2\mu(r)r,\quad \forall r \ge 0.
\]
\item[(ii)] There exist constants $0<C_3< C_4<+\infty$ depending only on $\smphi, \spphi$ such that
\[
C_3\Phi(r) \le \psi(r)^2 \le C_4\Phi(r),\quad \forall r \ge 0.
\]
\end{enumerate}
\end{lemma}

\begin{remark}\label{rem:ass_sindices_strict_not_necessary}
Assumption \ref{ass:sindices_strict} ensures that $\Phi$ satisfies Assumption \ref{ass:Phi_concace_or_convex_strictly}; in fact, %
\[
 1 < \smphi \Rightarrow \Phi \text{ strictly convex}, \quad \spphi < 1 \Rightarrow \Phi \text{ strictly concave},
\]
but the converse implication is not true in general.
The equivalence $c \in \LPhi \iff \Psi(c) \in \Ltwo$ can also be obtained by relaxing the lower bound condition in statement (ii) of Lemma \ref{lem:psi_bounds_mu} as
\[
\Phi(r)\le C_2(1+\psi^2(r)).
\]
\end{remark}

Using the bijective map $\Psi$ given in \eqref{eq:psi_map_def}
we can define a metric on $\LPhi$ as
\begin{align}\label{eq:X_distance}
d(c_0,c_1) := \|\Psi(c_1) - \Psi(c_0)\|_{\Ltwo}.
\end{align}
Hence, $d(\cdot,\cdot)$ is the usual distance in $\Ltwo$ under the change of variables $\eta = \Psi(c)$. Clearly, non-negativity, symmetry, and the triangle inequality hold; since $\Psi$ is a bijection, we can show that $d(u,v)=0$ implies $u=v$, and thus $d(\cdot,\cdot)$ is a metric.
A key tool in our continuous analysis is the existence of a geodesic, denoted by $\gamma_g(s)$, connecting two arbitrary points $c_0$ and $c_1$ of $\LPhi$, i.e.
\[
\gamma_g : [0,1] \rightarrow\LPhi,\quad \gamma_g(0) = c_0, \quad \gamma_g(1) = c_1,
\]
and such that
\begin{equation*}%
d(\gamma_g(s),\gamma_g(r)) = |s-r|\,d(c_0,c_1)\quad \forall s,r\in[0,1].
\end{equation*}
Since $\Psi$ is a bijection, this geodesic is obtained by pulling back the straight segment between 
$\Psi(c_0)$ and $\Psi(c_1)$ in $\Ltwo$
\begin{equation}\label{eq:geodesics}
\gamma_g(s) := \Psi^{-1}((1-s)\Psi(c_0) + s\,\Psi(c_1)), \quad 0\le s \le 1.
\end{equation}
The space $\LPhi$ endowed with $d$ is complete, and $\gamma_g$ in \eqref{eq:geodesics} is a geodesic (see Appendix Lemma~\ref{lem:C_d_complete_geodesic_proof}).

\begin{remark}
For the power-law case, $\phi(r)=\frac{1}{p}r^p$, $\mu(r) = r^{(p-2)/2}$, with $p>1$, and $\smphi=\spphi=p-1$. We have $\LPhi = L^{p/2}(\Omega)$, and $\psi(r) = C \text{sgn}(r) |r|^{p/4}$, with $C=2\sqrt{|p-2|}/p$ for $p\ne 2$. %
When $c_0,c_1\in L^{p/2}(\Omega)$, we have $|c_0|^{p/4},|c_1|^{p/4}\in \Ltwo$, and thus the metric is well defined. For non-negative $c_0$ and $c_1$, the geodesic is
\[
\gamma_g(s) = \left((1-s)c_0^{p/4} + s\,c_1^{p/4}\right)^{4/p}.
\]
\end{remark}

\section{Auxiliary energy reformulation}\label{sec:auxiliary_energy}

\subsection{The scalar case}

We consider the minimization problem 
\begin{equation}\label{eq:scalar_primal}
\argmin_{u\in W^{1,\phi}_0(\Omega)} \int_{\Omega}\phi(|\grad u|)\dx - F(u),
\end{equation}
where 
$W^{1,\phi}_0(\Omega)$ is the closure of $C^\infty_c(\Omega)$ with respect to the
$W^{1,\phi}(\Omega)$-norm and $F$ is a continuous linear functional in $W^{1,\phi}_0(\Omega)^*$.
The Euler--Lagrange equations associated with \eqref{eq:scalar_primal} read:
Find $u\in W^{1,\phi}_0(\Omega)$ such that
\begin{equation}\label{eq:scalar_EL}
\int_{\Omega}
\frac{\phi'(|\grad u|)}{|\grad u|}\grad u
\cdot \grad v \dx
= F(v)
\quad \forall v\in W^{1,\phi}_0(\Omega).
\end{equation}
For regularity results of the minimizers of \eqref{eq:scalar_primal} under Assumption \ref{ass:sindices}, see \cite{cianchi_grad_bound,cianchi2018second} and the references therein.
Using the definition of $\Phi$ from \eqref{eq:Phi}, we rewrite the energy in the equivalent form
\begin{equation}\label{eq:scalar_primal_Phi} 
E(u):=\int_{\Omega}\Phi(|\grad u|^2)\dx - F(u),
\end{equation}
yielding the Euler--Lagrange equations: Find $u\in W^{1,\phi}_0(\Omega)$ such that
\begin{equation}\label{eq:scalar_EL_mu}
\int_{\Omega}\mu(|\grad u|^2)\,\grad u\cdot\grad v\dx
= F(v)
\quad \forall v\in  W^{1,\phi}_0(\Omega).
\end{equation}

\begin{remark}\label{rem:mu_bounded}
Assumption \ref{ass:sindices_strict} used in the metric-space construction excludes the cases in which the viscosity law 
becomes asymptotically Newtonian as its argument approaches $0$ or $+\infty$. In such cases, since
\[
\frac{t\phi''(t)}{\phi'(t)}
=
1+\frac{2r\mu'(r)}{\mu(r)}, \quad t = \sqrt{r}
\]
we have
\[
\frac{2r\mu'(r)}{\mu(r)}\to 0
\Rightarrow
\frac{t\phi''(t)}{\phi'(t)}\to 1.
\]
Therefore, in the convex case, one obtains $\smphi=1$, whereas in the concave case $\spphi=1$. In this sense, the asymptotically Newtonian cases excluded by Assumption \ref{ass:sindices_strict} are often easier to study, since in the asymptotic regime the coefficient $\mu(r)$ approaches
a positive constant and the operator in \eqref{eq:scalar_EL_mu} behaves like a uniformly elliptic linear operator with bounded coefficients.
\end{remark}

To obtain weighted linear state problems, we introduce the extended functional
\begin{equation*} %
J:\Wphiz\times \LPhi \rightarrow [-\infty,+\infty], \; J(u,c)
:=
\begin{cases}
\displaystyle \int_{\Omega} \Phi(c)+\frac{1}{2}\,\mu(c)\big(|\grad u|^2-c\big)\dx- F(u), & c\in\LPhip,\\
+\infty,&\text{otherwise},
\end{cases}
\end{equation*}
where
\begin{equation}\label{eq:L_Phi}
\LPhip := \left\{ c \in \LPhi : c(x) \ge 0 \,\text{ a.e. in }\,\Omega \right\}.
\end{equation}
For a fixed $u\in \Wphiz$, since $2\Phi'(c) = \mu(c)$, stationarity of $J$ with respect to $c$ implies
\begin{equation*}%
\mu'(c)(|\grad u|^2-c) = 0,\quad\text{ a.e. in }\Omega,
\end{equation*}
and thus, from Assumption \ref{ass:Phi_concace_or_convex_strictly}, we obtain
\begin{equation*}%
|\grad u|^2 = c\, \text{ a.e. in }\{x\in\Omega \st c(x) > 0\}.
\end{equation*}
Stationarity of $J$ with respect to $u$ instead yields
\begin{equation}\label{eq:scalar_state_nwp}
a[c](u,v) = F(v),
\quad \forall v\in \Wphiz,
\end{equation}
where
\begin{equation}\label{eq:scalar_bilin}
a[c](u,v):=\int_{\Omega}\mu(c)\,\grad u\cdot\grad v\dx.
\end{equation}
For a general $c\in \LPhip$, the variational problem \eqref{eq:scalar_state_nwp} is not well-posed on $\Wphiz$ and results like the Minty--Browder theorem \cite[Theorem 9.14-1]{ciarlet} cannot be applied. In the concave case we can in fact prove the coercivity of the functional in $\Wphiz$, but the bilinear form \eqref{eq:scalar_bilin} is not well defined since $\mu(c)\grad\cdot$ does not map $\Wphiz$ into $\Wphiz^*$. On the other hand, when $\Phi$ is convex, the linear operator induced by the bilinear form is hemicontinuous, but we lose coercivity instead.

When $\Phi$ is concave,  we have 
\[
\Phi(r) \leq \Phi(c) + \frac{1}{2} \mu(c) (r - c)\quad \forall\, r,c \ge 0.
\]
Integrating  on $\Omega$ and subtracting $F(u)$ from both sides yields
\[
E(u) \le J(u,c), \quad \forall u\in \Wphiz,\, \forall c \in \LPhip,
\]
and thus
\begin{equation}\label{eq:Eu_bound_conc}
\inf_{u\in \Wphiz}E(u) \le \inf_{u\in \Wphiz} J(u,c) =  \Jc(c) + \inf_{u\in \Wphiz}  \Big( \int_\Omega  \frac{1}{2}\mu(c)|\grad u|^2 \dx  - F(u)\Big),
\end{equation}
where the functional
\begin{equation}\label{eq:Jc_def}
\Jc(c) := \int_\Omega \Phi(c)  -\frac{1}{2} \mu(c) c \dx
\end{equation}
is well defined and finite for all $c\in \LPhip$ by Lemma \ref{lem:integ_E_c}.
In the case of convex $\Phi$, the comparison with the tangent functional reverses
\[
\Phi(r) \geq \Phi(c) + \frac{1}{2} \mu(c) (r - c)\quad \forall\, r,c \ge 0,
\]
yielding
\[
E(u) \ge J(u,c), \quad \forall u\in \Wphiz,\, \forall c \in \LPhip,
\]
and
\begin{equation}\label{eq:Eu_bound_conv}
\inf_{u\in \Wphiz}E(u) \ge \inf_{u\in \Wphiz} J(u,c) =  \Jc(c) + \inf_{u\in \Wphiz}  \Big( \int_\Omega  \frac{1}{2}\mu(c)|\grad u|^2 \dx  - F(u)\Big).
\end{equation}

We then define an auxiliary energy in $c$ by replacing the space $\Wphiz$ in the minimization of the Dirichlet integral with $\Honez$ and constraining $c$ to the set  
\begin{equation}\label{eq:C_d_def}
\Ce := \left\{c \in \LPhi \st  \epsilon \le c(x) \le 1/\epsilon \,\text{ a.e. in } \Omega \right\},
\end{equation}
for an arbitrary $\epsilon \in (0,1)$, which ensures the following uniform bounds almost everywhere
\begin{equation}\label{eq:mu_unif_bounds} 
0 < \mu_-  \le \mu(c(x)) \le \mu_+ < +\infty, \quad \mu_- :=\min\{\mu(\epsilon), \mu(1/\epsilon)\}, \quad \mu_+ :=\max\{\mu(\epsilon), \mu(1/\epsilon)\}.
\end{equation}
For the problem to be well-posed, we need a further assumption.
\begin{assumption}\label{ass:F_Hm1}
The linear functional $F$ admits a continuous extension to $\Honez^*$.
\end{assumption}
Note that the above condition is always satisfied in the concave case but not necessarily in the convex case. We also point out that, for ease of discussion, from now on, all the inequalities (and equalities) involving pointwise bounds for functions must be understood in the almost everywhere sense.
Under Assumption \ref{ass:F_Hm1} we can consider the well-posed functional for a given $c\in\Ce$, see, e.g., \cite[Section 25.2]{ernguermondII},
\begin{equation}\label{eq:Ju_def}
\Ju(c) := \inf_{u\in \Honez} \left( \int_{\Omega}\frac{1}{2}\mu(c)|\grad u|^2\dx - F(u)\right),
\end{equation}
whose Euler--Lagrange equations lead to the weighted Poisson problem: Find $u[c] \in \Honez$ such that
\begin{equation}\label{eq:scalar_state}
\int_\Omega\mu(c)\grad u[c]\cdot\grad v \dx = F(v),
\quad \forall v\in \Honez,
\end{equation}
yielding the identities
\begin{equation}\label{eq:energy_u_c}
\int_{\Omega}\mu(c)|\grad u[c]|^2\dx = F(u[c]),\quad \Ju(c) = -\frac{1}{2}F(u[c]) = -\frac{1}{2}\int_{\Omega}\mu(c)|\grad u[c]|^2\dx.
\end{equation}
Combining \eqref{eq:Eu_bound_conc} with the embedding \eqref{eq:inclusions_concave} in the concave case, and \eqref{eq:Eu_bound_conv} with \eqref{eq:inclusions_convex} in the convex case, we find, for all $c\in\Ce$,
\[
\begin{aligned}
\inf_{u\in \Wphiz}E(u) &\le \Jc(c) + \Ju(c), &&\Phi\text{ concave},\\
\inf_{u\in \Wphiz}E(u) &\ge \Jc(c) + \Ju(c), &&\Phi\text{ convex}.
\end{aligned}
\]
Using the definition of $\omega$ given in \eqref{eq:omega_def} we arrive at the extended auxiliary energy functional
\begin{equation}\label{eq:E_c}
\Ee : \LPhi \to \Rinf, \quad \Ee(c) :=
\begin{cases}
-\omega \Jc(c) -\omega \Ju(c),& c\in\Ce,\\
+\infty,&\text{otherwise},
\end{cases}
\end{equation}
which is a proper, well-defined functional bounded from below by the value of the energy of the minimizer of \eqref{eq:scalar_primal_Phi} multiplied by $-\omega$, since we have
\begin{equation}\label{eq:Ee_bounded}
-\omega E(u[c]) \le -\omega J(u[c], c) =\Ee(c),\quad  c \in \Ce. \end{equation}

\begin{remark}
For power laws, the substitution $c=a^2$ identifies $J(u,c)$
with the extended energy of Diening--Fornasier--Tomasi--Wank
\cite[Section~2]{diening2020relaxed}.
For concave $\Phi$, the tangent inequality gives $E(u)\le J(u,c)$;
the same upper bound is used by Heid--S\"uli
\cite[Lemma~3.1]{heid2022convergence} and
Bartels--Diening--Nochetto
\cite[Lemma~3.3]{bartels2018unconditional}.
\end{remark}

\subsection{The incompressible vector-valued case}

The extension of the previous results to the vector-valued case
\[
\argmin_{\bu\in W^{1,\phi}_{0}(\Omega;\R^d)}\int_{\Omega}\phi(|\grad \bu|)\dx - F(\bu),
\]
is straightforward. In this section, we instead consider the  minimization problem
\begin{equation}\label{eq:vector_min}
\argmin_{\bu\in W^{1,\phi}_{0,\mathrm{div}}(\Omega;\R^d)}\int_{\Omega}\phi(|\eps(\bu)|)\dx - F(\bu),
\end{equation}
where
\[
W^{1,\phi}_{0,\mathrm{div}}(\Omega;\R^d)
:=\{\bv\in \vWphiz: \diver \bv=0\},
\quad
\eps(\bu):=\frac12(\grad \bu+\grad \bu^{\mathsf T}),
\]
and $F\in W^{1,\phi}_0(\Omega;\mathbb{R}^d)^\ast$.
For existence, uniqueness, and regularity of the solutions under Assumption \ref{ass:sindices} see the recent work \cite{giannetti2024fractional} and the references therein.
The auxiliary functional is constructed by following verbatim the scalar case, replacing 
$|\grad u|^2$ by $|\eps(\bu)|^2$, and
using
\begin{equation}\label{eq:Ju_def_vector}
\Jbu(c) := \inf_{\bu\in H^{1}_{0,\mathrm{div}}(\Omega;\R^d)}  \left(\int_\Omega  \frac{1}{2}\mu(c)|\eps(\bu)|^2 \dx -F(\bu)\right),
\end{equation}
for a given $c\in\Ce$ under the following assumption:
\begin{assumption}\label{ass:F_Hm1_vector}
The linear functional $F$ admits a continuous extension to $\vHonez^\ast$.
\end{assumption}
We enforce incompressibility using a pressure Lagrange multiplier in
\[
Q:=\left\{p\in \Ltwo : \int_\Omega p \dx = 0\right\},
\]
and obtain the Euler--Lagrange equations for \eqref{eq:Ju_def_vector} given by the weighted incompressible Stokes system: Find $\bu[c] \in \vHonez$ and $p[c] \in Q$ such that
\begin{align}
\ba[c](\bu[c], \bv)- b(\bv, p[c]) &= F(\bv), \quad \forall \bv\in \vHonez, \label{eq:vector_state}
\\
b(\bu[c],q) &= 0,  \quad \forall q\in Q, \label{eq:vector_incomp}
\end{align}
where
\begin{equation}\label{eq:vector_bilin}
\ba[c](\bu, \bv) := \int_{\Omega}\mu(c)\,\eps(\bu):\eps(\bv)\dx, \quad 
b(\bu,q) := \int_{\Omega} q\,\diver \bu\dx.
\end{equation}
The well-posedness of the weighted incompressible Stokes problem is standard for a given $c\in\Ce$, see, e.g., \cite[Section 53.2]{ernguermondII}. Testing \eqref{eq:vector_state} with $\bv=\bu[c]$ and using the incompressibility of $\bu[c]$ gives %
\begin{equation}\label{eq:energy_u_c_vector}
\int_{\Omega}\mu(c)|\eps(\bu[c])|^2\dx = F(\bu[c]),\quad \Jbu(c) = -\frac{1}{2}F(\bu[c])  = -\frac{1}{2}\int_{\Omega}\mu(c)|\eps(\bu[c])|^2\dx,
\end{equation}
and the auxiliary functional has the same structure as the scalar case
\begin{equation}\label{eq:E_c_vector}
\Ee : \LPhi \to \Rinf, \quad \Ee(c) :=
\begin{cases}
-\omega \Jc(c) -\omega \Jbu(c),& c\in\Ce,\\
+\infty,&\text{otherwise}.
\end{cases}
\end{equation}
In this case too, $\Ee$ is bounded from below by the minimum of \eqref{eq:vector_min} multiplied by $-\omega$ and we have $-\omega E(\bu[c]) \le  \Ee(c)$ for all $c \in \Ce$. Stationarity with respect to a given $c$ yields
\[
\mu'(c)(|\eps(\bu)|^2-c)=0 \quad \text{a.e. in }\Omega.
\]

\subsection{Gradient flow in metric space}

In this section, we study the minimization of the auxiliary energy $\Ee$ using the theory of gradient flows in metric spaces \cite{ambrosio2005gradient}. Using the notation $C:=\LPhi$, we consider the metric functional system as in \cite{daneri2010lecture}
\begin{equation}\label{eq:metsym}
\begin{aligned}
&(C,d,\Ee), \text{ where } (C,d) \text{ is a complete metric space},\\
&\Ee : C\rightarrow \Rinf \text{ is a proper, lower semicontinuous functional, bounded from below},
\end{aligned}
\end{equation}
and provide sufficient conditions for the existence of flows satisfying evolution variational inequalities (EVI).
Under Assumption~\ref{ass:sindices_strict}, $\Jc$ is continuous on $\LPhip$ in the $d$-topology (Appendix Theorem~\ref{thm:Jc_cont_proof}). With the additional Assumption~\ref{ass:F_Hm1} in the scalar case or Assumption~\ref{ass:F_Hm1_vector} in the incompressible vector-valued case, $\Ju$ and $\Jbu$, respectively, are continuous on $\Ce$ (Appendix Theorem~\ref{thm:Ju_cont_proof}). These continuity results and the closedness of $\Ce$ in the $d$-topology imply the lower semicontinuity of $\Ee$; see Appendix Theorem~\ref{thm:Ee_lsc_proof}.

We first recall the definition of $\mathrm{EVI}_\lambda$-gradient flows.
\begin{definition}[$\mathrm{EVI}_\lambda$-gradient flows]
Let $\lambda\in\R$ and let $(C,d,\Ee)$ be a metric functional system as in \eqref{eq:metsym}. An $\mathrm{EVI}_\lambda$-gradient flow of $\Ee$ in $\dom{\Ee}$ is a family of continuous maps $S_t : \dom{\Ee} \rightarrow \dom{\Ee}$ with $t\ge 0$, such that for every $c_0\in \dom{\Ee}$
\[
S_{t+s}(c_0) = S_s(S_t(c_0)),%
\quad
\lim_{t\downarrow 0}S_{t}(c_0) = S_0(c_0) = c_0,
\quad\forall t,s\ge0,
\]
and the curve $t\mapsto S_t(c_0)$ satisfies the $\mathrm{EVI}_\lambda(C,d,\Ee)$ inequality
\begin{equation}\label{eq:evil}
\frac{1}{2}\dd{{}^+}{t} d(S_t(c_0),c)^2 + \frac{\lambda}{2} \,d(S_t(c_0),c)^2 + \Ee(S_t(c_0)) \leq \Ee(c),\quad \forall t\in(0,+\infty), \,\forall c\in C,
\end{equation}
where $\dd{{}^+}{t}$ is the upper right Dini derivative, defined by
\[
\dd{{}^+}{t}\eta(t) := \limsup_{h\downarrow 0}\frac{\eta(t+h)-\eta(t)}{h}.
\]
\end{definition}
We recall here some of the properties of $\mathrm{EVI}_\lambda$-gradient flows taken from \cite[Theorem 3.5]{muratori2020gradient} when the space $(C,d)$ is complete
\begin{equation}\label{eq:evil_props}
\begin{aligned}
\textbf{contraction property}:\quad& d(S_t(c_0),S_t(c_1))\leq e^{-\lambda(t-s)}d(S_s(c_0),S_s(c_1)),\quad 0\leq s\leq t < +\infty,\\
\textbf{asymptotic behaviour}:\quad& \lambda > 0 \Rightarrow d(S_t(c),c^*) \leq e^{-\lambda(t-t_0)}d(S_{t_0}(c),c^*),\quad c^\ast = \argmin_{c\in C}\Ee(c).
\end{aligned}
\end{equation}
The contraction estimate gives uniqueness
and continuous dependence on the initial datum. If $\lambda\ge 0$, the distance between two solutions
is nonincreasing, while if $\lambda>0$, it decays exponentially. In the latter case, the asymptotic estimate shows that the flow selects the unique minimizer $c^*$ of $\Ee$. %
Here $c^*$ minimizes the auxiliary energy on the fixed truncated set $\Ce$, and need not satisfy $c^*=|Du[c^*]|^2$ on the active sets. An $\epsilon\to0$ limit is not addressed here. %

The evolution maps associated with the $\mathrm{EVI}_\lambda$-gradient flows are approximated by interpolants obtained by incremental minimizations using the \emph{method of minimizing movements} \cite{degiorgi}; at iteration $k$, we define iteratively
\begin{equation}\label{eq:mm}
c_{k+1} \in \argmin_{c\in C} \,\Ee(c) + \frac{1}{2\tau_k} d(c_k, c)^2,
\end{equation}
where $\tau_k$ is the current time step, and we define the left-continuous piecewise constant interpolant
\begin{equation*}%
\overline{c}_{\btau}(t) :=
\begin{cases}
c_0, &\quad t=0,\\
c_{k+1} &\quad t \in(t_k, t_{k+1}],
\end{cases}
\end{equation*}
where discrete times are given recursively by $t_{k+1} = t_k + \tau_k$, with $t_0 = 0$. For a fixed time step $\tau$ and $\lambda>0$, the discrete estimate of \cite[Lemma~3.4.1 and Corollary~3.4.2]{ambrosio2005gradient} gives
\begin{equation}\label{eq:apriori_lambda_asym}
\frac{\lambda}{2}d(\overline c_{\btau}(t),c^*)^2
\le \Ee(\overline c_{\btau}(t))-\Ee(c^*)
\le \bigl(\Ee(c_0)-\Ee(c^*)\bigr)e^{-2\lambda_\tau t},
\end{equation}
where $\lambda_\tau:=\log(1+\lambda\tau)/\tau$. Thus the minimizing movements converge to the minimizer of the truncated problem with rate $\lambda_\tau$.

Central to the proof of the existence of the solutions of evolution variational inequalities is the notion of geodesical $\lambda$-convexity of a functional, which is stated below.
\begin{definition}[Geodesical $\lambda$-convexity]
We say that a functional $\J:C\to\Rinf$ is geodesically $\lambda$-convex with respect to a metric $d(\cdot,\cdot)$ if for all $c_0,c_1\in C$, there exists a geodesic curve $\gamma_g$ such that
\begin{equation*} %
\J(\gamma_g(s))
\le (1-s)\J(c_0)+s\J(c_1)-\frac{\lambda}{2}\,s(1-s)\,d(c_0,c_1)^2,\quad \forall s\in[0,1].
\end{equation*}
\end{definition}
Under the conditions that the space $(C,d)$ is a complete geodesic space and that the map $c\mapsto\frac{1}{2}d(\bar{c},c)^2$ is geodesically $1$-convex for all $\bar{c}\in C$, proving the existence of $\mathrm{EVI}_\lambda$-gradient flows reduces to proving the geodesical $\lambda$-convexity of $\Ee$ \cite{muratori2020gradient,daneri2010lecture,mayer1998gradient}. The required geodesical $1$-convexity is proved in Appendix Lemma~\ref{lem:d_1_convex_proof}.
The following theorem gives sufficient conditions for the geodesical $\lambda$-convexity of $\Ee$ (see Appendix Theorem~\ref{thm:Ee_lambda_conv_proof}).
\begin{theorem}\label{thm:Ee_lambda_conv}
Let Assumption~\ref{ass:sindices_strict} hold, together with Assumption~\ref{ass:F_Hm1} in the scalar case or Assumption~\ref{ass:F_Hm1_vector} in the incompressible vector-valued case, and assume $\inf_{r>0}\left\{1+\frac{r\mu''(r)}{2\mu'(r)}\right\}>-\infty$. Then
\begin{enumerate}
    \item[(i)] if $\Phi$ is convex and $\mu''(r)\le0$, or  
    \item[(ii)] if $\Phi$ is concave and $\mu''(r) \ge \frac{4 \mu'(r)^2}{\mu(r)}$,
\end{enumerate}
$\Ee$ is geodesically $\lambda$-convex for all $\epsilon \in(0,1)$ with $\lambda=\inf_{r>0}\left\{1+\frac{r\mu''(r)}{2\mu'(r)}\right\}$.
\end{theorem}

\begin{remark}\label{rem:powerlaw_lam}
For the power-law case, a straightforward calculation shows that Theorem \ref{thm:Ee_lambda_conv} applies when $\frac{4}{3}\le p < 2$ or $2<p\le 4$ with $\lambda = \frac{p}{4}>0$. The cases of $p<4/3$ and $p>4$ do not fall under the current theoretical framework.
\end{remark}

\section{Finite element discretization}\label{sec:impl}

In this section, we describe our finite element discretization strategy to approximate the minimum of the auxiliary energy. We denote by $\Ch$, $\Uh$, $\bUh$, and $Q_h$ the finite element function spaces associated with the discretization of $\LPhi$, $\Honez$, $\vHonez$, and $Q$, respectively, and we assume the following.
\begin{assumption}\label{ass:fem_space}
$C_h\subset\Linf$, $|\grad u_h|^2\in C_h$ for all $u_h \in \Uh$, $|\eps(\bu_h)|^2\in C_h$ for all $\bu_h \in \bUh$, and $\bUh \times Q_h$ is an inf-sup stable finite element pair for the weighted Stokes system.
\end{assumption}
Moreover, we use the common notation
\[
W_h :=
\begin{cases}
\Uh, &\text{scalar case}\\
\bUhdiv:=\displaystyle\left\{\bu_h \in \bUh \st \int_\Omega q_h\diver \bu_h  \dx = 0,\, \forall q_h \in Q_h\right\}, &\text{vector-valued case}
\end{cases}.
\]

Instead of directly discretizing the method of minimizing movements \eqref{eq:mm} (see Remark \ref{rem:mm_fem}), we consider the metric gradient flow approach described in \cite{liero2013gradient} (see also \cite[Section 4]{mielke2023introduction}). The derivation is performed entirely at the finite-dimensional level; at the continuous level, identifying the first
variation with an $L^2$-gradient through the mapped functional \eqref{eq:E_l2}, defined in the appendix, would require additional regularity of the state variable. In the finite-dimensional setting, instead, derivatives are elements of $C_h^*$, and gradients are obtained from them only after choosing the appropriate Riesz map. Unless explicitly stated, we do not need Assumption \ref{ass:sindices_strict} to be satisfied.

\subsection{Metric gradient flow}
We work on the open positive cone
\[
C_h^+:=\left\{c_h\in C_h\st\essinf_{x\in\Omega} c_h(x)>0\right\},
\]
that avoids the need to introduce one-sided admissible directions associated with box constraints (see Theorem \ref{thm:ode_bound}), and consider the smooth discrete energy
\begin{equation}\label{eq:Eh_def}
\Eh(c_h):=-\omega \Jc(c_h)-\omega \mathcal{J}_{u_h}(c_h),\quad c_h\in C_h^+,
\end{equation}
where $\mathcal{J}_{u_h}$ is identical to $\Ju$ except that the infimum of the Dirichlet problem is sought in $W_h$. Throughout this section, we always assume the validity of Assumption \ref{ass:F_Hm1} or \ref{ass:F_Hm1_vector}.

We then derive the first variation of the discrete energy for a fixed $c_h\in C^+_h$. We consider the scalar case first. Let $w_h\in C_h$ be such that $\evar w_h$ is an admissible direction. Recalling that $\mu(r)=2\Phi'(r)$, the first variation of $\Eh$ in the direction $w_h$ is 
\begin{align}
\DEh(c_h)[w_h]
&:=\left.\dd{}{\evar}\Eh(c_h+\evar w_h)\right|_{\evar=0}\nonumber\\
&=\frac{\omega}{2}\int_\Omega\mu'(c_h)\bigl(c_h+|\grad u_h[c_h]|^2\bigr)w_h\dx
+\omega\int_\Omega\mu(c_h)\grad u_h[c_h]\cdot\grad\auc_h\dx,\label{eq:first_var_temp}
\end{align}
where $\auc_h$ is the variation of $u_h[c_h]$ in the direction $w_h$, i.e.
\[
\auc_h:=\left.\dd{}{\evar}u_h[c_h+\evar w_h]\right|_{\evar=0}.
\]
Differentiating \eqref{eq:scalar_state} yields 
\begin{equation}\label{eq:sens_scalar}
\int_\Omega \mu(c_h)\grad \auc_h\cdot \grad v_h\dx
+\int_\Omega \mu'(c_h) \grad u_h[c_h]\cdot \grad v_h \,w_h \dx = 0,
\quad \forall v_h\in \Uh,
\end{equation}
that, when tested with $v_h=u_h[c_h]$, gives the identity
\begin{equation}\label{eq:key_identity_scalar}
\int_\Omega \mu'(c_h)|\grad u_h[c_h]|^2 w_h\dx
=
-\int_\Omega \mu(c_h)\grad \auc_h \cdot \grad u_h[c_h]\dx.
\end{equation}
Substituting \eqref{eq:key_identity_scalar} in \eqref{eq:first_var_temp}, we obtain
\begin{equation}\label{eq:first_var_scalar}
\DEh(c_h)[w_h]
=
\int_\Omega \frac{-\omega\mu'(c_h)}{2} \left(|\grad u_h[c_h]|^2-c_h\right)\,w_h\dx.
\end{equation}

We then compute the first variation of the auxiliary energy in the incompressible vector-valued case
\[
\DEh(c_h)[w_h]
=\int_{\Omega}
-\omega\left(
-\frac{1}{2}\mu'(c_h)c_h
-\frac{1}{2}\mu'(c_h)\,|\eps(\bu_h[c_h])|^2
\right)w_h
+\omega\mu(c_h)\,\eps(\bu_h[c_h])\fro\eps(\bauc_h)\dx
\]
where
\[
\bauc_h:=\left.\dd{}{\evar}\bu_h[c_h+\evar w_h]\right|_{\evar=0}.
\]
Using the notation
$
\pauc_h:=\left.\dd{}{\evar}p_h[c_h+\evar w_h]\right|_{\evar=0},
$
a differentiation of the Stokes constraint \eqref{eq:vector_state}--\eqref{eq:vector_incomp}
yields
\begin{align}
\int_{\Omega}\mu(c_h)\,\eps(\bauc_h)\fro\eps(\bv_h)\dx
-\int_{\Omega} \pauc_h\,\diver\bv_h\dx
&= -\int_{\Omega}\mu'(c_h)\,\eps(\bu_h[c_h])\fro\eps(\bv_h)\,w_h\dx,\quad\forall \bv_h \in \bUh,\label{eq:sens_vector}\\
\int_{\Omega} q_h\,\diver \bauc_h  \,\dx &=0,\quad \forall q_h\in Q_h.\label{eq:sens_vector_inc}
\end{align}
Testing the first equation with $\bv_h=\bu_h[c_h]$ and using the incompressibility of $\bu_h[c_h]$, we arrive at
\begin{equation}\label{eq:first_var_vector}
\DEh(c_h)[w_h]
:=
\int_\Omega \frac{-\omega\mu'(c_h)}{2} \left(|\eps(\bu_h[c_h])|^2-c_h\right)\,w_h\dx.
\end{equation}

The same first variation is obtained if homogeneous essential boundary
conditions are replaced by homogeneous or non-homogeneous natural boundary conditions. In that case the state problem is posed on the corresponding
linear space $W_h$, normalized, as usual, when the operator has a non-trivial kernel, and the right-hand side is replaced by
\[
F_N(v_h):=F(v_h)+\int_{\partial\Omega} g_N \,v_h\ds.
\]
Since $F_N$ is independent of $c_h$, differentiating the state equation gives the same linearized problems as above, and therefore the same cancellation. We then show how to handle the case of non-homogeneous essential boundary conditions; we consider the scalar case, as the argument for the incompressible vector-valued case is analogous. With non-homogeneous essential boundary conditions, $\Uh$ is an affine space and the state problem is 
\[
\int_\Omega \mu(c_h)\grad u_h[c_h]\cdot \grad v_h\dx=F(v_h),
\quad \forall v_h\in U_{h,0},
\]
where $U_{h,0}$ denotes the space with homogeneous boundary conditions. In this case %
\begin{equation*} %
\mathcal{J}_{u_h}(c_h)
=
\frac12\int_\Omega \mu(c_h)|\grad u_h[c_h]|^2\dx-F(u_h[c_h]),
\end{equation*}
with first variation
\[
\mathcal{DJ}_{u_h}(c_h)[w_h]
= \frac12\int_\Omega \mu'(c_h)|\grad u_h[c_h]|^2w_h\dx
+
\int_\Omega \mu(c_h)\grad u_h[c_h] \cdot \grad \auc_h\dx
-
F(\auc_h),
\]
where $\auc_h\in U_{h,0}$ since the lifting is independent of $c_h$. Therefore, $\auc_h$ is an admissible test function in the state equation, giving the identity
\[
\int_\Omega \mu(c_h)\grad u_h[c_h]\cdot \grad\zeta_h\,dx
=
F(\auc_h),
\]
and thus yielding the same first variation \eqref{eq:first_var_scalar}.

\begin{remark}\label{rem:mm_fem}
The discretization of the method of minimizing movements given in \eqref{eq:mm} reads
\[
c_{h,n+1}\in\argmin_{c_h\in C_{h,\epsilon}}\left\{\Eh(c_h)+\frac{1}{2\tau}\int_\Omega|\psi(c_h)-\psi(c_{h,n})|^2\dx\right\}, \quad
C_{h,\epsilon}=\{c_h\in C_h:\epsilon\leq c_h\leq \epsilon^{-1}\}.
\]
A formal derivation of the Euler--Lagrange equations, ignoring the box constraints, leads to
\[
\int_\Omega\psi'(c_{h,n+1})\left[\frac{\psi(c_{h,n+1})-\psi(c_{h,n})}{\tau}-\psi'(c_{h,n+1})\bigl(|D u_h[c_{h,n+1}]|^2-c_{h,n+1}\bigr)\right]w_h\dx=0, \quad \forall w_h\in C_h, 
\]
which can be read as an implicit Euler scheme in the metric coordinate $\eta_h=\psi(c_h)$
\[
\frac{\eta_{h,n+1}-\eta_{h,n}}{\tau}=\psi'(c_{h,n+1})\bigl(|D u_h[c_{h,n+1}]|^2-c_{h,n+1}\bigr),\quad
\eta_{h,n+1}=\psi(c_{h,n+1}),
\]
tested against variations induced by $w_h\in C_h$.
\end{remark}

Following \cite{liero2013gradient}, the differential $\DEh(c_h)\in C_h^*$ is turned into a gradient $\grad_{g_{c_h}} \Eh(c_h) \in C_h$ by solving the variational problem
\[
g_{c_h}(\grad_{g_{c_h}} \Eh(c_h), w_h) = \DEh(c_h)[w_h], \quad \forall w_h \in C_h,
\]
where the inner product arising from the discrete metric tensor is
\begin{equation*}\label{eq:gch_inner_product}
g_{c_h}(z_h,w_h):=\bigl(\mathcal{D}\Psi(c_h)z_h,\mathcal{D}\Psi(c_h)w_h\bigr)_{\Ltwo} = \int_\Omega\frac{\omega}{2}\mu'(c_h)z_h w_h\dx.
\end{equation*}
Note that, since $C_h$ is a finite-dimensional linear space, its tangent space at
each point is canonically $C_h$ itself. Moreover, under the sole Assumption \ref{ass:Phi_concace_or_convex_strictly}, $\omega\mu'(c_h)>0$ on $C_h^+$, and therefore $g_{c_h}(\cdot,\cdot)$ defines an inner product on $C_h$ for all $c_h \in C^+_h$. Note that this distinction between derivatives and gradients is standard in
Hilbert-space PDE-constrained optimization
\cite{HinzePinnauUlbrichUlbrich2009} and  
Bayesian inverse problems
\cite{Stuart2010}.

In the present setting, the metric tensor, equivalently the Riesz map
\begin{equation}\label{eq:def_G}
\mathcal{G}(c_h) : C_h\rightarrow C^\ast_h,\quad\langle\mathcal{G}(c_h)w_h,z_h\rangle = g_{c_h}(z_h,w_h), %
\quad c_h \in C^+_h,\,w_h,z_h \in C_h 
\end{equation}
is explicitly known, and we do not need to solve any additional variational problem to define the gradient. This is in contrast with many gradient-flow equations, where the so-called Onsager operator $K_h(c_h):C_h^*\to C_h$, mapping forces to velocities, is specified first and the metric tensor is only implicitly defined via its inverse. %

Using the first variation of the functional given in \eqref{eq:first_var_scalar} for the scalar case, \eqref{eq:first_var_vector} for the vector-valued case, and the shorthand notation $Du$ to denote $\grad u$ in the scalar case and $\eps(\bu)$ in the vector-valued case, the metric gradient flow approach therefore gives
\[
\grad_{g_{c_h}} \Eh(c_h) = c_h - |D u_h[c_h]|^2,
\]
and it allows us to consider the ordinary differential equation 
\begin{equation}\label{eq:ode_fem}
\left\{
\begin{aligned}
\dd{c_h}{t} &= |D u_h[c_h]|^2 - c_h,\\
c_h(0) &= c_{h,0} \in C^+_h.
\end{aligned}
\right.
\end{equation}
Assumption \ref{ass:fem_space} ensures that the pointwise ODE is equivalent to the semi-discrete problem
\begin{equation}\label{eq:ode_semi}
\left\{
\begin{aligned}
\int_\Omega \dd{c_h}{t} w_h\dx &= \int_\Omega(|D u_h[c_h]|^2 - c_h) w_h\dx, \quad \forall w_h \in C_h,\\
c_h(0) &= c_{h,0} \in C^+_h.
\end{aligned}
\right.
\end{equation}
Moreover, $\Eh$ is a Lyapunov function, since
\[
\dd{}{t}\Eh(c_h(t)) = -\int_\Omega\frac{\omega}{2}\mu'(c_h(t))\,\dot{c}_h^{\,2}(t)\dx\le 0, \quad \dot{c}_h :=\dd{c_h}{t}.
\]

Box constraints are only needed to define the continuous functional $\Ee$. At the discrete level, however, the following result shows that they can be neglected for finite-time evolutions, and no explicit projection onto the box constraint is needed in the finite element implementation. The proof is given in Appendix Theorem~\ref{thm:ode_bound_proof}.
\begin{theorem}\label{thm:ode_bound}
Under Assumptions \ref{ass:Phi_concace_or_convex_strictly}, \ref{ass:sindices}, and \ref{ass:fem_space}, problem \eqref{eq:ode_fem} has a unique global solution $c_h\in C^1([0,+\infty);C_h)$ and $c_h(t) \in C^+_h$.
Moreover, for each fixed $T <+\infty$, there exists $\epsilon(T)> 0$ such that
\[
\epsilon(T) \le c_h(t)\le 1/\epsilon(T), \quad \forall\, 0\le t\le T.   
\]
\end{theorem}

We now state conditions under which the solution of the ODE \eqref{eq:ode_fem} converges to the FEM solution minimizing the energy \eqref{eq:scalar_primal} or \eqref{eq:vector_min}, so that approximation properties are maintained. The convergence argument has two ingredients: a metric Hessian bound \cite[Theorem 3.5]{liero2013gradient} that controls the
auxiliary residual, and natural-distance estimates \cite{DieningKreuzer2008} that control the state error by this residual.
The metric Hessian lower bound is proved under sufficient conditions in the next theorem. Assumptions on $\mu$ and estimates for $\lambda$ are identical to those in the continuous case in Theorem \ref{thm:Ee_lambda_conv}, and Remark \ref{rem:powerlaw_lam} translates directly to the discrete case. See Appendix Theorem~\ref{thm:mielke_lambda_discrete_proof} for the proof.
\begin{theorem}\label{thm:mielke_lambda_discrete}
Let the assumptions of Theorem \ref{thm:ode_bound} hold; assume further that
$\inf_{r>0}\left\{1+\frac{r\mu''(r)}{2\mu'(r)}\right\}>-\infty$. Then  
\begin{enumerate}
    \item[(i)] if $\Phi$ is convex and $\mu''(r)\le0$, or  
    \item[(ii)] if $\Phi$ is concave and $\mu''(r) \ge \frac{4 \mu'(r)^2}{\mu(r)}$,
\end{enumerate}
the following inequality holds
with $\lambda=\inf_{r>0}\left\{1+\frac{r\mu''(r)}{2\mu'(r)}\right\}$
\[
\mathcal{B}(c_h,z_h) \ge \lambda \mathcal{A}(c_h,z_h) \quad \forall c_h\in C^+_h,\quad \forall z_h \in C_h,
\]
where
\begin{align}
\mathcal{A}(c_h,z_h)&:=\langle \mathcal{G}(c_h)z_h,z_h\rangle
=\int_\Omega \frac{\omega}{2}\,\mu'(c_h)\,z_h^2\,\dx, \label{eq:Adef}\\
\mathcal{B}(c_h,z_h)&:=\langle \mathcal{G}(c_h)z_h,\mathcal{DR}_h(c_h)z_h\rangle+\frac12\langle \mathcal{DG}(c_h)[\mathcal{R}_h(c_h)]z_h,z_h\rangle,\label{eq:Bdef}
\end{align}
with $\mathcal{G}(c_h)$ defined by \eqref{eq:def_G} and $\mathcal{R}_h(c_h) = c_h-|Du_h[c_h]|^2$.
\end{theorem}

\begin{remark}\label{rem:powerlaw_lam_fem}
Remark \ref{rem:powerlaw_lam} for the power-law case translates directly to the FEM problem. Making the assumption that, for large $t$, $|D u_h[c_h(t)]|^2 \approx c_h(t)$ holds, we can obtain \emph{asymptotic} estimates for $\lambda$ for $p>4$ and $p<4/3$. For the case $p>4$, from Equation \eqref{eq:convex_lambdach}, we obtain $\lambda=1$ since
\[
\pos{\frac{\mu''(r)}{2\mu'(r)}}r=\frac{p-4}{4}.
\]
When $p<4/3$, using Equation \eqref{eq:concave_lambdach}, we find $\lambda=p-1$, since
\[
\pos{\frac{\mu''(r)}{2\mu'(r)}-\frac{2\mu'(r)}{\mu(r)}}r = \frac{4-3p}{4}.
\]
\end{remark}

We next collect the natural-distance estimates needed to pass from the
auxiliary residual to the state error. We use results from Diening and Kreuzer \cite{DieningKreuzer2008}, which generalized the quasi-norm convergence estimates initially developed for the $p$-Laplacian by Ebmeyer and Liu \cite{ebmeyer2005quasi}. In particular, we assume that the error of the Galerkin approximation of the Euler--Lagrange equations 
\begin{equation}\label{eq:galerkin_fem}
\int_\Omega A(D\bar u_h):D v_h\dx= F(v_h), \quad \forall v_h\in W_h,
\end{equation}
is measured in terms of
\begin{equation}\label{eq:errornorm_fem}
\|V(D\bar{u}) - V(D\bar{u}_h)\|_\Ltwo,
\end{equation}
where $\bar{u}$ is the minimizer of problem \eqref{eq:scalar_primal} or \eqref{eq:vector_min}, and
\[
A(P):=\mu(|P|^2)P=\frac{\phi'(|P|)}{|P|}P,\quad V(P):=\sqrt{\mu(|P|^2)}\,P, \quad P\in\mathbb{R}^{Nd},
\]
with continuous extension $A(0)=0$ and $V(0)=0$, and $N$ is $1$ or $d$ depending on whether we consider the scalar or the vector-valued case.
We match the notation of \cite{DieningKreuzer2008} and use $f\lesssim g$ when $f,g \ge0$ and
$f\le C g$, where $C>0$ is a constant that only depends on $\smphi$, $\spphi$ and the $\Delta_2$ constants of $\phi$ and $\phi^*$. We write $f\simeq g$ if both $f\lesssim g$ and $g\lesssim f$.

For $a\ge0$, we also define the shifted N-function $\phi_a$ by
\begin{equation}\label{eq:shifted_nfunc}
\frac{\phi_a'(s)}{s}:=\frac{\phi'(a+s)}{a+s},\quad s>0,\quad\phi_a(0)=0,
\end{equation}
and first state the following lemma, which recalls results from \cite{DieningKreuzer2008}. %
\begin{lemma}\label{lem:DK-shifted}
Let $\phi_a$ be defined by \eqref{eq:shifted_nfunc} for a given $a\ge 0$. Then, under Assumption \ref{ass:sindices}
\begin{enumerate}
\item[(i)]
$\phi_a'(s)\simeq s\phi''(a+s),$
\item[(ii)]
$\phi_a(s)\simeq s^2\frac{\phi'(a+s)}{a+s}\simeq s^2\phi''(a+s),$
\item[(iii)] 
$(A(P)-A(Q)):(P-Q) \simeq |V(P)-V(Q)|^2 \simeq \phi_{|P|}(|P-Q|), \quad \forall\, P, Q \in \mathbb{R}^{Nd},$
\item[(iv)]  For every $\delta>0$, there exists $C_\delta$ depending only on $\delta$, $\smphi$, $\spphi$ and the uniform $\Delta_2$ constant of the family $\{\phi_a^*\}_{a\ge0}$ such that
\[
tr \le \delta \phi_a(t)+C_\delta \phi_a^*(r), \quad t,r\ge0.
\]
\end{enumerate}
\end{lemma}
Statements (i) and (ii) of the previous lemma yield the following intermediate result (see Appendix Lemma~\ref{lem:coefficient-residual_proof}).
\begin{lemma}\label{lem:coefficient-residual}
Let the assumptions of Lemma \ref{lem:DK-shifted} hold. Then, if %
\[
\frac13\le\smphi\le\spphi<1, \quad \text{or}\quad 1<\smphi\le\spphi\le3,
\]
we have
\begin{equation}\label{eq:phistar_P_bound}
\phi^*_{|P|}\left(|\mu(r)-\mu(|P|^2)|\,|P|\right)
\lesssim|\mu'(r)|\,(r-|P|^2)^2,\quad \forall r>0,\, P \in \mathbb{R}^{Nd}.
\end{equation}
\end{lemma}

\begin{remark}
The curvature conditions in Theorem~\ref{thm:mielke_lambda_discrete}
are sufficient, but not necessary in general, for the corresponding
index thresholds $\smphi\ge1/3$ and $\spphi\le3$ in
Lemma~\ref{lem:coefficient-residual}.
In fact, observe that any positive differentiable
concave function $g$ on $(0,+\infty)$ satisfies $rg'(r)\le g(r)$.
If $\mu''\le0$, then $\mu$ is positive and concave. The preceding
argument therefore gives $r\mu'(r)/\mu(r)\le1$.
If $\mu\mu''\ge4(\mu')^2$, we instead consider the positive concave function
$\mu^{-3}$. %
The same argument now gives $-3 r\mu'(r)/\mu(r)\le1$.
Rearranging \eqref{eq:mup_exact_final} gives
\[
\frac{t\phi''(t)}{\phi'(t)}
=1+2\frac{r\mu'(r)}{\mu(r)},\qquad t=\sqrt r.
\]
Taking the supremum over $t>0$ in the first case yields $\spphi\le3$,
while taking the infimum in the second case yields $\smphi\ge1/3$.
For power laws, the curvature conditions and the corresponding thresholds are equivalent.
\end{remark}

Introducing the squared metric speed of the auxiliary flow,
\begin{equation}\label{eq:Dht_def}
d_h(t):=g_{c_h(t)}(\dot{c}_h(t),\dot{c}_h(t))=\frac{1}{2}\int_\Omega |\mu'(c_h(t))|(c_h(t)-|Du_h[c_h(t)]|^2)^2\dx,
\end{equation}
the preceding results give
\[
d_h(t)\le e^{-2\lambda t}d_h(0),
\qquad
\|V(Du_h[c_h(t)])-V(D\bar u_h)\|_\Ltwo^2\lesssim d_h(t).
\]
The first estimate follows by differentiating $d_h$ and applying
Theorem~\ref{thm:mielke_lambda_discrete}; the second follows from
Lemmas~\ref{lem:DK-shifted} and~\ref{lem:coefficient-residual}
and it holds for any $c_h\in C_h^+$ and its reconstructed state $u_h[c_h]$, independently of the evolution.
Therefore, the weighted residual controls the natural error relative to the nonlinear Galerkin solution.
Combining them,
we obtain the following convergence result (see Appendix Theorem~\ref{thm:ode_fem_convergence_proof}).
\begin{theorem}\label{thm:ode_fem_convergence}
Let Assumption~\ref{ass:sindices_strict} and the assumptions of Theorems~\ref{thm:ode_bound} and~\ref{thm:mielke_lambda_discrete} be satisfied. Then there exists a constant $C>0$, depending only on $\smphi$, $\spphi$, and the $\Delta_2$ constants of $\phi$ and $\phi^*$, such
that
\[
\|V(Du_h[c_h(t)])-V(D\bar u_h)\|_\Ltwo \le C e^{-\lambda t}\sqrt{d_h(0)},
\]
where $c_h(t)$ is the solution of the ODE \eqref{eq:ode_fem}, $\bar u_h$ is the solution of \eqref{eq:galerkin_fem} and $\lambda$ is given in Theorem \ref{thm:mielke_lambda_discrete}. Therefore, if $\lambda>0$, we have
\[
\lim_{t\to+\infty}\|V(Du_h[c_h(t)])-V(D\bar u_h)\|_\Ltwo=0.
\] 
\end{theorem}

The natural distance $\|V(Du)-V(Dv)\|_\Ltwo$ compares gradients (resp. symmetric gradients) of scalar (resp. vector) fields and measures the state error in Theorem~\ref{thm:ode_fem_convergence}. To relate it to our coefficient metric $d$, consider
\begin{equation}\label{eq:chi}
\chi(c):=\sqrt{c\mu(c)},\quad c>0,\qquad \chi(0):=0.
\end{equation}
Under Assumption~\ref{ass:sindices_strict}, $\psi'\simeq\chi'$, so $d(c_0,c_1)\simeq\|\chi(c_1)-\chi(c_0)\|_\Ltwo$ for $c_0,c_1\in\LPhip$. Since $\chi(|P|^2)=|V(P)|$, this gives
\begin{equation}\label{eq:metric_natural_magnitudes}
d(|Du|^2,|Dv|^2)\simeq\bigl\||V(Du)|-|V(Dv)|\bigr\|_\Ltwo
\le\|V(Du)-V(Dv)\|_\Ltwo.
\end{equation}
Thus the equivalence concerns gradient magnitudes, while the full natural distance also accounts for their directions.
We choose $\psi$ because its weight $\psi'(c_h)^2=|\mu'(c_h)|/2$ cancels the factor in the first variation of the auxiliary energy, yielding the simple evolution~\eqref{eq:ode_fem}. Replacing $\psi$ by $\chi$ in \eqref{eq:X_distance} would formally give
\[
\dd{c_h}{t}=\frac{2c_h\mu(c_h)|\mu'(c_h)|}{(\mu(c_h)+c_h\mu'(c_h))^2}\bigl(|Du[c_h]|^2-c_h\bigr).
\]
For power laws with $p\ne2$, the prefactor reduces to the constant $4|p-2|/p^2$, so these two coefficient flows differ only by a rescaling of time.
For general viscosity laws, the inclusion $|Du_h|^2\in C_h$ in Assumption \ref{ass:fem_space} alone does not ensure equivalence between the pointwise ODE
induced by $\chi$ defined in~\eqref{eq:chi}
and its weak formulation: the additional prefactor may take the right-hand side outside $C_h$. %

\subsection{Space and time discretizations}\label{sec:impl_time}

We first consider the space discretization. In the case of a simplicial element $T\subset\R^d$, we use $\Pk(T)$ to denote the space of polynomials in $d$ variables with degree at most $k$.
Given a polynomial degree $k$, the choice of the finite element space for the discretization of functions in $\Honez$ is standard
\[
U^{k,T}_h := \big\{ u_h \in\Honez \st u_h \in C^0(\Omega), {u_h}_{|T} \in \Pk(T)\big\}.
\]
Given $u_h \in U^{k,T}_h$, we have $\grad {u_h}_{|T} \in \Pkm(T)$, and thus $|\grad {u_h}|^2$ is a polynomial of degree $2(k-1)$ on each $T$; a suitable discretization space for functions in $\LPhi$ that satisfies Assumption \ref{ass:fem_space} is then
\[
C^{k,T}_h := \big\{ c_h \in\LPhi \st {c_h}_{|T} \in \Pkc(T)\big\}.
\]
Similar considerations hold using tensor-product cells $K$. In such a case, we use
\[
U^{k,K}_h := \big\{ u_h \in\Honez \st u_h \in C^0(\Omega), {u_h}_{|K} \in \Qk(K)\big\},
\]
where $\Qk(K)$ denotes the space of polynomials in $d$ variables with degree at most $k$ for each variable. Since $\grad {u_h}_{|K} \in \Qk(K)$, a sufficient choice for the discretization space for functions in $\LPhi$ is
\[
C^{k,K}_h := \big\{ c_h \in\LPhi \st {c_h}_{|K} \in \Q{2k}(K)\big\}.
\]

The same matching conditions hold true when $\grad u_h$ is replaced by the symmetric gradient $\eps(\bu_h)$; in the case of vector-valued spaces, we will consider the Taylor--Hood pairs with $k\ge2$
\[
\begin{aligned}
\bm{U}_h^{k,T} &:= \big\{ \bu_h \in \vHonez \st \bu_h \in C^0(\Omega;\mathbb{R}^d), {\bu_h}_{|T} \in \Pk(T)^d\big\},\\
Q_h^{k,T} &:= \big\{ p_h \in Q \st p_h \in C^0(\Omega), \,{p_h}_{|T} \in \P{k-1}(T)\big\}.
\end{aligned}
\]
or, in the case of tensor-product cells
\[
\begin{aligned}
\bm{U}_h^{k,K} &:= \big\{ \bu_h \in \vHonez \st \bu_h \in C^0(\Omega;\mathbb{R}^d), {\bu_h}_{|K} \in \Qk(K)^d\big\},\\
Q_h^{k,K} &:= \big\{ p_h \in Q \st p_h \in C^0(\Omega),\, {p_h}_{|K} \in \Q{k-1}(K)\big\}.
\end{aligned}
\]

\begin{remark}
The finite element spaces described above are only one possible choice. 
Any inf--sup stable pair $\bUh\times Q_h$ may be used for the
Stokes discretization, provided that $C_h$ is chosen rich enough to contain
$|\eps(\bu_h)|^2$ for all $\bu_h\in \bUh$.
\end{remark}

We then describe the time discretization for the scalar problem. 
The modifications required for the incompressible vector-valued case are given at the end of the section. Throughout, the subscript $h,n$ denotes finite element functions at the $n$-th iteration. Also, we drop the cell and order superscripts and always write $\Uh$, $\Ch$, $\bUh$, $Q_h$ for the finite element spaces, assuming matching order conditions.

For a given initial condition $c_{h,0}\in \Ch^+$, we first compute the initial state $u_{h,0}$ by solving the weighted Poisson problem: Find $u_{h,0} \in  \Uh$ such that for all $v_h \in \Uh$
\[
a[c_{h,0}](u_{h,0}, v_h) = F(v_h),
\]
where $a[\cdot](\cdot,\cdot)$ is given by \eqref{eq:scalar_bilin}.
As in our previous work on biological transportation networks \cite{zampiniBio}, we use the backward Euler method to discretize \eqref{eq:ode_semi} in time and consider the nonlinear system of equations: Given $(u_{h,n},c_{h,n})$ and a current time step $\tau_n$, find $(u_{h,n+1},c_{h,n+1})$ such that  for all $(v_h,w_h) \in \Uh \times\Ch$
\begin{equation}\label{eq:be_mm}
\begin{aligned}
\int_\Omega\frac{c_{h,n+1} - c_{h,n}}{\tau_n}\,w_h\dx &= \int_\Omega(|\grad u_{h,n+1}|^2 - c_{h,n+1})\,w_h \dx,\\
a[c_{h,n+1}](u_{h,n+1}, v_h) &= F(v_h).
\end{aligned}
\end{equation}
Note that the backward Euler method ensures the strict positivity of $c_{h}$ throughout the time iteration, since
\begin{equation}\label{eq:be_positive}
c_{h,n+1} = \frac{\tau_n}{1+\tau_n}|\grad u_{h,n+1}|^2 + \frac{1}{1+\tau_n}c_{h,n}. %
\end{equation}
Although strict positivity of $c_h$ is preserved throughout the time iteration, the nonlinear system for time advancement can be difficult to solve directly, because Newton iterates need not preserve positivity before convergence.

When only the steady state is of interest, one may instead consider the split iteration
\begin{equation}\label{eq:be_split}
\begin{aligned}
\int_\Omega\frac{c_{h,n+1} - c_{h,n}}{\tau_n}\,w_h\dx &= \int_\Omega(|\grad u_{h,n}|^2 - c_{h,n+1})\,w_h \dx,\\
a[c_{h,n+1}](u_{h,n+1}, v_h) &= F(v_h),
\end{aligned}
\end{equation}
which completely decouples the two equations and, given our choice of finite element spaces, allows us to replace the first variational equation for the time advancement of $c_{h}$ with the exact interpolation formula
\begin{equation}\label{eq:be_interp}
c_{h,n+1} = \frac{\tau_n}{1+\tau_n}|\grad u_{h,n}|^2 + \frac{1}{1+\tau_n}c_{h,n}.
\end{equation}

\Needspace{5\baselineskip}
\begin{remark}\label{rem:kac_split}
Using $\tau_n=+\infty$ in \eqref{eq:be_interp}, or equivalently a forward Euler scheme with $\tau_n=1$, leads to the \kac iteration scheme \cite{han1997kavcanov, kavcur1968convergence}, which, as noted for example in \cite[Remark 21]{diening2020relaxed}, can fail for power laws with $p>2$; furthermore, using the update $c_{h,n+1} = |\grad u_{h,n}|^2$ requires dynamic clipping of the weights of the Poisson equation since the strict positivity of $|\grad u_{h,n}|^2$ cannot be guaranteed.
\end{remark}

In numerical experiments, we found that the split-iteration scheme does not always converge to the system's steady state and may require small time steps to avoid divergence.
A robust alternative we have found, for which we do not currently have any proof of reliable convergence, borrows ideas from the pseudo-transient continuation scheme \cite{kelley1998convergence} and its extension to differential-algebraic equations \cite{coffey2003pseudotransient}. In our pseudo iteration scheme, we first compute an intermediate value for $c_h$, denoted by $c_{h,n+1/2}$, by using the interpolation formula \eqref{eq:be_interp}, and then perform a single Newton step by linearizing \eqref{eq:be_mm} around $(c_{h,n+1/2},u_{h,n})$
\[
\begin{bmatrix}
(\tau_n^{-1}+1) M & -2 B^n\\
C^n & K^n
\end{bmatrix}
\begin{bmatrix}
\delta c_h\\
\delta u_h
\end{bmatrix}
=
\begin{bmatrix}
0\\
g
\end{bmatrix},
\]
where the matrix and vector entries are
\[
\begin{aligned}
M_{i,j} &= \int_\Omega \phi_i\phi_j \dx, & B^n_{i,j} &= \int_\Omega \phi_i \grad u_{h,n}\cdot \grad \varphi_j\dx, \\
C^n_{i,j} &= \int_\Omega \mu'(c_{h,n+1/2})\,\phi_j \grad u_{h,n}\cdot \grad \varphi_i\dx, &
K^n_{i,j} &= \int_\Omega \mu(c_{h,n+1/2})\,\grad \varphi_i \cdot \grad \varphi_j\dx,\\
g_i &=  F(\varphi_i) - \int_\Omega \mu(c_{h,n+1/2})\,\grad u_{h,n} \cdot \grad \varphi_i\dx,
\end{aligned}
\]
with $\phi_i$ and $\varphi_i$ the basis functions for $\Ch$ and $\Uh$, respectively. Note that the first component of the residual is zero because of the choice of the partial update.
The coupled correction $(\delta c_h,\delta u_h)$ is damped with a standard residual-based backtracking line search applied to the coupled backward-Euler system \eqref{eq:be_mm}, including its time-derivative term.
Unlike the fully implicit \eqref{eq:be_mm} and split \eqref{eq:be_split} steps, the new iterate is not guaranteed to remain positive unless negative values produce non-finite residuals.
We also note here that, in the fully nonlinear case \eqref{eq:be_mm}, the Newton iteration always starts from $(c_{h,n+1/2},u_{h,n})$.

\begin{remark}\label{rem:schur}
The pseudo-time stepping strategy enjoys a very favorable property, since the Schur complement obtained after eliminating the $c_h$ variable
\begin{equation}\label{eq:schur}
S^n:=K^n+2\alpha_n C^nM^{-1}B^n, \quad \alpha_n:=\tau_n/(1+\tau_n),
\end{equation}
is symmetric and positive definite provided $c_{h,n}\in C_h^+$ (see Appendix Lemma~\ref{lem:schur}). This implies that we can leverage scalable preconditioning strategies, such as Algebraic Multigrid \cite{xu2017algebraic} or Balancing Domain Decomposition methods \cite{zampini2016pcbddc}, to efficiently solve the linear problem in the pseudo-time stepping algorithm. %
\end{remark}

Time advancement for the incompressible vector-valued case is analogous to that for the scalar case, with the difference that the weighted Poisson problem is replaced by a weighted Stokes problem, and the equations for $c_h$ involve the symmetric gradient. We state here only the fully implicit backward Euler step:  Given $(\bu_{h,n},p_{h,n},c_{h,n})$ and a current time step $\tau_n$, find $(\bu_{h,n+1},p_{h,n+1},c_{h,n+1})$ such that for all $(\bv_h,q_h,w_h) \in \bUh \times Q_h\times\Ch$
\[
\begin{aligned}
\int_\Omega\frac{c_{h,n+1} - c_{h,n}}{\tau_n}\,w_h\dx &= \int_\Omega(|\eps(\bu_{h,n+1})|^2 - c_{h,n+1})\,w_h \dx,\\
\ba[c_{h,n+1}](\bu_{h,n+1}, \bv_h) - b(\bv_h,p_{h,n+1})&= F(\bv_h),\\
b(\bu_{h,n+1},q_h)&=0,
\end{aligned}
\]
where $\ba[\cdot](\cdot,\cdot)$ and $b(\cdot,\cdot)$ are given in \eqref{eq:vector_bilin}. The matrices and the right-hand side defining the pseudo update are the same, except that inner products are replaced by Frobenius inner products, symmetric gradients are used in place of plain gradients, and the matrix blocks accounting for the incompressibility condition of the weighted Stokes problem are included. Under the same positivity assumption, the velocity Schur complement obtained by eliminating the $c_h$ variable in the pseudo iteration scheme is also positive definite, and a BDDC method can be employed \cite{bevilacqua2024bddc} to solve the modified constrained equations.

\subsection{Adaptive time step selection}\label{sec:time_adapt}
Since we are only interested in the steady state of the system, in our numerical experiments, we will also make use of a heuristic adaptive time step strategy. The time step is chosen adaptively using numerical residual norms
\[
\rho_h^c(c_h,u_h):=\|\mathcal{R}_h^{c,*}(c_h,u_h)\|_{\ell^2},
\quad
\rho^f_h(c_h,u_h):=\|\mathcal{R}_h^{f,*}
(c_h,u_h)\|_{\ell^2}
\]
where $\mathcal{R}_h^{c,*}(c_h,u_h)$ is the assembled residual vector corresponding to the right-hand side of the ODE, and $\mathcal{R}_h^{f,*}
(c_h,u_h)$ the assembled vector including also the constraint residual.
We reject the time step if
\[
E(u_{h,n+1})> E(u_{h,n})+10^{-4}\text{ and }
    \rho_{1,n}^c>\rho_{1/2,n}^c,
\]
where $E$ is the original energy \eqref{eq:scalar_primal_Phi} and $\rho_{j,n}^c:=\rho_h^c(c_{h,n+j},u_{h,n+j})$, with the convention $u_{h,n+1/2} = u_{h,n}$.
This means that we accept an increase in the energy $E$ if it reduces the residual norm of the right-hand side of the ODE and, vice versa, accept the step if an increase in the ODE right-hand side leads to a decrease in the original energy. In case of rejection, we restore the previous iterate and set $\tau_{n+1}=0.1\,\tau_n$.

If the step is accepted, the next time step value is selected according to the rule
\[
\tau_{n+1}=
\begin{cases}
\tau_n,& \rho_{1,n}^c\leq \rho_{1/2,n}^c,\\
\displaystyle\tau_n \min\left\{\frac{\rho^f_{1/2,n}}{\rho^f_{1,n}},2\right\},&\rho_{1,n}^c>\rho_{1/2,n}^c,
\end{cases}
\]
where $\rho_{j,n}^f:=\rho_h^f(c_{h,n+j},u_{h,n+j})$.
The $c$-residual measures violation of $c_h=|Du_h|^2$ and determines whether rescaling is activated. As in conventional pseudo-transient continuation, the second indicator is the fully coupled residual, which includes the residuals of both the auxiliary equation and the discrete state equation. Its change sets the rescaling factor, with time-step growth limited to a factor of two.

\section{Numerical experiments}\label{sec:experiments}

In this section, we present numerical results using the {\tt Firedrake} Python package \cite{FiredrakeUserManual}; the code used for the experiments and the generation of figures and tables is public\footnote{\url{https://github.com/stefanozampini/generalized\_newtonian}}.
The initial condition of the ODE \eqref{eq:ode_fem} is always the constant function $c_{h,0}(x)=1$. Nonlinear problems are solved using damped Newton methods in {\tt PETSc} \cite{petsc-user-ref, petsc-web-page} via the {\tt petsc4py} Python bindings \cite{dalcinpazklercosimo2011}. For the nonlinear problem \eqref{eq:galerkin_fem}, Newton's method starts from $u_h[c_{h,0}]$ \cite{luo2016effective} and uses a cubic backtracking line search. Linear systems are solved using the {\tt MUMPS} factorization package \cite{MUMPS}.
Results are reported for the unscaled auxiliary energy
\[
\mathcal{E}(c_h) := \int_\Omega \Phi(c_h) -\frac{1}{2}\mu(c_h)c_h + \frac{1}{2}\mu(c_h)|D u_h[c_h]|^2\dx - F(u_h[c_h]),
\]
which is equivalent, modulo $-\omega$ scaling, to $\mathcal{E}_h$ given in \eqref{eq:Eh_def}. %

\subsection{Validation of the auxiliary energy and \texorpdfstring{$\lambda$}{lambda}-convexity estimates}

\subsubsection{Power-law model: scalar case}

In the first set of experiments, we validate the auxiliary energy and our geodesical $\lambda$-convexity estimates using the power-law model. For the scalar case, we consider the \emph{bump} example from \cite[Section 6]{diening2020relaxed}, which uses the polynomial exact solution $u_e(x) = (x_1^2 - 1)(x_2^2 - 1)$ on the domain $\Omega=(-1,1)^2$. As noted in \cite{diening2020relaxed}, $\diver(|\grad u_e|^{p-2}\grad u_e) \in \Lpp$ if and only if $p>\sqrt{2}$, with $p' := \frac{p}{p-1}$ the H\"older conjugate of $p$. Still, we can obtain a well-posed problem for all $p > 1$  by using the datum
\[
F(v) = \int_\Omega |\grad u_e|^{p-2}\grad u_e \cdot \grad v \dx,
\]
since $|\grad u_e|^{p-2}\grad u_e \in C^{0,\alpha}(\Omega;\mathbb{R}^d)$ with $\alpha=\min\{p-1,1\}$, given that $|\grad u_e|$ vanishes at $(\pm 1, \pm1)$ and $(0,0)$, and around each such point $x_0$ we have $|\grad u_e(x)|\approx|x-x_0|$, and thus the right-hand side satisfies the regularity requirements from Assumption \ref{ass:F_Hm1}.

We use a structured $32\times 32$ grid of tensor-product cells and approximation order $k=2$.
This choice allows us to exactly represent the solution and to compare the expected convergence rate of the method of minimizing movements \eqref{eq:apriori_lambda_asym}
with the observed convergence rates of the time-advancement methods introduced in Section \ref{sec:impl}, namely the backward Euler method (BE) \eqref{eq:be_mm}, the split iteration method (Split) \eqref{eq:be_split}, and the pseudo-transient continuation method (Pseudo). From Remark \ref{rem:powerlaw_lam}
we use the reference rate
\[
\lambda_\tau = \frac{\log{(1+\lambda \tau)}}{\tau}, \quad
\lambda =
\begin{cases}
p-1,&p<4/3\\
\frac{p}{4},& \frac{4}{3}\le p \le 4,p \ne 2,\\
1,&p>4.
\end{cases}
\]
We recall here that our current theory does not cover the cases $p<4/3$ or $p>4$; however, from Remark \ref{rem:powerlaw_lam_fem}, we expect to obtain the corresponding values for $\lambda$ asymptotically. %

\begin{figure}[hbt]
\centering
\includegraphics[width=0.85\textwidth]{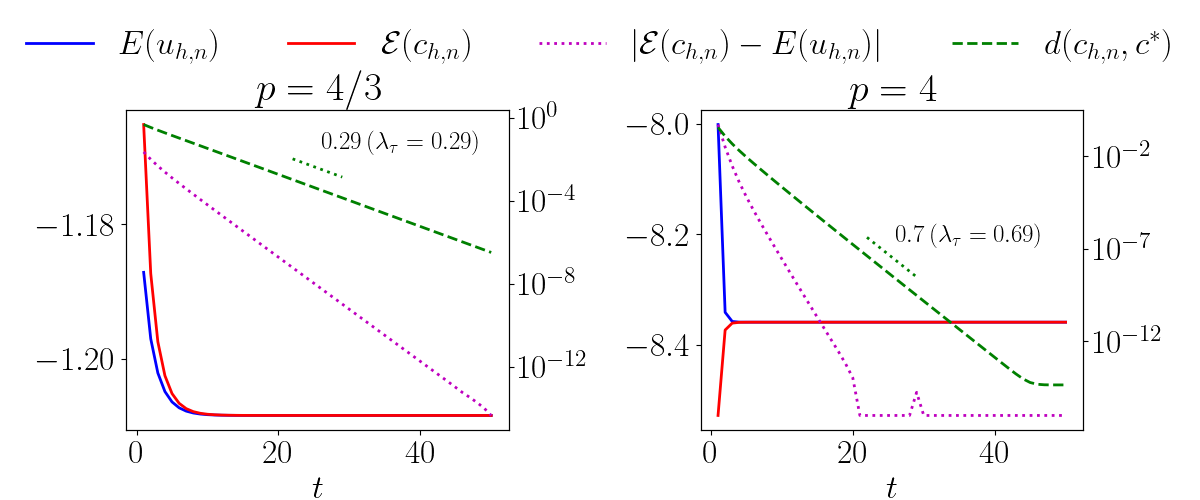}
\caption{Representative convergence curves for the $p$-Laplacian problem using the BE method with $\tau=1$ for different values of $p$ as shown on top. Values for $E(u_{h,n})$ (blue solid lines) and $\mathcal{E}(c_{h,n})$ (red solid lines) are given on the left $y$-axis, while $d(c_{h,n},c^*)$ (green dashed lines) and $|\mathcal{E}(c_{h,n}) -E(u_{h,n})|$ (magenta dotted lines) are given on a logarithmic scale on the right axis. Experimental convergence rates and reference rates (in parentheses) are annotated in each panel.}
\label{fig:tconv_repr_scalar}
\end{figure}
Figure \ref{fig:tconv_repr_scalar} shows representative convergence curves for the energies $E(u_{h,n})$ (blue solid lines, see eq. \eqref{eq:scalar_primal}) and $\mathcal{E}(c_{h,n})$ (red) as a function of time $t$ using the BE method with a constant time step $\tau=1$ for the cases $p=4/3$ (left panel) and $p=4$ (right panel) up to final time $T=50$.
As expected, $E(u_{h,n}) \le \mathcal{E}(c_{h,n})$ for $p=4/3$, while the inequality is reversed for $p=4$, see Equation \eqref{eq:Ee_bounded}. The time variation of the distance of the discrete conductivity from the steady state $d(c_{h,n},c^*)$ (green dashed lines) and the difference of the energies $|\mathcal{E}(c_{h,n}) -E(u_{h,n})|$ (magenta dotted lines) are reported on a logarithmic scale using the right $y$-axes, confirming convergence to the exact solution. The slopes of the experimental convergence curves in the metric are annotated in the panels, indicating perfect agreement with the reference rates (given in parentheses). %

\begin{figure}[hbt]
\centering
\includegraphics[width=0.95\textwidth]{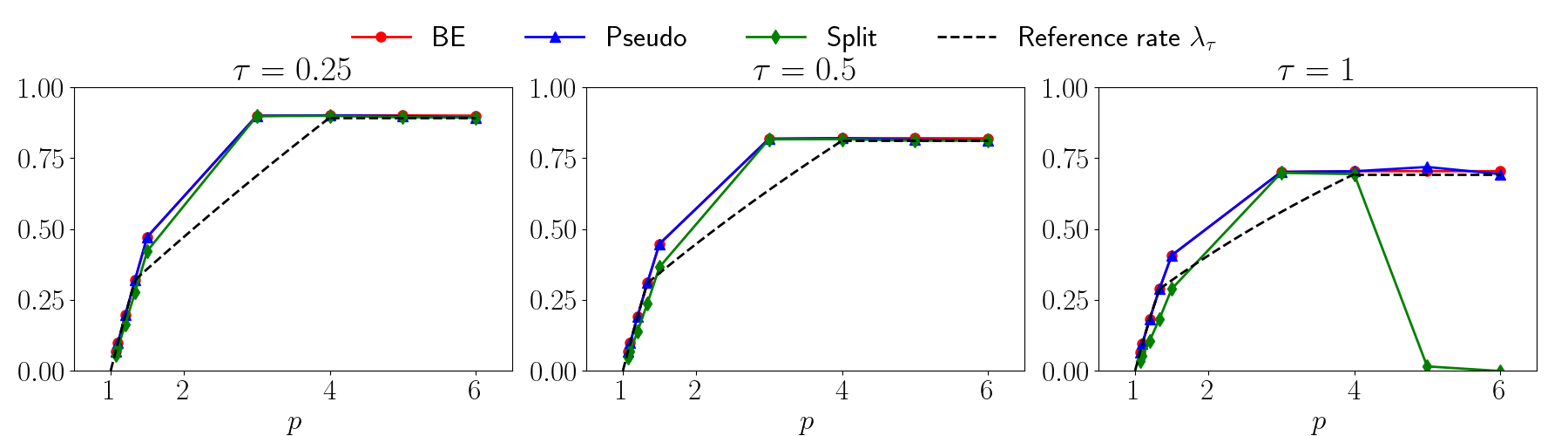}
\caption{Comparison of experimental convergence rates (solid lines) with reference rates (dashed) for the BE (red), Pseudo (blue), and Split (green) methods for the $p$-Laplacian problem as a function of $p$ for different values of $\tau$ as shown on top of each panel. For $p<4/3$ or $p>4$, the dashed curve represents the heuristic asymptotic rates from Remark \ref{rem:powerlaw_lam_fem}.}
\label{fig:tconv_lambda_scalar}
\end{figure}

Convergence rates using the BE, Pseudo, and Split methods are given for different values of $p$, ranging from $16/15$ to $6$, and for different constant time step values in Figure \ref{fig:tconv_lambda_scalar}. The solvers are run up to final time $T=50$. The Split method fails to converge for $\tau=1$ and $p=5,6$, while the observed convergence rates of the BE and Pseudo methods are very close to each other and are no smaller than the reference rates in these experiments.

\subsubsection{Power-law model: incompressible vector-valued case}\label{sec:powerlaw-vector-test}

\begin{figure}[hbt]
\centering
\includegraphics[width=0.85\textwidth]{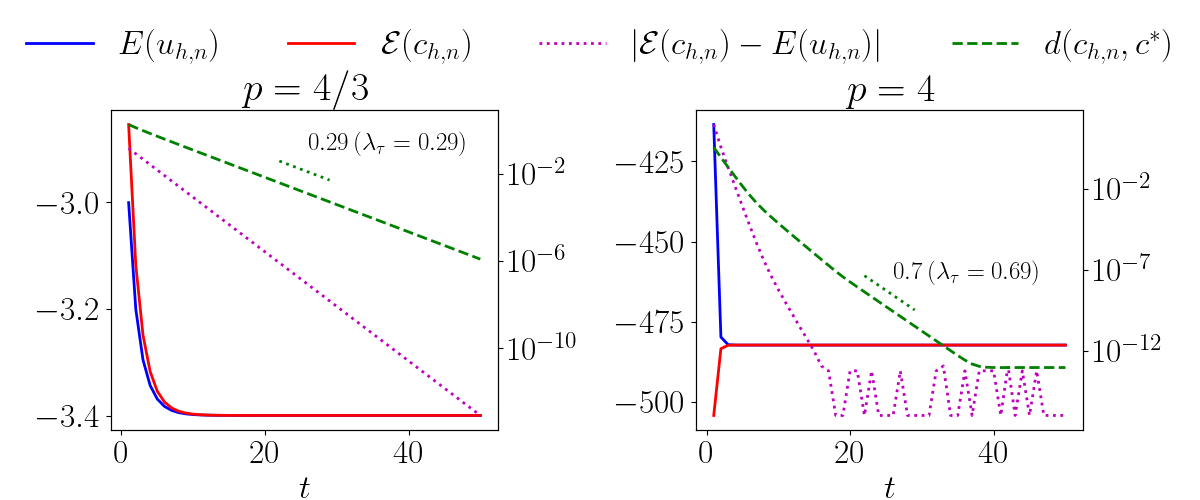}
\caption{Representative convergence curves for the $p$-Stokes problem using the BE method with $\tau=1$ for different values of $p$ as shown on top. Values for $E(u_{h,n})$ (blue solid lines) and $\mathcal{E}(c_{h,n})$ (red solid lines) are given on the left $y$-axis, while $d(c_{h,n},c^*)$ (green dashed lines) and $|\mathcal{E}(c_{h,n}) -E(u_{h,n})|$ (magenta dotted lines) are given on a logarithmic scale on the right axis. Experimental convergence rates and reference rates (in parentheses) are annotated in each panel.}
\label{fig:tconv_repr_vector}
\end{figure}

\begin{figure}[hbt]
\centering
\includegraphics[width=0.95\textwidth]{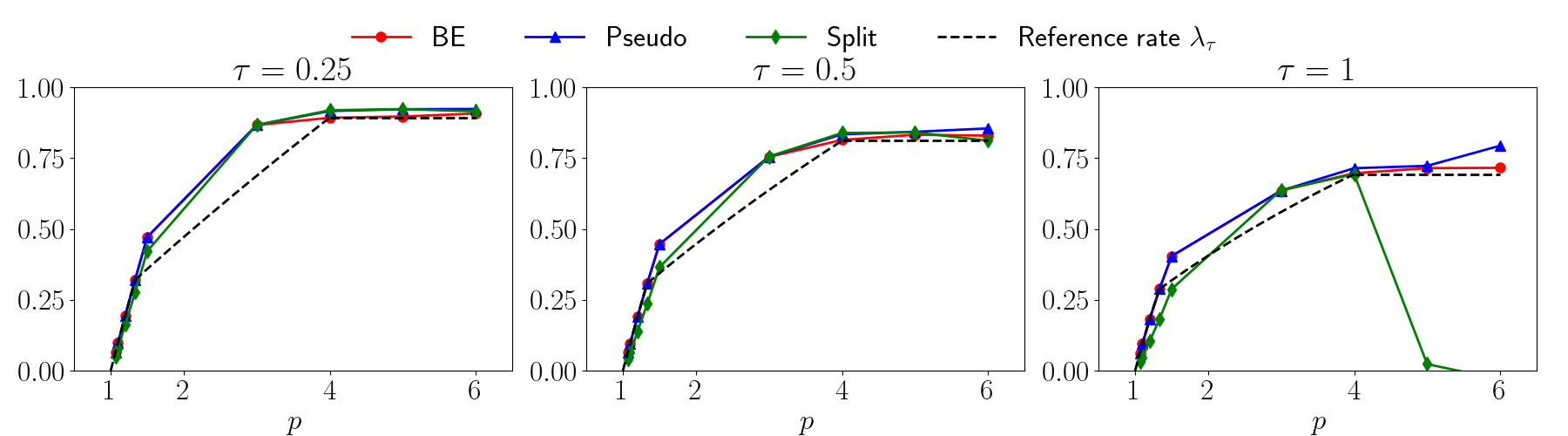}
\caption{Comparison of experimental convergence rates (solid lines) with reference rates (dashed) for the BE (red), Pseudo (blue), and Split (green) methods for the $p$-Stokes problem as a function of $p$ for different values of $\tau$ as shown on top of each panel. For $p<4/3$ or $p>4$, the dashed curve represents the heuristic asymptotic rates from Remark \ref{rem:powerlaw_lam_fem}.}
\label{fig:tconv_lambda_vector}
\end{figure}

For the vector-valued case, we instead consider the exact solution pair
\[
\bu_e(x) = 4(1-x_1^2 - x_2^2)(x_2, -x_1), \quad p_e(x) = x^2_1+ x^2_2 - \frac{1}{|\Omega|}\int_\Omega  x^2_1+ x^2_2 \dx,
\]
on $\Omega=B_{1,h}(0)$, a discretization of the unit disk with four thousand linear simplices, and we use $k=3$ for the Taylor--Hood finite element pair; essential boundary conditions are imposed by sampling $\bu_e(x)$ on the boundary points of the linear mesh. The magnitude of the symmetric gradient is $|\eps(\bu_e)|=4\sqrt{2}(x^2_1+ x^2_2)$ and $\diver(|\eps(\bu_e)|^{p-2}\eps(\bu_e)) \in \Lpp$ if and only if $p>\frac{1+\sqrt{17}}{4}$; as before, we identify the right-hand side as
\[
F(\bv) = \int_\Omega |\eps(\bu_e)|^{p-2}\eps(\bu_e) \fro \eps (\bv) \dx -\int_\Omega p_e \diver \bv\dx,
\]
and find that $|\eps(\bu_e)|^{p-2}\eps(\bu_e) \in C^{0,\alpha}(\Omega;\mathbb{R}^{d\times d})$ with $\alpha=\min\{2p-2,1\}$ and thus Assumption \ref{ass:F_Hm1_vector} is satisfied. Representative convergence curves for the cases $p=4/3$ and $p=4$ are given in Figure \ref{fig:tconv_repr_vector}, while experimental convergence rates are compared with the reference rates in Figure \ref{fig:tconv_lambda_vector}. We do not discuss the results in detail as they are qualitatively similar to those in the scalar case. The observed rates agree with the reference rates, and the Split method fails to converge for $\tau=1$ and $p=5,6$.

\subsubsection{Carreau--Yasuda model for incompressible fluids}

For the next experiment, we use the Carreau--Yasuda model \cite{yasuda1979investigation,yasuda1981shear} (see also
\cite[Section 4.1.a]{bird1986dynamics}), for which the viscosity law is
\[
\mu(r)=\mu_\infty+(\mu_0-\mu_\infty)\left(1+\rho^a r^{a/2}\right)^{\frac{n-1}{a}},
\quad
a>0,\quad \rho>0,\quad \mu_0>\mu_\infty\ge0,\quad n>0,\quad n\neq 1.
\]
Here $\mu_0$ is the zero-shear viscosity, $\rho$ sets the time scale, $a$ controls the transition width, and $n$ is the power-law index. The model is shear-thinning for $0<n<1$, with limiting viscosity $\mu_\infty$, and shear-thickening for $n>1$.

A direct calculation shows %
\[
\phi(t)
=
\frac{t^2}{2}
\left[
\mu_\infty
+
\Delta\mu\,
{}_2F_1
\left(
\frac{1-n}{a},
\frac{2}{a};
1+\frac{2}{a};
-\rho^a t^a
\right)
\right],\quad \Delta\mu:=\mu_0-\mu_\infty,
\]
where ${}_2F_1$ is the Gauss hypergeometric function. Moreover, setting $z:=\rho^a r^{a/2}$ and $q:=\frac{n-1}{a}$, we find
\[
\mu'(r)
=
\frac{\Delta\mu(n-1)}{2r}
z(1+z)^{q-1},
\quad
\mu''(r)
=
\frac{\Delta\mu(n-1)}{4r^2}
z(1+z)^{q-2}
\bigl[(a-2)+(n-3)z\bigr],
\]
which gives us, using $t^2:=r$ and \eqref{eq:mup_exact_final}
\[
\frac{t\phi''(t)}{\phi'(t)} = 1 + \frac{2r\mu'(r)}{\mu(r)}
=
1+
(n-1)
\frac{
\Delta\mu\,z(1+z)^{q-1}
}{
\mu_\infty+\Delta\mu(1+z)^q
}.
\]

Therefore, for $n<1$, $\Phi$ is strictly concave, while for $n>1$, $\Phi$ is strictly convex, and we have
\[
1+\frac{r\mu''(r)}{2\mu'(r)}
=
\frac{(a+2)+(n+1)z}{4(1+z)},\quad
\lambda_e:=\inf_{r>0}
\left\{
1+\frac{r\mu''(r)}{2\mu'(r)}
\right\}
=
\min\left\{
\frac{a+2}{4},
\frac{n+1}{4}
\right\}.
\]

We consider the case $\mu_\infty=0$ and, without loss of generality, simplify the model using $\rho=1$. From $\mu_\infty=0$, since $0<z/(1+z)<1$, we obtain the indices $\smphi=n>0$ and $\spphi=1$ in the concave case, and $\smphi=1$ and $\spphi=n<+\infty$ in the convex case; Assumption \ref{ass:sindices_strict} is thus not satisfied.

The structural conditions of Theorem \ref{thm:mielke_lambda_discrete} are satisfied for $a\le 2$, $1/3\le n\le3$, with $n\ne1$. Assuming that asymptotically $|D u_h[c_h(t)]|^2\le c_h(t)$ holds (see the discussion in Remark \ref{rem:powerlaw_lam_fem}), we find 
\begin{equation}\label{eq:lambda_carreau_yasuda}
\lambda =
\begin{cases}
\lambda_e,&a\le 2, 1/3\le n\le3,\\
1,&a>2, n>3,\\
\min\left\{n,\lambda_e\right\},&\text{otherwise}.
\end{cases}
\end{equation}
The requirements of Lemma \ref{lem:coefficient-residual} are not satisfied; nevertheless, the time-discrete iterations converge in the experiments reported below.
We use the same experimental setting as in Section \ref{sec:powerlaw-vector-test} in terms of manufactured solution, domain, and finite-element spaces. In this case, we use the right-hand side
\[
F(\bv) = \int_\Omega {\bf f}_e\cdot \bv\dx, \quad {\bf f}_e :=-\diver \left(\mu(|\eps(\bu_e)|^2) \eps(\bu_e)\right) + \grad p_e \in L^\infty(\Omega;\mathbb{R}^d).
\]

In Table \ref{tab:lambdaconv_carreau_yasuda}, we report the experimental convergence rates of the BE and Pseudo methods with constant time
steps $\tau=1$ measured using the square root of \eqref{eq:Dht_def}, and compare them with the reference rate $\lambda_\tau$,
with $\lambda$ given in \eqref{eq:lambda_carreau_yasuda}. 
In all cases, the observed rates for the BE and Pseudo methods are very close to each other and are no smaller than the reference rates. The reference rates capture the dependence on the model parameters and are sharpest near the limiting regimes. %
These tests examine convergence rates for varying model parameters, without investigating mesh dependence.

\begin{table}[hbt]
\centering
\caption{Experimental convergence rates for the BE (top left value) and Pseudo (top right value) methods compared against the reference rate $\lambda_\tau$ (bottom) for the Carreau--Yasuda model for different values of model parameters $n$ and $a$.}
\label{tab:lambdaconv_carreau_yasuda}
\begin{tabular}{lcccccc}
\toprule
 & $n=0.20$ & $n=0.33$ & $n=0.66$ & $n=0.90$ & $n=1.10$ & $n=2.00$\\
\midrule
\multirow{ 2}{*}{$a=1.00$} & 0.29, 0.29 & 0.37, 0.37 & 0.54, 0.54 & 0.65, 0.65 & 0.70, 0.70 & 0.67, 0.67\\
 & 0.18 & 0.29 & 0.35 & 0.39 & 0.42 & 0.56\\
\toprule
 & $a=0.25$ & $a=0.50$ & $a=1.00$ & $a=1.50$ & $a=2.00$ & $a=4.00$\\
\midrule
\multirow{ 2}{*}{$n=2.50$} & 0.56, 0.56 & 0.60, 0.60 & 0.67, 0.67 & 0.70, 0.70 & 0.70, 0.70 & 0.71, 0.71\\
 & 0.45 & 0.49 & 0.56 & 0.63 & 0.63 & 0.63\\
\toprule
 & $n=0.20$ & $n=0.33$ & $n=0.66$ & $n=1.10$ & $n=3.00$ & $n=4.00$\\
\midrule
\multirow{ 2}{*}{$a=3.00$} & 0.19, 0.19 & 0.29, 0.29 & 0.51, 0.51 & 0.74, 0.74 & 0.71, 0.71 & 0.71, 0.72\\
 & 0.18 & 0.29 & 0.35 & 0.42 & 0.69 & 0.69\\
\toprule
 & $a=0.25$ & $a=0.50$ & $a=1.00$ & $a=1.50$ & $a=2.00$ & $a=4.00$\\
\midrule
\multirow{ 2}{*}{$n=4.00$} & 0.58, 0.62 & 0.59, 0.69 & 0.67, 0.72 & 0.70, 0.72 & 0.70, 0.71 & 0.71, 0.73\\
 & 0.45 & 0.49 & 0.56 & 0.63 & 0.69 & 0.69\\
\bottomrule
\end{tabular}
\end{table}

\subsection{Comparison against Newton's method}

In this section, we compare the convergence of Newton's method for the solution of \eqref{eq:galerkin_fem} with the Pseudo iteration scheme endowed with adaptive time step selection as described in Section \ref{sec:time_adapt}. As a common stopping criterion, we use the residual of \eqref{eq:galerkin_fem}, i.e., the linear form
\[
\mathcal{R}^*[u_h](v_h):=\int_\Omega \mu(|D u_h|^2)D u_h:D v_h\dx - F(v_h),
\]
and declare convergence when the Euclidean norm satisfies $\|\mathcal{R}^*[u_h]\|_{\ell^2} < 10^{-8}$ within a maximum of 200 iterations. When reporting results, we use the symbol ``-'' to indicate failure to converge.
The initial time step $\tau_0$ is chosen empirically to obtain convergence with few nonlinear iterations.

\subsubsection{Power-law model}

The first test concerns the power-law model and the scalar benchmark problem considered in
\cite[Section 7.3]{balci2023relaxed}. We take a constant source term and impose
homogeneous essential boundary conditions on the L-shaped domain $\Omega=(-1,1)^2 \setminus [0,1)\times (-1,0]$.
Table \ref{tab:comp_newton_power_law} reports, for several values of $p$ ranging
from $1.5$ to $100$, and for successive refinement levels $l$, the number of
iterations required by the Pseudo method and by Newton's method, with the latter shown
in parentheses. The initial mesh consists of approximately $46\mathrm{K}$
triangular elements, and the polynomial degree is fixed at $k=1$. The table also reports the initial time step $\tau_0$; for small values of $p$, relatively large initial time steps can be employed,
whereas increasingly smaller values of $\tau_0$ are needed to obtain a
convergent iteration as $p$ grows. The iteration counts of the two methods are
comparable for $p \le 4$. At $p=8$, Newton's method fails on the two refined meshes; for $p\ge10$, it fails at all reported refinement levels. The Pseudo method remains convergent even for $p=100$. Moreover, the Pseudo method exhibits an essentially
mesh-independent and $p$-independent number of iterations for this test case.

\begin{table}[hbt]
\centering
\caption{Iteration counts for the Pseudo method compared against Newton's method (in parentheses) for the power-law model for different model parameters (columns) and refinement levels (rows). The initial time step $\tau_0$ is also reported.}
\label{tab:comp_newton_power_law}
\begin{tabular}{cccccccccc}
\toprule
 & $p = 1.5$ & $p = 3$ & $p = 4$ & $p = 8$ & $p = 10$ & $p = 20$ & $p = 40$ & $p = 80$ & $p = 100$\\
\midrule
$\tau_0$ & 10 & 9 & 4 & 0.5 & 0.2 & 0.1 & 0.05 & 0.02 & 0.015\\
\midrule
$l = 0$ & 15 (15) & 8 (8) & 9 (13) & 19 (25) & 28 (-) & 25 (-) & 25 (-) & 30 (-) & 31 (-)\\
$l = 1$ & 18 (16) & 8 (8) & 10 (15) & 18 (-) & 28 (-) & 25 (-) & 28 (-) & 33 (-) & 34 (-)\\
$l = 2$ & 18 (20) & 9 (10) & 12 (20) & 18 (-) & 31 (-) & 28 (-) & 29 (-) & 34 (-) & 38 (-)\\
\bottomrule
\end{tabular}

\end{table}

\begin{figure}[hbt]
\centering
\begin{tabular}{@{}cc@{}}
$u_h$ & $c_h$\\
\includegraphics[width=0.47\textwidth]{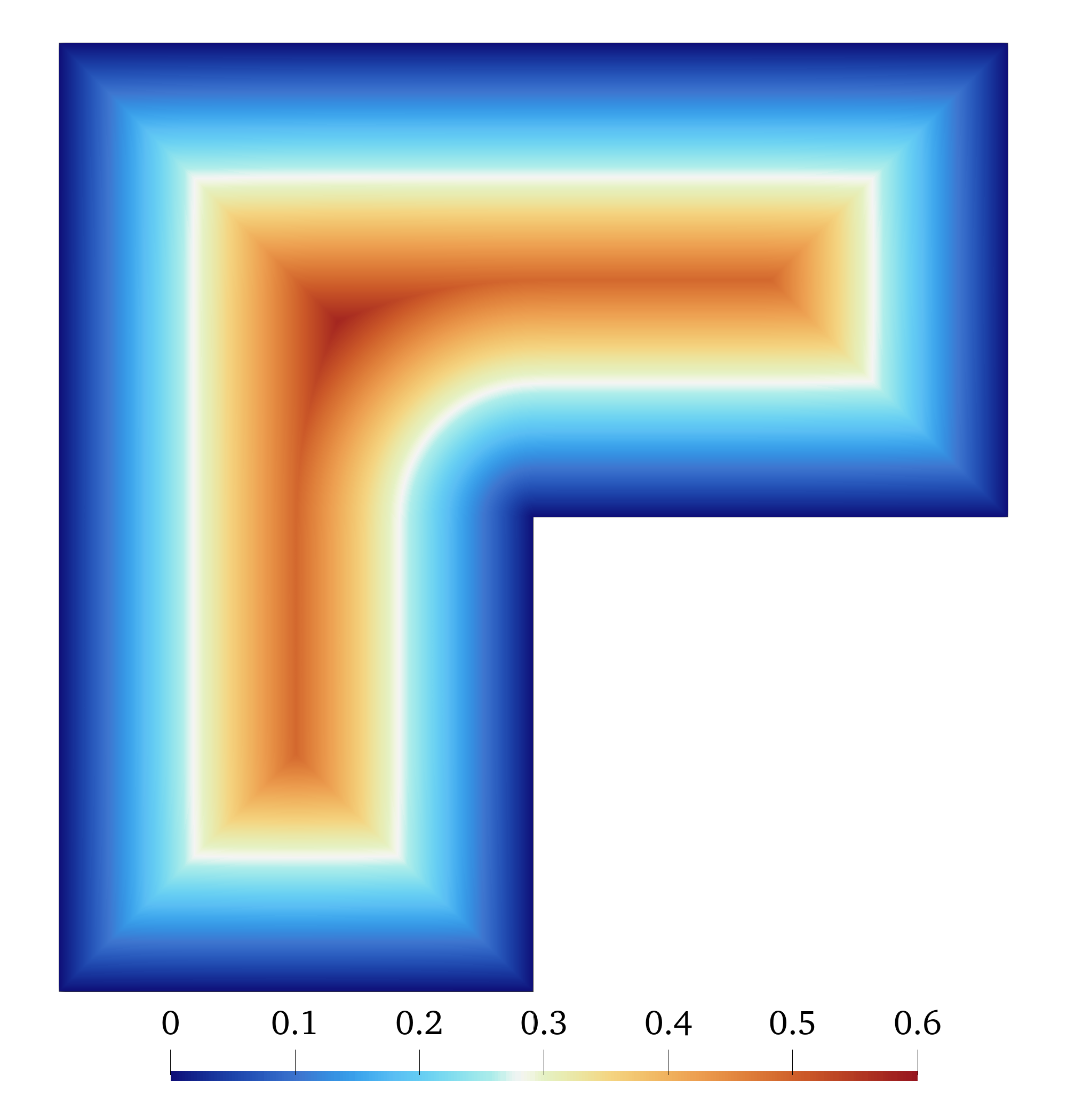} &
\includegraphics[width=0.47\textwidth]{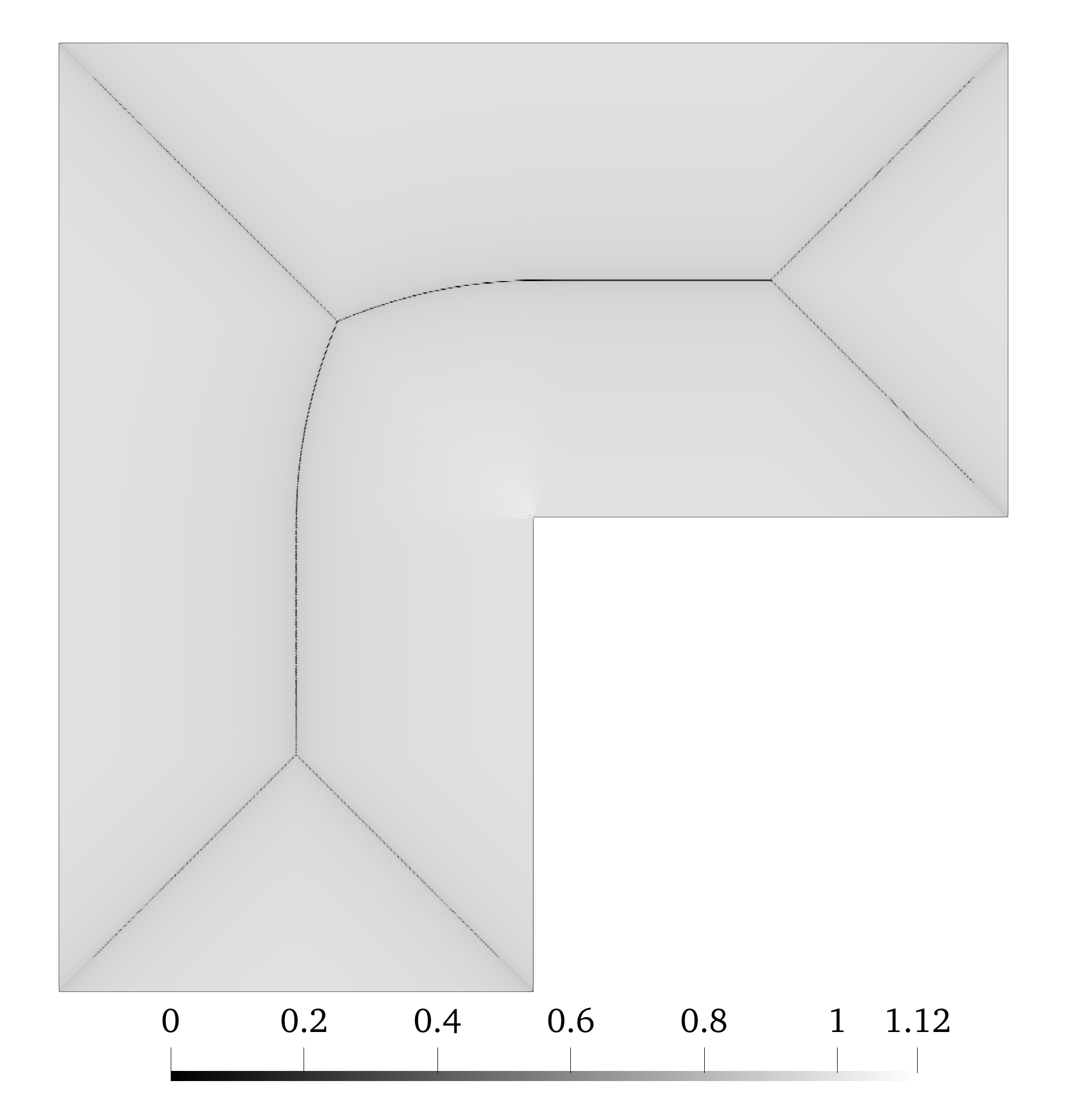}
\end{tabular}
\caption{Converged fields for the $p$-Laplacian benchmark on the L-shaped domain with $p=100$: $u_h$ (left) and $c_h$ (right).}
\label{fig:p100_solution}
\end{figure}

Figure~\ref{fig:p100_solution} shows the converged fields $u_h$ and $c_h$ for $p=100$.
The narrow bands in $c_h$ resemble the refinement pattern in \cite[Figure~4]{balci2023relaxed}, up to reflection of the L-shaped domain.

\subsubsection{Regularized Bingham model}

Here we consider the Bercovier--Engelman regularization \cite {BERCOVIER1980313} of Bingham fluids \cite{bingham1922fluidity} for which
\[
\phi(t) = \nu t^2+\sigma(\sqrt{t^2+b_\epsilon^2}-b_\epsilon), \quad
\mu(r) = 2\nu+\frac{\sigma}{\sqrt{r+b_\epsilon^2}},
\]
where $\sigma>0$ is the so-called yield stress, $\nu>0$ the fluid viscosity, and $b_\epsilon>0$ a regularization parameter. The model satisfies the hypotheses of Theorem \ref{thm:ode_bound} in the concave regime since
\[
\mu'(r) =-\frac{\sigma}{2}(r+b_\epsilon^2)^{-3/2}<0, \quad\forall r>0.
\]
However, it satisfies neither Assumption \ref{ass:sindices_strict} nor the requirements of Theorem \ref{thm:ode_fem_convergence} since
\[
\frac{t\phi''(t)}{\phi'(t)}
=\frac{2\nu+\sigma b_\epsilon^2(t^2+b_\epsilon^2)^{-3/2}}{
 2\nu+\sigma(t^2+b_\epsilon^2)^{-1/2}},
\]
and thus $\smphi>0$ but $\spphi=1$ (see also Remark \ref{rem:mu_bounded}). For this model,
\[
\mu''(r)=\frac{3\sigma}{4}(r+b_\epsilon^2)^{-5/2},
\quad
1+\frac{r\mu''(r)}{2\mu'(r)}
=
\frac{r+4b_\epsilon^2}{4(r+b_\epsilon^2)} .
\]
Hence Theorem~\ref{thm:mielke_lambda_discrete} is applicable with
$\lambda=1/4$ precisely when the structural condition in statement (ii) is satisfied. Such a condition reduces to $6\nu b_\epsilon \ge \sigma$, and it is not satisfied in the experiment we report.

We consider a manufactured Bingham problem for the incompressible vector-valued energy \eqref{eq:vector_min} on $\Omega=(0,1)^2$, using the velocity profile from \cite[Section 5.1]{heid2022adaptive}:
\[
\bu_e = (u_\infty,0), \quad
u_\infty(x) :=
\begin{cases}
\frac18\left(0.4^2 - (0.4 - 2x_2)^2\right), & 0 \le x_2 \le 0.2, \\
0.02, & 0.2 < x_2 < 0.8, \\
\frac18\left(0.4^2 - (2x_2 - 1.6)^2\right), & 0.8 \le x_2 \le 1,
\end{cases}
\]
with pressure $p_e=0$, boundary velocity $\bu_e$, and the manufactured load
\[
F(\bv)=\int_\Omega\mu(|\eps(\bu_e)|^2)\,\eps(\bu_e):\eps(\bv)\dx.
\]
Here $\mu$ is the regularized viscosity, so the load depends on $b_\epsilon$ while the exact velocity remains fixed. The experiment is performed with $\sigma=0.3$ and $\nu=1$, for
regularization parameters $b_\epsilon \in \{10^{-2},10^{-3},10^{-4}\}$,
on a sequence of uniformly refined triangular meshes starting from an
initial mesh with approximately $8\,\mathrm{K}$ elements. We use the Taylor--Hood pair with $k=2$ to discretize the problem in space, and the initial time step of the Pseudo method is set to $\tau_0=1$. The iteration counts are reported in Table \ref{tab:comp_newton_bingham} (left). The Pseudo method exhibits mesh-independent
convergence and remains robust with respect to the regularization
parameter. Moreover, it requires a comparable or, in several cases,
smaller number of iterations than Newton's method.
\begin{table}[hbt]
\centering
\caption{Iteration counts for the Pseudo method compared with Newton's method (in parentheses): Bingham (left) and optimal design (right), for the parameters shown and refinement levels $l$.}
\label{tab:comp_newton_bingham}
\label{tab:comp_newton_optimal_design}
\begin{tabular}{@{}c@{\qquad}c@{}}
\begin{tabular}[t]{cccc}
\toprule
 & $b_\epsilon = 10^{-2}$ & $b_\epsilon = 10^{-3}$ & $b_\epsilon = 10^{-4}$\\
\midrule
$l = 0$ & 13 (8) & 19 (27) & 40 (-)\\
$l = 1$ & 12 (8) & 20 (18) & 35 (86)\\
$l = 2$ & 12 (7) & 19 (22) & 37 (92)\\
\bottomrule
\end{tabular}
&
\begin{tabular}[t]{ccc}
\toprule
 & $\lambda_d$ = 0.0084 & $\lambda_d$ = 0.0145\\
\midrule
$l = 0$ & 25 (46) & 30 (56)\\
$l = 1$ & 37 (91) & 41 (155)\\
$l = 2$ & 50 (-) & 58 (-)\\
\bottomrule
\end{tabular}
\end{tabular}
\end{table}

\subsubsection{Optimal design}

In the last experiment, we consider the optimal design problem described in \cite{bartels2008convergent} where
\[
    \phi(r)
    =
    \begin{cases}
    \mu_2 r^2/2,
    & 0\leq r\leq \xi_1,
    \\
    \xi_1\mu_2\left(r-\xi_1/2\right),
    & \xi_1\leq r\leq \xi_2,
    \\
    \mu_1r^2/2
    -
    \xi_1\mu_2\left(\xi_1/2-\xi_2/2\right),
    & \xi_2\leq r,
    \end{cases}
\]
with %
$
    0<\xi_1<\xi_2,
    0<\mu_1<\mu_2,
    \xi_1\mu_2=\xi_2\mu_1.
$
Then
\[
\phi'(r)=
\begin{cases}
\mu_2r,& 0<r<\xi_1,\\
\xi_1\mu_2,& \xi_1<r<\xi_2,\\
\mu_1r,& r>\xi_2,
\end{cases}
,\quad
\mu(r)=
\begin{cases}
\mu_2,& 0<r<\xi_1^2,\\
\xi_1\mu_2 r^{-1/2},& \xi_1^2<r<\xi_2^2,\\
\mu_1,& r>\xi_2^2.
\end{cases}
\]
The potential is $C^1$. On its three smooth branches, $\phi''$ equals $\mu_2,0,\mu_1$, respectively, giving $\smphi=0$ and $\spphi=1$. Moreover, $\mu'$ vanishes for $0<r<\xi_1^2$ and $r>\xi_2^2$, while $\mu'(r)=-\xi_1\mu_2/(2r^{3/2})$ for $\xi_1^2<r<\xi_2^2$. Hence the Riesz map is degenerate.
This model is an N-function, but falls outside our convergence theory. We reproduce two experimental settings from \cite[Sections 5.3.2--5.3.3]{carstensen2021unstabilized}
using the scalar model with $k=1$, with $\mu_1=1$, $\mu_2=2$, and $\xi_1 = \sqrt{2\lambda_d\mu_1/\mu_2}$ for a given $\lambda_d>0$ and constant right-hand side; the first experiment uses $\Omega=(0,1)^2$ and $\lambda_d=0.0084$, while the second experiment uses an L-shaped domain $\Omega=(-1,1)^2 \setminus [0,1)\times (-1,0]$ and $\lambda_d = 0.0145$. Iteration counts using the Pseudo method with adaptive time stepping with $\tau_0=1$ and Newton's method are reported in Table \ref{tab:comp_newton_optimal_design} (right), and confirm the robustness of the Pseudo method. In both tests, the Pseudo method converges, but its iteration count increases with mesh refinement.

\section{Conclusions}\label{sec:conclusions}

In this work, we introduced an auxiliary gradient-flow framework for generalized Newtonian-type variational problems. At the continuous level, we studied the auxiliary energy in a metric space adapted to the growth of the N-function $\phi$. Under suitable assumptions, we proved lower semicontinuity, geodesical $\lambda$-convexity, and exponential convergence of the corresponding minimizing-movement scheme. At the discrete level, we derived a finite-dimensional metric gradient flow using an explicit Riesz map to identify first variations with metric gradients. For the resulting evolution \eqref{eq:ode_fem}, we proved global existence and uniqueness of the solution, and convergence to the finite element solution under additional structural assumptions.

For power laws, the convergence theory covers $4/3\le p\le4$, $p\ne2$. Outside this range, the discrete formulation suggests asymptotic convergence rate estimates near equilibrium that are confirmed by numerical experiments.

The time discretizations include a splitting scheme related to \kac iteration and an adaptive iteration inspired by pseudo-transient continuation and compatible with scalable solvers. Comparisons with Newton's method show competitive performance, including convergence in cases where Newton fails. We observe essentially mesh-independent iteration counts for the Pseudo method in the power-law and Bingham tests, whereas the counts increase with mesh refinement in the optimal-design tests.

Future work includes refining the ODE theory, analyzing fully discrete convergence, adaptive mesh refinement, and extensions to less regular models, including dependence on $\eps(\bu)$ \cite{ruuvzivcka2007non}. Adaptive time stepping could also benefit from extending \cite[Section 6.4]{deuflhard} and \cite{amrein2017adaptive} to differential-algebraic systems.

\section*{Acknowledgments}
SZ gratefully acknowledges Giuseppe Savar\'e for fruitful discussions. DB is a member of the INdAM research group GNCS. The work of SZ, DB, and GD was partially supported by the KAUST CRG13-2025 grant (6911).

\bibliographystyle{plain}
\bibliography{bibliography}

\clearpage

\appendix
\renewcommand{\theHsection}{\Alph{section}}

\section{Appendix}\label{sec:appendix}

\begin{lemma}\label{lem:phiprime_bounds_proof}
Let $0<s\le t$ and $\smphi, \spphi$ be given in Assumption \ref{ass:sindices}. Then
\begin{enumerate}
\item[(i)] $\phi'$ satisfies the two-point growth:
$\displaystyle \left(\frac{s}{t}\right)^{\spphi}\le \frac{\phi'(s)}{\phi'(t)}\le \left(\frac{s}{t}\right)^{\smphi},$
\item[(ii)] $\mu$ satisfies the two-point growth:
$\displaystyle \left(\frac{s}{t}\right)^{\frac{\spphi-1}{2}}\le \frac{\mu(s)}{\mu(t)}\le \left(\frac{s}{t}\right)^{\frac{\smphi-1}{2}}.$
\end{enumerate}
\end{lemma}
\mybegproof
Let $g(r):=\log(\phi'(r))$, $g'(r)=\phi''(r)/\phi'(r)$; from Assumption \ref{ass:sindices} we have $\smphi/r\le g'(r)\le \spphi/r$.
Integrating these inequalities from $s$ to $t$ with $0<s\le t$ yields
\[
\smphi(\log(t) - \log(s))\le g(t)-g(s)\le \spphi(\log(t) - \log(s)).
\]
Exponentiating gives 
\[
\left(\frac{t}{s}\right)^{\smphi} \le \frac{\phi'(t)}{\phi'(s)}\le \left(\frac{t}{s}\right)^{\spphi},
\]
which is equivalent to the statement (i) for $\phi'$.

In order to prove the statement for $\mu$, we replace $s$ and $t$ with their square roots and use the definition of $\mu$ from \eqref{eq:Phi_mu}, i.e.
\[
\frac{\mu(s)}{\mu(t)} = \frac{\phi'(\sqrt{s})}{\phi'(\sqrt{t})}\,\frac{\sqrt{t}}{\sqrt{s}},
\]
then apply the bounds for $\phi'$ with $\sqrt{s}\le \sqrt{t}$.

\myendproof

\begin{lemma}\label{lem:integ_E_c_proof}
Let $\smphi, \spphi$ be given as in Assumption \ref{ass:sindices}. Then, for all $r\ge0$,
\[
\frac{1}{\spphi+1}\,r\phi'(r) \le \phi(r) \le \frac{1}{\smphi+1}\,r\phi'(r),\quad \frac{1}{\spphi+1}\,r\mu(r) \le \Phi(r) \le \frac{1}{\smphi+1}\,r\mu(r),
\]
where $r\phi'(r)$ and $r\mu(r)$ at $r=0$ are understood through their continuous extensions $\lim_{s\to 0^+} s\phi'(s)=0$ and $\lim_{s\to 0^+} s\mu(s)=0$.
\end{lemma}
\mybegproof
Let $r>0$ and set $t=\sqrt r$. Then, since $r\mu(r)=r\frac{\phi'(t)}{t}=t\phi'(t)$, it suffices to prove the statement
\[
\frac{1}{\spphi+1}\,t\phi'(t) \le \phi(t) \le \frac{1}{\smphi+1}\,t\phi'(t).
\]
Since $\phi'(t) > 0$, from (i) in Lemma \ref{lem:phiprime_bounds},
\[
\phi'(t)\left(\frac{s}{t}\right)^{\spphi} \le \phi'(s)\le \phi'(t)\left(\frac{s}{t}\right)^{\smphi},\quad 0<s\le t.
\]
Integrating from $0$ to $t$ yields
\[
\frac{1}{\spphi+1}\,t\phi'(t) = \int_0^t \phi'(t)\left(\frac{s}{t}\right)^{\spphi}\ds
 \le \phi(t)\le \int_0^t \phi'(t)\left(\frac{s}{t}\right)^{\smphi}\ds
=\frac{1}{\smphi+1}\,t\phi'(t).
\]

It remains to consider the case $r=0$. Although $\phi'(r)$ may only be defined for $r>0$, the quantity $r\phi'(r)$ admits a continuous extension at $r=0$, since
\[
0 \le r\phi'(r) \le \sipphi\phi(r) \to 0 \quad \text{as } r\to 0^+.
\]
where $\sipphi$ is given in Assumption \ref{ass:phi_assum_reflex}. Similar reasoning holds for $r\mu(r)$.
\myendproof

\begin{lemma}\label{lem:psi_bounds_mu_proof}
Let $\psi$ be defined as in \eqref{eq:psi_def} and let Assumption \ref{ass:sindices_strict} hold. Then
\begin{enumerate}
\item[(i)] There exist constants $0<C_1< C_2<+\infty$ depending only on $\smphi, \spphi$ such that for all $r\ge0$,
\[
C_1\mu(r)r \le \psi(r)^2 \le C_2\mu(r)r.
\]
\item[(ii)] There exist constants $0<C_3< C_4<+\infty$ depending only on $\smphi, \spphi$ such that for all $r\ge0$,
\[
C_3\Phi(r) \le \psi(r)^2 \le C_4\Phi(r).
\]
\end{enumerate}
\end{lemma}
\mybegproof
We prove statement (i); statement (ii) can be obtained by combining the chain of inequalities of (i) with Lemma \ref{lem:integ_E_c_proof}. We define for $t>0$
\[
\theta(t):=\omega\left(\frac{t\phi''(t)}{\phi'(t)}-1\right).
\]
From Assumption \ref{ass:sindices_strict}, we have
\begin{equation}\label{eq:theta_bounds}
0 < A \leq \theta(t) \leq B,
\end{equation}
with $A = \smphi-1$ and $B = \spphi-1$ when $\Phi$ is convex, while  $A = 1-\spphi$ and $B = 1-\smphi$ when $\Phi$ is concave.
From \eqref{eq:mup_exact_final} we obtain
\[
\psi'(s)=\sqrt{\frac{\omega}{2}\mu'(s)}=\sqrt{\frac{\theta(\sqrt s)}{4}}\sqrt{\frac{\mu(s)}{s}},
\]
and thus, integrating from $0$ to $r$, and using \eqref{eq:theta_bounds} we obtain
\begin{equation}\label{eq:psi_temp_bounds}
\frac{\sqrt{A}}{2} \int^r_0 \sqrt{\frac{\mu(s)}{s}}\ds \le \psi(r) \le \frac{\sqrt{B}}{2} \int^r_0 \sqrt{\frac{\mu(s)}{s}}\ds.
\end{equation}

Using the upper bound for $\mu$ from statement (ii) in Lemma \ref{lem:phiprime_bounds} we have
\[
\sqrt{\frac{\mu(s)}{s}}
\le \sqrt{\mu(r)}\,r^{\frac{1-\smphi}{4}}\, s^{\frac{\smphi-3}{4}}, \quad 0<s\le r.
\]
Inserting the above inequality into the upper bound in \eqref{eq:psi_temp_bounds}, and noting that $(\smphi-3)/4>-1$, yields
\[
\psi(r)\le \frac{\sqrt{B}}{2}\sqrt{\mu(r)}\,r^{\frac{1-\smphi}{4}}\int_0^r s^{\frac{\smphi-3}{4}}\ds = \frac{2\sqrt{B}}{ \smphi+1}\sqrt{\mu(r)}\,r^{\frac{1-\smphi}{4}} r^{\frac{\smphi+1}{4}} = \frac{2\sqrt{B}}{ \smphi+1}\sqrt{\mu(r)r}.
\]
Squaring yields 
\[
\psi(r)^2\le C_2 r\mu(r),\quad C_2 :=
\begin{cases}
4\frac{\spphi-1}{(\smphi+1)^2},&\Phi \text{ convex},\\
4\frac{1-\smphi}{(\smphi+1)^2},&\Phi \text{ concave}.
\end{cases}
\]

For the proof of the lower bound, we proceed similarly, using the lower bound for $\mu$ in (ii) of Lemma \ref{lem:phiprime_bounds}
\[
\sqrt{\frac{\mu(s)}{s}}
\ge \sqrt{\mu(r)}\,r^{\frac{1-\spphi}{4}}\, s^{\frac{\spphi-3}{4}}, \quad 0<s\le r,
\]
which yields
\[
\psi(r)\ge \frac{\sqrt{A}}{2}\sqrt{\mu(r)}\,r^{\frac{1-\spphi}{4}}\int_0^r s^{\frac{\spphi-3}{4}}\ds = \frac{2\sqrt{A}}{ \spphi+1}\sqrt{\mu(r)}\,r^{\frac{1-\spphi}{4}} r^{\frac{\spphi+1}{4}} = \frac{2\sqrt{A}}{ \spphi+1}\sqrt{\mu(r)r}.
\]
Squaring yields
\[
\psi(r)^2\ge C_1 r\mu(r),\quad C_1 :=
\begin{cases}
4\frac{\smphi-1}{(\spphi+1)^2},&\Phi \text{ convex},\\
4\frac{1-\spphi}{(\spphi+1)^2},&\Phi \text{ concave}.
\end{cases}
\]
\myendproof

\begin{lemma}\label{lem:C_d_complete_geodesic_proof}
Let Assumption \ref{ass:sindices_strict} hold. Then $\LPhi$, defined as in \eqref{eq:LPhi_def}, and endowed with the metric $d$ defined in \eqref{eq:X_distance}, is a complete geodesic space.
\end{lemma}

\mybegproof
Let $\{c_n\}_{n \in \mathbb{N}} \subset \LPhi$ be a Cauchy sequence with respect to the metric $d$. By the definition of the distance in \eqref{eq:X_distance}, we have
\begin{equation*}
    d(c_n, c_m) = \|\Psi(c_n) - \Psi(c_m)\|_\Ltwo.
\end{equation*}
This implies that the sequence $\eta_n = \Psi(c_n)$ is a Cauchy sequence in $\Ltwo$,  which is a complete metric space; thus there exists a limit function $\eta^* \in \Ltwo$ such that $\eta_n \to \eta^*$ in $\Ltwo$ and $\eta_n$ admits a subsequence (not relabeled) for which $\eta_n\to \eta^*$ almost everywhere in $\Omega$ as $n \to +\infty$.

We define the candidate limit
\begin{equation*}
    c^*(x) := \psi^{-1}(\eta^*(x)) \quad \text{a.e. in } \Omega,
\end{equation*}
which is well defined since $\psi: \mathbb{R} \to \mathbb{R}$ is a homeomorphism. The identity $\Psi(c^*) = \eta^* \in \Ltwo$ ensures that $c^* \in \LPhi$ according to Lemma \ref{lem:psi_bounds_mu_proof}. 
Finally, by the definition of the metric, we have
\begin{equation*}
    \lim_{n \to +\infty} d(c_n, c^*) = \lim_{n \to +\infty} \|\Psi(c_n) - \Psi(c^*)\|_{\Ltwo} = \lim_{n \to +\infty} \|\eta_n - \eta^*\|_{\Ltwo} = 0,
\end{equation*}
and thus, the Cauchy sequence $\{c_n\}$ converges to $c^* \in \LPhi$, and the space is complete.

To prove that $\gamma_g(s)$ defined in \eqref{eq:geodesics} is in fact a geodesic, a simple substitution shows
\[
\begin{aligned}
d(\gamma_g(s),\gamma_g(r))^2
&= \int_\Omega |\Psi(\gamma_g(s)) - \Psi(\gamma_g(r))|^2\dx \\
&= \int_\Omega |(1-s)\Psi(c_0) + s\,\Psi(c_1) - (1-r)\Psi(c_0) - r\,\Psi(c_1)|^2\dx \\
&= \int_\Omega |(r-s)\Psi(c_0) -(r-s)\Psi(c_1)|^2\dx \\
&= |r-s|^2\int_\Omega |\Psi(c_0) -\Psi(c_1)|^2\dx. \\
\end{aligned}
\]
\myendproof

\begin{theorem}\label{thm:Jc_cont_proof}
Let Assumption \ref{ass:sindices_strict} hold. Then, the functional $\Jc$  defined in \eqref{eq:Jc_def} is continuous on $\LPhip$ in the topology induced by the metric \eqref{eq:X_distance}.
\end{theorem}
\mybegproof
Set $h(r) := \Phi(r) -\frac{1}{2}\mu(r)r$. Since $2\Phi'(r)=\mu(r)$, we have $h'(r) = - \frac{1}{2}\mu'(r) r$, and hence
\[
h(r) = - \frac{1}{2}\int^r_0 \mu'(s)s\ds.
\]
Combining \eqref{eq:mup_exact_final} with Assumption \ref{ass:sindices} yields
\begin{equation}\label{eq:temp_muprime_bound}
A \mu(s) \le \omega s \mu'(s) \le B \mu(s),
\end{equation}
with $A = \frac{\smphi-1}{2}$ and $B = \frac{\spphi-1}{2}$ when $\Phi$ is convex, while  $A = \frac{1-\spphi}{2}$ and $B = \frac{1-\smphi}{2}$ when $\Phi$ is concave.
Integrating both sides from zero to $r$, we obtain
$A \Phi(r) \le -\omega h(r) \le B \Phi(r)$,
and thus $-\omega \Jc(c)$ is positive and controlled by the modular of $c$ in $\LPhip$.

Defining $q(\eta) := -\omega h(\psi^{-1}(\eta))$ with $\eta := \psi(r)$ and $\psi$ given by \eqref{eq:psi_def}, straightforward calculations yield
\begin{equation}\label{eq:k_eta_prime_temp}
\begin{aligned}
q'(\eta) = -\omega \frac{h'(r)}{\psi'(r)} =  r\frac{\omega \mu'(r)/2}{\sqrt{\omega \mu'(r)/2}} = \sqrt{\frac{r}{2}}\sqrt{\omega \mu'(r)r}  \leq \sqrt{\frac{B}{2}}\sqrt{\mu(r) r} \leq L \psi(r) = L\eta,
\end{aligned}
\end{equation}
where the first inequality follows from the upper bound in \eqref{eq:temp_muprime_bound} while the second inequality follows from the lower bound in statement (i) of Lemma \ref{lem:psi_bounds_mu}. Thus the constant is $L=\sqrt{\frac{B}{2C_1}}$.
From \eqref{eq:k_eta_prime_temp}, we obtain
\begin{equation}\label{eq:k_eta_prime_temp_2}
\left|q(\eta_2)-q(\eta_1)\right| = \left|\int^{\eta_2}_{\eta_1}q'(s)\ds\right|\leq \frac{L}{2}\left|\eta^2_2 - \eta^2_1\right| =  \frac{L}{2}\left|\eta_2 - \eta_1\right|\left|\eta_2 + \eta_1\right|,\quad \forall \eta_1,\eta_2 \ge 0.
\end{equation}
Let $\{c_n\}_{n\in\mathbb{N}}\subset \LPhip$ be such that $c_n\to c^*$ in the  $d$-topology; then, since $\LPhip$ is closed, $c^\ast \in \LPhip$. Using the notation $\eta_n := \Psi(c_n)$, $\eta^* := \Psi(c^*)$, we have, using \eqref{eq:k_eta_prime_temp_2}, the Cauchy--Schwarz inequality, the triangle inequality, and the definition of the metric \eqref{eq:X_distance}
\[
\begin{aligned}
|\Jc(c_n) - \Jc(c^*)|
&\leq \frac{L}{2} \int_\Omega \left|\eta_n - \eta^*\right| \left|\eta_n + \eta^*\right| \dx\\
&\leq \frac{L}{2} \|\eta_n +\eta^*\|_\Ltwo \|\eta_n -\eta^*\|_\Ltwo  \leq \frac{L}{2} (d(c_n,c^*) + 2\|\eta^*\|_\Ltwo) d(c_n,c^*) \to 0,
\end{aligned}
\]
and thus $\Jc$ is continuous on $\LPhip$.
\myendproof

\begin{theorem}\label{thm:Ju_cont_proof}
Let Assumption  \ref{ass:sindices_strict} hold. Under Assumptions \ref{ass:F_Hm1} and \ref{ass:F_Hm1_vector}, the functionals $\Ju$  defined in \eqref{eq:Ju_def} and $\Jbu$  defined in \eqref{eq:Ju_def_vector} are continuous on $\Ce$ in the $d$-topology for all $\epsilon \in (0,1)$.
\end{theorem}
\mybegproof
Let $\{c_n\}_{n\in\mathbb{N}}\subset \Ce$ be such that $c_n\to c^*$ in the  $d$-topology; then, since $\Ce$ is closed in the $d$-topology, $c^\ast \in \Ce$. By Chebyshev's inequality, $\Psi(c_n) \to \Psi(c^\ast)$ in measure on $\Omega$. Since $\Psi^{-1}$ is continuous, we also have $c_n \to c^\ast$ in measure, as well as $\mu(c_n) \to \mu(c^\ast)$ since $\mu$ is continuous and finite on $[\epsilon,1/\epsilon]$ by \eqref{eq:mu_unif_bounds}.

We first give the proof for the scalar case, and show that $u[c_n]$ converges to $u[c^\ast]$ in $\Honez$; this allows us to prove the statement since we obtain $F(u[c_n])\to F(u[c^\ast])$ given that $F\in\Honez^\ast$ and, using the identity in \eqref{eq:energy_u_c}, $\Ju(c_n)\to  \Ju(c^\ast)$. For convenience, we write $a_n:=\mu(c_n)$, $ a^\ast:=\mu(c^\ast)$, $u_n:=u[c_n]$,  $u^\ast:=u[c^\ast]$.
Then, from \eqref{eq:scalar_state}
\[
\int_\Omega a_n\grad u_n\cdot\grad v\dx=F(v)=\int_\Omega a^\ast\grad u^\ast\cdot\grad v \dx,  \quad\forall v\in\Honez.
\]
Subtracting $\int_\Omega a_n\grad u^\ast\cdot\grad v\dx$ from both sides yields
\[
\int_\Omega a_n\grad(u_n-u^\ast)\cdot\grad v\dx=\int_\Omega (a^\ast-a_n)\grad u^\ast\cdot\grad v\dx.
\]
Choosing $v=u_n-u^\ast$ we obtain
\[
\int_\Omega a_n|{\grad(u_n-u^\ast)}|^2\dx=\int_\Omega (a^\ast-a_n)\grad u^\ast\cdot\grad(u_n-u^\ast)\dx.
\]
Using the lower bound $a_n\ge \mu_-$ from \eqref{eq:mu_unif_bounds}, and a Cauchy--Schwarz inequality argument we get
\begin{equation*}%
\|\grad(u_n-u^\ast)\|_\Ltwo\le\frac{1}{\mu_-}\|(a_n-a^\ast)\grad u^\ast\|_\Ltwo.
\end{equation*}
The right-hand side of %
the above inequality is uniformly integrable since $u^\ast\in \Honez$ and
\[
|(a_n-a^\ast)\grad u^\ast|^2\le 4\mu_+^2|\grad u^\ast|^2.
\]
Therefore, by the Lebesgue--Vitali theorem \cite[Theorem 4.5.4]{bogachev}, $\|\grad(u_n-u^\ast)\|_\Ltwo\to 0$  since $a_n\to a^\ast$ in measure. By Poincar\'e's inequality, convergence of the gradients implies $u[c_n]\to u[c^\ast]$ in $\Honez$.

The proof for the incompressible vector-valued case is very similar. Using the shorthand notation
$\bu_n:=\bu[c_n]$, $\bu^\ast:=\bu[c^\ast]$,
$p_n:=p[c_n]$,  and $p^\ast:=p[c^\ast]$,
from \eqref{eq:vector_state} and \eqref{eq:vector_incomp} we obtain
\[
\begin{aligned}
\int_\Omega a_n\eps(\bu_n-\bu^\ast)\fro\eps(\bv)\dx -\int_\Omega\diver\bv\,(p_n - p^\ast)\dx&=\int_\Omega (a^\ast-a_n)\,\eps (\bu^\ast)\fro\eps (\bv)\dx,\quad \forall\bv\in\vHonez,\\
\int_\Omega\diver(\bu_n - \bu^\ast)\,q\dx &= 0, \quad \forall q\in Q.
\end{aligned}
\]
Choosing the divergence-free test function $\bv = \bu_n - \bu^\ast$ yields
\[
\int_\Omega a_n|{\eps(\bu_n-\bu^\ast)}|^2\dx=\int_\Omega (a^\ast-a_n)\eps(\bu^\ast)\fro\eps(\bu_n-\bu^\ast)\dx.
\]
Repeating the same arguments, using the symmetric gradient instead of the plain gradient, we arrive at $\|\eps(\bu_n-\bu^\ast)\|_\Ltwo\to 0$ and thus, by Korn's and Poincar\'e's inequalities, $\bu[c_n]\to \bu[c^\ast]$ in $\vHonez$.
\myendproof

\begin{theorem}\label{thm:Ee_lsc_proof}
Under the assumptions of Theorems \ref{thm:Jc_cont_proof} and \ref{thm:Ju_cont_proof}, the functional $\Ee$  defined in \eqref{eq:E_c} for the scalar case and  \eqref{eq:E_c_vector} in the incompressible vector-valued case is lower semicontinuous in the $d$-topology for all $\epsilon \in (0,1)$.
\end{theorem}
\mybegproof
Lower semicontinuity amounts to proving that
\[
\Ee(c^\ast)\le \liminf_{n\to+\infty}\Ee(c_n)
\]
for all sequences $\{c_n\}_{n\in\mathbb{N}}\subset \LPhi$ that converge to $c^\ast \in \LPhi$. We split the proof into two cases.

Suppose $c^\ast\notin \Ce$. Then $\Ee(c^\ast)=+\infty$. Since $\Ce$ is closed with respect to $d$, eventually $c_n\notin \Ce$ because the complement of $\Ce$ is open, and thus $\liminf_{n\to+\infty}\Ee(c_n)=+\infty$.

Suppose now that $c^\ast\in \Ce$. If $c_n\in\Ce$ only for finitely many terms, then $\liminf_{n\to+\infty}\Ee(c_n)=+\infty$, and thus, for such sequences, the functional is lower semicontinuous. It only remains to consider the case where $c_n\in\Ce$ for infinitely many terms, having indices $I:=\left\{n_k \in \mathbb{N} \st c_{n_k} \in \Ce\right\}$. Since the terms of the sequence for which $n\notin I$ do not lower the limit inferior, we have
\[
\liminf_{n\to+\infty}\Ee(c_n) = \liminf_{k\to+\infty}\Ee(c_{n_k}).
\]
From Theorems \ref{thm:Jc_cont_proof} and \ref{thm:Ju_cont_proof}, $\Ee$ is continuous on $\Ce$, and thus $\Ee(c_{n_k})\to \Ee(c^\ast)$, proving that the functional is lower semicontinuous.

\myendproof

\begin{lemma}\label{lem:d_1_convex_proof}
Let Assumption \ref{ass:sindices_strict} hold and the metric be given by \eqref{eq:X_distance}. For all $\bar{c}\in \LPhi$, the map $c\rightarrow\frac{1}{2}d(\bar{c},c)^2$ is geodesically $1$-convex.
\end{lemma}
\mybegproof
Let
\[
\J(c):=\frac{1}{2}\,d(\bar{c},c)^2 = \frac{1}{2}\|\Psi(c) - \Psi(\bar{c})\|^2_\Ltwo,
\]
where $\Psi$ is given by \eqref{eq:psi_map_def}.
We have
\[
\J(\gamma_g(s)) = \frac{1}{2}\|(1-s)\Psi(c_0) + s\,\Psi(c_1)- \Psi(\bar{c})\|^2_\Ltwo,
\]
and thus, setting $u=\Psi(c_0)-\Psi(\bar{c})$ and $v=\Psi(c_1)-\Psi(\bar{c})$, we get
\[
\J(\gamma_g(s))=\frac12\|(1-s)u+s v\|_{L^2(\Omega)}^2.
\]
Expanding the square in the Hilbert space
$L^2(\Omega)$, using $(\cdot,\cdot)$ to denote the inner product, gives
\[
\|(1-s)u+s v\|_\Ltwo^2
=(1-s)^2\|u\|_\Ltwo^2+s^2\|v\|_\Ltwo^2+2s(1-s) ( u,v).
\]
Using the identity $2(u,v) = \|u\|_\Ltwo^2+\|v\|_\Ltwo^2-\|u-v\|_\Ltwo^2$ yields
\[
\|(1-s)u+s v\|_\Ltwo^2
=(1-s)\|u\|^2_\Ltwo+s\|v\|_\Ltwo^2-s(1-s)\|u-v\|_\Ltwo^2.
\]
The result follows by multiplying by $1/2$ and noting that
\[
\J(c_0) =\frac{1}{2}\|u\|_\Ltwo^2,\quad \J(c_1) = \frac{1}{2}\|v\|_\Ltwo^2, \quad
\|u-v\|_2^2=d(c_0,c_1)^2.
\]
\myendproof

\begin{theorem}\label{thm:Ee_lambda_conv_proof}
Let Assumption~\ref{ass:sindices_strict} hold, together with Assumption~\ref{ass:F_Hm1} in the scalar case or Assumption~\ref{ass:F_Hm1_vector} in the incompressible vector-valued case, and assume $\inf_{r>0}\left\{1+\frac{r\mu''(r)}{2\mu'(r)}\right\}>-\infty$. Then
\begin{enumerate}
    \item[(i)] if $\Phi$ is convex and $\mu''(r)\le0$, or  
    \item[(ii)] if $\Phi$ is concave and $\mu''(r) \ge \frac{4 \mu'(r)^2}{\mu(r)}$,
\end{enumerate}
$\Ee$ is geodesically $\lambda$-convex for all $\epsilon \in(0,1)$ with $\lambda=\inf_{r>0}\left\{1+\frac{r\mu''(r)}{2\mu'(r)}\right\}$.
\end{theorem}
\mybegproof
Given that $\Ee=+\infty$ on $\Ce^c$, we only need to prove geodesical $\lambda$-convexity of $\Ee$ on $\Ce$; instead of proving it directly on $\Ce$, we leverage the isometry of $\LPhi$ and $\Ltwo$ via the map $\Psi$ given in \eqref{eq:psi_map_def} and prove the equivalent statement that the functional
\begin{equation}\label{eq:E_l2}
\wEe : \Ltwo\longrightarrow\Rinf, \quad \wEe(\eta):=(\Ee \circ\Psi^{-1})(\eta) ,
\end{equation}
is $\lambda$-convex on the image of $\Ce$ under $\Psi$, that is, the closed convex interval $$\Le:=\left\{\eta \in \Ltwo \st \psi(\epsilon)\le \eta\le \psi(1/\epsilon)\right\},$$ where $\psi$ is given in \eqref{eq:psi_def}.
Specifically,
we derive sufficient conditions under which the function
\[
H(t):=\wEe(\eta_t)-\frac{\lambda}{2}\|\eta_t\|_{L^2(\Omega)}^2,\quad \eta_t := (1-t) \eta_0 + t\eta_1,
\]
is convex, which allows us to prove the $\lambda$-convexity of $\wEe$ in $\Le$. In fact, if 
\[
H(t) \le (1-t) H(0) + tH(1) = (1-t) \wEe(\eta_0) + t\wEe(\eta_1) -\frac{\lambda}{2}\left((1-t)\|\eta_0\|_\Ltwo^2 + t\|\eta_1\|_\Ltwo^2\right),
\]
using the identity
\[
\|\eta_t\|_\Ltwo^2 = (1-t)\|\eta_0\|_\Ltwo^2 + t\|\eta_1\|_\Ltwo^2 -t(1-t)\|\eta_1-\eta_0\|_\Ltwo^2,
\]
yields
\[
\wEe(\eta_t) \le  (1-t) \wEe(\eta_0) + t\wEe(\eta_1)-\frac{\lambda}{2}t(1-t)\|\eta_1 - \eta_0\|_\Ltwo^2,
\]
and thus the geodesical $\lambda$-convexity of $\Ee$
\[
\Ee(c_t) \le  (1-t) \Ee(c_0) + t\Ee(c_1)-\frac{\lambda}{2}t(1-t)d(c_1,c_0)^2,
\]
where $c_t := \Psi^{-1}(\eta_t)$ is the geodesic in $\Ce$.

In the sequel, we use the notation $K_\epsilon:=[\epsilon,1/\epsilon]$, $I_\epsilon:=[\psi(\epsilon), \psi(1/\epsilon)]$, together with $\xi:=\eta_1-\eta_0$.
Note that $\xi\in\Linf$; moreover, since $\mu\in C^2(0,+\infty)$ and $\omega\mu'$ is strictly positive on $K_\epsilon$, we have that $\psi$ is a $C^2$ diffeomorphism from $K_\epsilon$ onto $I_\epsilon$.

We first consider the functional
\[
\wJc:=-\omega(\Jc \circ\Psi^{-1})(\eta),
\]
and define, as in the proof of Theorem \ref{thm:Jc_cont_proof}, $h(r) := \Phi(r) -\frac{1}{2}\mu(r)r$, $q(s):=-\omega h(\psi^{-1}(s))$, $s := \psi(r)$,
so that
\[
\wJc(\eta)=\int_\Omega q(\eta)\dx.
\]
Since $\eta_t$ takes values in $I_\epsilon$, and $q\in C^2(I_\epsilon)$, the chain
rule gives
\begin{equation}\label{eq:J_c_2}
\dd{^2}{t^2} \wJc(\eta_t)=\int_\Omega q''(\eta_t)\xi^2\dx.
\end{equation}
It remains to compute $q''$: since $\displaystyle \dd{r}{s} = \frac{1}{\sqrt{\omega\mu'(r)/2}}$, we have
\begin{align}
q'(s) &= -\omega h'(r) \dd{r}{s} =\frac{1}{2}\omega \mu'(r)r \frac{1}{\sqrt{\omega\mu'(r)/2}} = r\sqrt{\omega\mu'(r)/2},\nonumber\\
q''(s) &= \dd{r}{s}\sqrt{\omega\mu'(r)/2} + r  \frac{1}{2\sqrt{\omega\mu'(r)/2}}\frac{\omega \mu''(r)}{2}\dd{r}{s} = 1 +\frac{r\mu''(r)}{2\mu'(r)}\label{eq:temp_q_2}.
\end{align}

We then consider the term
\begin{equation}\label{eq:Jueta}
\wJu(\eta) = \frac{\omega}{2}\int_\Omega \widetilde{\mu}(\eta)|\grad u[\eta]|^2\dx, \quad \widetilde{\mu}(\eta) := (\mu \circ\Psi^{-1})(\eta), 
\end{equation}
under the constraint
\begin{equation}\label{eq:Jueta_state}
\int_\Omega \widetilde{\mu}(\eta)\grad u[\eta]\cdot\grad v\dx = F(v),\quad \forall v \in \Honez.
\end{equation}
We first treat the scalar case in detail and outline the changes needed to extend the proof to the incompressible case later.
For $t\in[0,1]$, we first prove that the map $t\mapsto u[\eta_t]$
is continuous from $[0,1]$ into $V:=\Honez$ and twice differentiable from $(0,1)$ into $V$. Since $\xi\in L^\infty(\Omega)$ and $\widetilde\mu\in C^2(I_\epsilon)$, the map $t\mapsto \widetilde\mu(\eta_t)$ is $C^2((0,1);L^\infty(\Omega))$, with
\[
\dd{}{t}\widetilde\mu(\eta_t)=\widetilde\mu'(\eta_t)\xi,
\quad
\dd{^2}{t^2}\widetilde\mu(\eta_t)=\widetilde\mu''(\eta_t)\xi^2.
\]
For $t\in[0,1]$, define 
\[
A_t:V\to V^\ast,\quad\langle A_t w,v\rangle:=\int_\Omega\widetilde\mu(\eta_t)\grad w\cdot\grad v\dx,
\]
with coefficient bounds $\mu_-\le \widetilde\mu(\eta_t)\le\mu_+$ from \eqref{eq:mu_unif_bounds}.
Therefore $A_t$ is an isomorphism from $V$ onto $V^\ast$, uniformly in $t$.
Consequently, $t\mapsto A_t$ is $C^2((0,1);\mathcal L(V,V^\ast))$. Since $A_t$ is uniformly coercive and invertible for every $t$, the map $t\mapsto A_t^{-1}$ is $C^2((0,1);L(V^\ast,V))$, and hence $t\mapsto u[\eta_t]\in C^2((0,1);V)$.

Differentiating \eqref{eq:Jueta} along the geodesic gives 
\begin{equation}\label{eq:Jueta_fvar_temp}
\dd{}{t}\wJu(\eta_t)=\omega\int_\Omega\widetilde\mu(\eta_t)\grad u[\eta_t]\cdot\grad\auc_t\dx+\frac{\omega}{2}\int_\Omega\widetilde\mu'(\eta_t)|\grad u[\eta_t]|^2\,\xi\dx,
\end{equation}
where we have used the notation $\auc_t:=\dd{}{t}u[\eta_t]$.
Differentiating \eqref{eq:Jueta_state} yields 
\begin{equation}\label{eq:Jueta_state_var}
\int_\Omega \widetilde{\mu}(\eta_t)\grad \auc_t \cdot\grad v\dx = -\int_\Omega \widetilde{\mu}'(\eta_t) \grad u[\eta_t]\cdot \grad v \,\xi \dx,
\quad \forall v\in V.
\end{equation}
Using the above identity with $v=u[\eta_t]$, and substituting the first term in \eqref{eq:Jueta_fvar_temp} yields
\begin{equation*}%
\dd{}{t}\wJu(\eta_t) = -\frac{\omega}{2}\int_\Omega\widetilde{\mu}'(\eta_t)|\grad u[\eta_t]|^2 \xi\dx.
\end{equation*}
Differentiating once more
\[
\dd{^2}{t^2}\wJu(\eta_t) = -\omega\int_\Omega \widetilde{\mu}'(\eta_t)\grad u[\eta_t]\cdot \grad \auc_t \,\xi\dx-\frac{\omega}{2}\int_\Omega\widetilde{\mu}''(\eta_t)|\grad u[\eta_t]|^2 \,\xi^2\dx.
\]
Using $v=\auc_t$ in \eqref{eq:Jueta_state_var} and substituting yields
\begin{equation}\label{eq:J_u_2}
\dd{^2}{t^2}\wJu(\eta_t) = \omega\int_\Omega \widetilde{\mu}(\eta_t)|\grad \auc_t|^2 \dx-\frac{\omega}{2}\int_\Omega\widetilde{\mu}''(\eta_t)|\grad u[\eta_t]|^2 \,\xi^2\dx.
\end{equation}
Using the notation $\widetilde{\mu}(s) = \mu(r)$, with $s=\psi(r)$, straightforward calculations show
\begin{align}
\widetilde{\mu}'(s) &= \frac{\mu'(r)}{\sqrt{\omega \mu'(r)/2}} = \omega\sqrt{2 \omega \mu'(r)},\label{eq:temp_mut_1}\\
\widetilde{\mu}''(s) &= \frac{\omega}{2\sqrt{2 \omega \mu'(r)}}2\omega \mu''(r)\frac{1}{\sqrt{\omega \mu'(r)/2}} = \omega\frac{\mu''(r)}{\mu'(r)}.\label{eq:temp_mut_2}
\end{align}

We can now prove sufficient conditions for the convexity of $H(t)$; from \eqref{eq:J_c_2} and \eqref{eq:J_u_2} we obtain
\[
H''(t)=\int_\Omega (q''(\eta_t) -\lambda)\xi^2\dx + \omega\int_\Omega\widetilde\mu(\eta_t)|\grad\auc_t|^2\dx
-\frac{\omega}{2}\int_\Omega\widetilde\mu''(\eta_t)|\grad u[\eta_t]|^2\xi^2\dx.
\]
We separate the cases of convex and concave $\Phi$. In the convex case, $\omega=1$, and thus we can drop the second integral in the formula of $H''(t)$ and find the pointwise sufficient condition
\[
\lambda \le q''(\eta_t) -\frac{1}{2}\pos{\widetilde{\mu}''(\eta_t)}|\grad u[\eta_t]|^2,
\]
where $\pos{x} = \max(x,0)$.
From \eqref{eq:temp_q_2} and \eqref{eq:temp_mut_2}, and using the fact that $\grad u[\eta_t] =\grad u[c_t]$ with $c_t=\Psi^{-1}(\eta_t)$ the geodesic in $\Ce$, we find the sufficient condition
\begin{equation}\label{eq:lambda_geo_temp_conv}
\lambda_{c_t} = \inf_{t\in[0,1]}\essinf_{x\in\Omega}\left\{1+\frac{c_t(x)\mu''(c_t(x))}{2\mu'(c_t(x))}-\pos{\frac{\mu''(c_t(x))}{2\mu'(c_t(x))}}|\grad u[c_t](x)|^2\right\} > -\infty.
\end{equation}
Since $\mu'(r)\ge0$ when $\Phi$ is convex, then, if $\mu''(r)\le 0$, we can obtain a uniform expression for $\lambda$ that holds true for all geodesics $c_t$
\[
\lambda = \inf_{r\in [\epsilon,1/\epsilon]}\left\{1+\frac{r\mu''(r)}{2\mu'(r)}\right\},
\]
which proves statement (i) after taking the infimum over all $\epsilon>0$. %

In the concave case, $\omega=-1$, and thus we must estimate the second term in $H''(t)$.
Testing \eqref{eq:Jueta_state_var} with $\auc_t$ and using a Cauchy--Schwarz argument we obtain
\begin{equation}\label{eq:cs_lambdaconv_conc}
\begin{aligned}
\int_\Omega \widetilde{\mu}(\eta_t)|\grad \auc_t|^2\dx &= \left|\int_\Omega \widetilde{\mu}'(\eta_t) \grad u[\eta_t]\cdot \grad \auc_t \,\xi \dx\right|\\
&=\left|\int_\Omega \frac{\widetilde{\mu}'(\eta_t)}{\sqrt{\widetilde{\mu}(\eta_t)}} \grad u[\eta_t] \xi \cdot \sqrt{\widetilde{\mu}(\eta_t)}\grad \auc_t  \dx\right|\\
&\leq\left(\int_\Omega \frac{\widetilde{\mu}'(\eta_t)^2}{\widetilde{\mu}(\eta_t)} |\grad u[\eta_t]|^2 \xi^2\dx\right)^{1/2} \left(\int_\Omega\widetilde{\mu}(\eta_t)|\grad \auc_t|^2  \dx\right)^{1/2},
\end{aligned}
\end{equation}
and thus
\[
\int_\Omega \widetilde{\mu}(\eta_t)|\grad \auc_t|^2\dx \le \int_\Omega \frac{\widetilde{\mu}'(\eta_t)^2}{\widetilde{\mu}(\eta_t)} |\grad u[\eta_t]|^2 \xi^2\dx.
\]
Plugging the above inequality into the formula for $H''(t)$, we arrive at
\[
H''(t)\ge \int_\Omega (q''(\eta_t) -\lambda)\xi^2\dx +\int_\Omega \left(\frac{1}{2}\widetilde{\mu}''(\eta_t)-\frac{\widetilde{\mu}'(\eta_t)^2}{\widetilde{\mu}(\eta_t)}\right) |\grad u[\eta_t]|^2 \xi^2\dx.
\]
Therefore, a sufficient pointwise condition is
\[
\lambda \le q''(\eta_t) - \pos{\frac{\widetilde{\mu}'(\eta_t)^2}{\widetilde{\mu}(\eta_t)}-\frac{1}{2}\widetilde{\mu}''(\eta_t)} |\grad u[\eta_t]|^2.
\]
Substituting the expressions for $q''(\eta_t)$,  $\widetilde{\mu}'(\eta_t)$ and $\widetilde{\mu}''(\eta_t)$ from \eqref{eq:temp_q_2}, \eqref{eq:temp_mut_1} and \eqref{eq:temp_mut_2} we find the sufficient condition
\begin{equation}\label{eq:lambda_geo_temp_conc}
\lambda_{c_t} = \inf_{t\in[0,1]}\essinf_{x\in\Omega}\left\{1+\frac{c_t(x)\mu''(c_t(x))}{2\mu'(c_t(x))}-\pos{\frac{\mu''(c_t(x))}{2\mu'(c_t(x))}-\frac{2\mu'(c_t(x))}{\mu(c_t(x))}}|\grad u[c_t](x)|^2\right\} > -\infty.
\end{equation}
Since in the concave case $\mu'(r)\le0$, if $\mu''(r)\mu(r) \ge 4\mu'(r)^2$ we can obtain the simpler expression
\[
\lambda = \inf_{r\in [\epsilon,1/\epsilon]}\left\{1+\frac{r\mu''(r)}{2\mu'(r)}\right\},
\]
and statement (ii) is proven.

The incompressible vector-valued case is identical, up to replacing the gradient by
the symmetric gradient. The differentiability of $t\mapsto u[\eta_t]$ in the Stokes case follows
in the same way, by working on the divergence-free space
$H^1_{0,\mathrm{div}}(\Omega;\mathbb R^d)$, where the weighted Stokes
operator is uniformly coercive by Korn's inequality and the bounds
$0<\mu_-\le \widetilde\mu(\eta_t)\le \mu_+$ hold. Let
$(\bu[\eta_t],p[\eta_t])\in H^1_0(\Omega;\mathbb R^d)\times Q$ solve the weighted
Stokes system. For a direction $\xi$, let
$(\bauc_t,\pauc_t)\in \vHonez\times Q$ be the solution of
the linearized Stokes system
\[
\begin{aligned}
\int_\Omega
\widetilde\mu(\eta_t)\,\eps(\bauc_t):\eps(\bv)\dx
-
\int_\Omega \pauc_t \grad\cdot \bv\dx
&=
-
\int_\Omega
\widetilde\mu'(\eta_t)
\eps(\bu[\eta_t]):\eps(\bv)\,\xi\dx ,\\
\int_\Omega q\,\diver\bauc_t\dx &=0 .
\end{aligned}
\]
The pressure terms vanish when testing with $\bu[\eta_t]$ or $\bauc_t$  since both are divergence-free. Therefore
\[
\dd{^2}{t^2}\wJbu(\eta_t)
=
\omega\int_\Omega
\widetilde\mu(\eta_t)|\eps(\bauc_t)|^2\dx
-
\frac{\omega}{2}
\int_\Omega
\widetilde\mu''(\eta_t)|\eps(\bu[\eta_t])|^2\xi^2\dx.
\]
The same estimates used in the scalar case imply that $H$ is convex.
\myendproof

\begin{theorem}\label{thm:ode_bound_proof}
Under Assumptions \ref{ass:Phi_concace_or_convex_strictly}, \ref{ass:sindices}, and \ref{ass:fem_space}, problem \eqref{eq:ode_fem} has a unique global solution $c_h\in C^1([0,+\infty);C_h)$ and $c_h(t) \in C^+_h$.
Moreover, for each fixed $T <+\infty$, there exists $\epsilon(T)> 0$ such that
\[
\epsilon(T) \le c_h(t)\le 1/\epsilon(T), \quad \forall\, 0\le t\le T.
\]
\end{theorem}

\mybegproof
We prove the result in the scalar case. The vector-valued case is identical after replacing
$\grad u_h$ by $\eps(\bu_h)$ and working in the discrete divergence-free space $\bUhdiv$. It is omitted for brevity.
We first prove that the function
\begin{equation}\label{eq:fode}
\mathcal{F}_h:C_h^+\to C_h,\quad \mathcal{F}_h(c_h):=|\grad u_h[c_h]|^2-c_h,
\end{equation}
is locally Lipschitz. To this end, for $0<m<M$ we define
\[
K_{m,M}:=\{c_h\in C_h:\ m\le c_h\le M\},
\]
which is compact in $C_h$. Moreover, for each $c_{h,0}\in C^+_h$, the elliptic problem that defines $\grad u_h[c_h]$ for each $c_h$ in the neighborhood of $c_{h,0}$ is well-posed, since  $c_{h,0} \in K_{m,M}$ for some $m,M$. We take two functions $c_i\in K_{m,M}$ and use the notation $u_i:=u_h[c_i]$, $i=1,2$, i.e.
\[
\int_\Omega\mu(c_i)\grad u_i\cdot\grad v_h\dx=F(v_h),\quad \forall v_h\in U_h.
\]
Combining the equations, we obtain
\[
\int_\Omega\mu(c_1)\grad(u_1-u_2)\cdot\grad v_h\dx=\int_\Omega(\mu(c_2)-\mu(c_1))\grad u_2\cdot\grad v_h\dx.
\]
Choosing $v_h=u_1-u_2$, and using the lower bound
\[
        \mu(c_1)\ge \mu_{m,M}^{\min}
        :=
        \min_{r\in[m,M]}\mu(r)>0,
\]
yields
\[
\mu_{m,M}^{\min}\|\grad u_1-\grad u_2\|_{\Ltwo}\le\|\mu(c_2)-\mu(c_1)\|_{\Linf}\|\grad u_2\|_{\Ltwo}.
\]
Since $\mu\in C^1([m,M])$,
\[
\|\mu(c_2)-\mu(c_1)\|_{\Linf}\le L^\mu_{m,M}\|c_2-c_1\|_{\Linf},\quad L^\mu_{m,M}:=\max_{r\in[m,M]}|\mu'(r)|< +\infty.
\]
Combining with the standard estimate
\[
\|\grad u_2\|_{\Ltwo}\le \frac{\|F\|_{H^{-1}(\Omega)}}{\mu_{m,M}^{\min}},
\]
we obtain
\begin{equation*}%
\|\grad u_1-\grad u_2\|_{\Ltwo}\le L_u\|c_1-c_2\|_{\Linf}, \quad L_u :=   L^\mu_{m,M}\frac{\|F\|_{H^{-1}(\Omega)}}{(\mu_{m,M}^{\min})^2},
\end{equation*}
for all $c_1,c_2\in K_{m,M}$. For the fixed mesh, let $C_{\mathrm{inv},h}$ be a constant in the inverse inequality $\|\grad v_h\|_{\Linf}\le C_{\mathrm{inv},h}\|\grad v_h\|_{\Ltwo}$ for all $v_h\in U_h$, obtained by applying \cite[Lemma~12.1]{ernguermondI} on each element. The local Lipschitz continuity of $\mathcal{F}_h$ then follows:
\[
\begin{aligned}
\|\mathcal{F}_h(c_1)-\mathcal{F}_h(c_2)\|_{\Linf}
&\le \||\grad u_1|^2-|\grad u_2|^2\|_{\Linf}+\|c_1-c_2\|_{\Linf}\\
&\le(\|\grad u_1\|_{\Linf}+\|\grad u_2\|_{\Linf})\|\grad u_1-\grad u_2\|_{\Linf}+\|c_1-c_2\|_{\Linf}\\
&\le C_{\mathrm{inv},h}(\|\grad u_1\|_{\Linf}+\|\grad u_2\|_{\Linf})\|\grad u_1-\grad u_2\|_{\Ltwo}+\|c_1-c_2\|_{\Linf}\\
&\le \bigg(C_{\mathrm{inv},h}L_u(\|\grad u_1\|_{\Linf}+\|\grad u_2\|_{\Linf}) + 1\bigg)\|c_1-c_2\|_{\Linf},\\
&\le L\,\|c_1-c_2\|_{\Linf},
\end{aligned}
\]
where
\[
L:= 2 C_{\mathrm{inv},h}^2L_u\frac{\|F\|_{H^{-1}(\Omega)}}{\mu_{m,M}^{\min}} + 1.
\]

Since $C_h^+$ is open in the finite-dimensional space $C_h$, and
$\mathcal{F}_h:C_h^+\to C_h$ is locally Lipschitz, the Picard--Lindel\"of theorem (see e.g. \cite[Theorem 2.2]{teschl}) gives a
unique maximal solution
\[
c_h\in C^1([0,T_{\max});C_h), \quad c_h(t)\in C_h^+,
\]
for some $0<T_{\max}\le+\infty$. %
For $t < T_{\max}$, the explicit solution of the ODE is
\begin{equation}\label{eq:ode_expl}
c_h(t)=c_{h,0}\,e^{-t}+\int_0^t e^{-(t-s)}|\grad u_h[c_h(s)]|^2\,ds.
\end{equation}
Since $c_{h,0}\in C^+_h$, we have $m_0 := \essinf_{x\in\Omega}c_{h,0}(x) > 0$, $M_0 := \esssup_{x\in\Omega}c_{h,0}(x) < +\infty$, and %
\[
c_h(t)\ge m_T:=m_0\,e^{-T} > 0,\quad \forall\, 0\le t \le T< T_{\max},
\]
and thus we cannot have finite-time extinction since the trajectories cannot approach the boundary of $C^+_h$ in finite time.
We then estimate an upper bound on $[0,T_{\max})$, and use the notation $y(t) := \|c_h(t)\|_\Linf$. By coercivity of the bilinear form and the inverse estimate we get
\[
\|\grad u_h[c_h(t)]\|_{\Linf}^2\le\left(\frac{C_{\mathrm{inv},h}\|F\|_{H^{-1}(\Omega)}}{\mu_{\min}(t)}\right)^2,\quad
\mu_{\min}(t):=\essinf_{x\in\Omega}\mu(c_h(x,t)) > 0.
\]
In the case of $\Phi$ convex, $\mu$ is increasing and thus $\mu_{\min}(t) \ge \mu(m_T)$ for $t\le T$; from  \eqref{eq:ode_expl} we obtain
\[
y(t) \le K_{T,\mathrm{conv}} + (M_0 - K_{T,\mathrm{conv}})e^{-t}, \quad 0\le t \le T, \quad K_{T,\mathrm{conv}}:=\left(\frac{C_{\mathrm{inv},h}\|F\|_{H^{-1}(\Omega)}}{\mu(m_T)}\right)^2.
\]
In the case of concave $\Phi$,  $\mu$ is decreasing; from Lemma \ref{lem:phiprime_bounds} we have
\[
\mu(c_h(t))\ge \mu(y(t)) \ge \mu(m_T)\left(\frac{m_T}{y(t)}\right)^\frac{q}{2}, \quad q:=1-\smphi\in[0,1).
\]
Therefore
\[
\|\grad u_h[c_h(t)]\|_{\Linf}^2\le K_{T,\mathrm{conc}}\,y(t)^q,\quad
K_{T,\mathrm{conc}}:=\frac{C_{\mathrm{inv},h}^2\|F\|_{H^{-1}(\Omega)}^2}{\mu(m_T)^2m^q_T}.
\]
Plugging the above estimate into \eqref{eq:ode_expl} yields
\[
y(t)\le M_0 + K_{T,\mathrm{conc}}\int_0^t y(s)^q\,ds.
\]
Since $0\le q<1$, Bihari's inequality \cite[Eq. (7)]{bihari} gives
\[
y(t)\le\left(M_0^{1-q}+K_{T,\mathrm{conc}}(1-q)t\right)^{1/(1-q)},
\quad 0\le t\le T.
\]
Therefore, in both regimes, for every $T<T_{\max}$, there exists $B_T<+\infty$ such that
\[
m_T\le c_h(x,t)\le B_T,\quad 0\le t\le T,
\]
and we obtain
\[
\epsilon(T):=\min\left\{m_T,\frac{1}{B_T}\right\}.
\]

We then prove global existence, and assume for contradiction that $T_{\max}<+\infty$. Defining $m_*:=m_0e^{-T_{\max}} > 0$, since $m_*\le m_T$, in the convex case we have
\[
K_{T,\mathrm{conv}} \le K_{*,\mathrm{conv}} := \left(\frac{C_{\mathrm{inv},h}\|F\|_{H^{-1}(\Omega)}}{\mu(m_*)}\right)^2 < +\infty,
\]
uniformly for $T<T_{\max}$. In the concave case, using the upper bound of statement (ii) in Lemma \ref{lem:phiprime_bounds} with $m_*\le m_T$, we have $\mu(m_T)m_T^{q/2}\ge \mu(m_*) m_*^{q/2}$, and thus
\[
K_{T,\mathrm{conc}}\le K_{*,\mathrm{conc}}:=\frac{C_{\mathrm{inv},h}^2\|F\|_{H^{-1}(\Omega)}^2}{\mu(m_*)^2m^q_*} < +\infty.
\]
Therefore, from the previous estimates
we obtain uniform lower and upper bounds
\[
0 < m_*\le c_h(x,t)\le B_* < +\infty,\quad 0\le t<T_{\max},
\]
where
\[
B_* := \begin{cases}
\max\left\{M_0,K_{*,\mathrm{conv}}\right\}, &\Phi \text{ convex}\\
\left(M_0^{1-q}+K_{*,\mathrm{conc}}(1-q)T_{\max}\right)^\frac{1}{1-q}, &\Phi \text{ concave }
\end{cases}
\]
Therefore $c_h(t)$ remains in the compact set $K_{m_*,B_*}\subset C_h^+ $. This contradicts the finite-dimensional continuation alternative for ODEs \cite[Corollary 2.16]{teschl}, which says that a maximal solution with finite maximal time must leave every compact subset of the open domain $C_h^+$. Hence $T_{\max}=+\infty$.
\myendproof

\begin{theorem}\label{thm:mielke_lambda_discrete_proof}
Let the assumptions of Theorem \ref{thm:ode_bound} hold; assume further that
$\inf_{r>0}\left\{1+\frac{r\mu''(r)}{2\mu'(r)}\right\}>-\infty$. Then  
\begin{enumerate}
    \item[(i)] if $\Phi$ is convex and $\mu''(r)\le0$, or  
    \item[(ii)] if $\Phi$ is concave and $\mu''(r) \ge \frac{4 \mu'(r)^2}{\mu(r)}$,
\end{enumerate}
the following inequality holds
with $\lambda=\inf_{r>0}\left\{1+\frac{r\mu''(r)}{2\mu'(r)}\right\}$
\[
\mathcal{B}(c_h,z_h) \ge \lambda \mathcal{A}(c_h,z_h) \quad \forall c_h\in C^+_h,\quad \forall z_h \in C_h,
\]
where
\[
\begin{aligned}
\mathcal{A}(c_h,z_h)&:=\langle \mathcal{G}(c_h)z_h,z_h\rangle
=\int_\Omega \frac{\omega}{2}\,\mu'(c_h)\,z_h^2\,\dx,\\
\mathcal{B}(c_h,z_h)&:=\langle \mathcal{G}(c_h)z_h,\mathcal{DR}_h(c_h)z_h\rangle+\frac12\langle \mathcal{DG}(c_h)[\mathcal{R}_h(c_h)]z_h,z_h\rangle,
\end{aligned}
\]
with $\mathcal{G}(c_h)$ defined by \eqref{eq:def_G} and $\mathcal{R}_h(c_h)=c_h-|Du_h[c_h]|^2$. In the scalar case, $\mathcal{R}_h(c_h) = -\mathcal{F}_h(c_h)$, with $\mathcal{F}_h$ given by \eqref{eq:fode}.
\end{theorem}
\mybegproof
We derive the estimate without distinguishing between the scalar and the vector-valued case, and, with an abuse of notation, we will always use $Du:Dv$ to denote the (Frobenius) inner product between (symmetric) gradients. Differentiating $\mathcal{R}_h(c_h)=c_h-|Du_h[c_h]|^2$ in the direction $z_h\in C_h$ we get
\[
\mathcal{DR}_h(c_h)\,z_h = \left.\dd{}{\evar}\mathcal{R}_h(c_h+\evar z_h)\right|_{\evar=0}= z_h - 2\,Du_h[c_h]:D\auc_h,\quad
\auc_h:=\left.\dd{}{\evar}u_h[c_h+\evar z_h]\right|_{\evar=0},
\]
and hence
\begin{align}
\langle \mathcal{G}(c_h)z_h,\mathcal{DR}_h(c_h)z_h\rangle
&=\int_\Omega \frac{\omega}{2}\mu'(c_h)\,z_h\,(z_h-2\,Du_h[c_h]:D\auc_h)\dx \notag\\
&=\int_\Omega \frac{\omega}{2}\mu'(c_h)\,z_h^2\dx
-\omega\int_\Omega \mu'(c_h)\,z_h\,Du_h[c_h]:D\auc_h\dx. \label{eq:first_part}
\end{align}
Furthermore
\begin{equation}\label{eq:second_part}
\frac12\langle \mathcal{DG}(c_h)[\mathcal{R}_h(c_h)]z_h,z_h\rangle
=\int_\Omega \frac{\omega}{4}\,\mu''(c_h)\,\mathcal{R}_h(c_h)\,z_h^2\dx.
\end{equation}
Combining \eqref{eq:Adef}, \eqref{eq:first_part}, and \eqref{eq:second_part} yields
\begin{equation}\label{eq:B_expand}
\mathcal{B}(c_h,z_h)=\mathcal{A}(c_h,z_h)
-\omega\int_\Omega \mu'(c_h)\,z_h\,Du_h[c_h]:D\auc_h\dx
+\int_\Omega \frac{\omega}{4}\,\mu''(c_h)\,\mathcal{R}_h(c_h)\,z_h^2\dx.
\end{equation}
Using \eqref{eq:sens_scalar} or \eqref{eq:sens_vector} (depending on the scalar or vector case) tested with $\auc_h$, we get 
\begin{equation}\label{eq:auc}
-\int_\Omega \mu'(c_h)\,z_h\,Du_h[c_h]:D\auc_h\dx
=
\int_\Omega \mu(c_h)\,|D\auc_h|^2\dx,
\end{equation}
and arrive at
\begin{equation}\label{eq:B_decomp}
\mathcal{B}(c_h,z_h)=\mathcal{A}(c_h,z_h)
+\omega\int_\Omega \mu(c_h)\,|D\auc_h|^2\dx
+\int_\Omega \frac{\omega}{4}\,\mu''(c_h)\,\mathcal{R}_h(c_h)\,z_h^2\dx.
\end{equation}

Now consider the convex case for $\Phi$, for which $\omega=1$. Then, plugging
\[
\omega\int \mu(c_h)|D\auc_h|^2\dx\ge 0,\quad\omega\mu''(c_h)\,\mathcal{R}_h(c_h) = \mu''(c_h)\,(c_h-|Du[c_h]|^2),
\]
into \eqref{eq:B_decomp}
yields
\[
\begin{aligned}
\mathcal{B}(c_h,z_h)&\ge\int_\Omega \frac{\mu'(c_h)}{2}\left(1+\frac{c_h\,\mu''(c_h)}{2\mu'(c_h)} -\pos{\frac{\mu''(c_h)}{2\mu'(c_h)}}|Du_h[c_h]|^2\right)z_h^2\dx\\
&\ge\lambda_{c_h}\int_\Omega \frac{\mu'(c_h)}{2}\,z_h^2\dx\\
&=\lambda_{c_h}\,\mathcal{A}(c_h,z_h),
\end{aligned}
\]
where
\begin{equation}\label{eq:convex_lambdach}
\lambda_{c_h} = \essinf_{x\in\Omega}\left\{1+\frac{c_h(x)\,\mu''(c_h(x))}{2\mu'(c_h(x))} -\pos{\frac{\mu''(c_h(x))}{2\mu'(c_h(x))}}|Du_h[c_h](x)|^2\right\}.
\end{equation}
We can then use the same arguments as in the proof of Theorem \ref{thm:Ee_lambda_conv_proof} and hence statement (i) follows.

We then prove (ii). To this end, we need the following inequality \begin{equation}\label{eq:CS_phi_prime}
\int_\Omega \mu(c_h)\,|D\auc_h|^2\dx
\le
\int_\Omega \frac{\mu'(c_h)^2}{\mu(c_h)}\,|Du_h[c_h]|^2\,z_h^2\dx,
\end{equation}
that can be proven with the same arguments used in \eqref{eq:cs_lambdaconv_conc}.
Plugging \eqref{eq:CS_phi_prime} into \eqref{eq:B_decomp}, using $\omega=-1$,
\[
\begin{aligned}
\mathcal{B}(c_h,z_h)
&\ge
\int_\Omega
\frac{-\mu'(c_h)}{2}\left(
1+\frac{c_h\mu''(c_h)}{2\mu'(c_h)}
\right)
z_h^2\,dx
+
\int_\Omega
\left(
\frac14\mu''(c_h)
-
\frac{\mu'(c_h)^2}{\mu(c_h)}
\right)
|D u_h[c_h]|^2z_h^2\dx\\
&=
\int_\Omega
\frac{-\mu'(c_h)}{2}\left[
1+\frac{c_h\mu''(c_h)}{2\mu'(c_h)}
+
\left(
\frac{2\mu'(c_h)}{\mu(c_h)}
-
\frac{\mu''(c_h)}{2\mu'(c_h)}
\right)|D u_h[c_h]|^2
\right]
z_h^2\dx\\
&\ge
\int_\Omega
\frac{-\mu'(c_h)}{2}\left[
1+\frac{c_h\mu''(c_h)}{2\mu'(c_h)}
-
\pos{
\frac{\mu''(c_h)}{2\mu'(c_h)}
-
\frac{2\mu'(c_h)}{\mu(c_h)}
}|D u_h[c_h]|^2
\right]
z_h^2\dx\\
&\ge
\lambda_{c_h}\mathcal{A}(c_h,z_h),
\end{aligned}
\]
where
\begin{equation}\label{eq:concave_lambdach}
\lambda_{c_h} = \essinf_{x\in\Omega}\left\{
1+\frac{c_h(x)\mu''(c_h(x))}{2\mu'(c_h(x))}
-
\pos{
\frac{\mu''(c_h(x))\mu(c_h(x))-4\mu'(c_h(x))^2}{2\mu'(c_h(x))\mu(c_h(x))}
}
|D u_h[c_h](x)|^2\right\}.
\end{equation}
Again, using the same arguments as in the proof of Theorem \ref{thm:Ee_lambda_conv_proof}, statement (ii) then follows. 
\myendproof

\begin{lemma}\label{lem:coefficient-residual_proof}
Let the assumptions of Lemma \ref{lem:DK-shifted} hold. Then, if %
\[
\frac13\le\smphi\le\spphi<1, \quad \text{or}\quad 1<\smphi\le\spphi\le3,
\]
we have
\begin{equation}\label{eq:phistar_P_bound_proof}
\phi^*_{|P|}\left(|\mu(r)-\mu(|P|^2)|\,|P|\right)
\lesssim|\mu'(r)|\,(r-|P|^2)^2,\quad \forall r>0,\,P \in \mathbb{R}^{Nd}.
\end{equation}
\end{lemma}
\mybegproof
If $P=0$, then the left-hand side is zero. Assume $P\neq0$, and use the notation $a:=|P|$, $b:=\sqrt r$, $z:=|a-b|$, and $m(t):=\mu(t^2)=\phi'(t)/t$.
Then
\[
|\mu(r)-\mu(|P|^2)|\,|P|=|m(b)-m(a)|\,a,
\quad
|\mu'(r)|\,(r-|P|^2)^2=|\mu'(b^2)|\,(b^2-a^2)^2,
\]
and the desired estimate is equivalent to
\[
\phi_a^*(|m(b)-m(a)|\,a)\lesssim |\mu'(b^2)|\,(b^2-a^2)^2.
\]
Let
\[
s(t):=\frac{t\phi''(t)}{\phi'(t)}, \quad \smphi\le s(t)\le \spphi.
\]
Since
\[
m'(t)=\frac{t\phi''(t)-\phi'(t)}{t^2}, 
\quad
\mu'(t^2)=\frac{m'(t)}{2t},
\]
we have 
\[
(\smphi -1)\frac{\phi'(t)}{t}
\le t m'(t)
\le
(\spphi -1)\frac{\phi'(t)}{t},
\quad
\frac{\smphi -1}{2}\frac{\phi'(t)}{t^3}
\le
\mu'(t^2)
\le
\frac{\spphi -1}{2}\frac{\phi'(t)}{t^3}.
\]
The extra assumptions $\smphi>1$ in the convex case and $\spphi<1$ in the concave case ensure that $s(t) -1$ is uniformly separated from zero and allow us to consider the equivalences
 \begin{align}
|t m'(t)|&\simeq\frac{\phi'(t)}{t}
,\label{eq:elem_mprime}\\ 
|\mu'(t^2)|&\simeq\frac{\phi'(t)}{t^3}.\label{eq:elem_muprime}
\end{align}

We split the proof into two cases. First we assume $b\ge a/2$ and prove that
\begin{equation}\label{eq:temp_m_bound_1}
|m(b)-m(a)|\,a\lesssim \phi_a'(z).
\end{equation}
When $a/2\le b\le2a$, $z = |a-b|\le a$, and $\rho \simeq a \simeq a+z$ for all $\min\{a,b\}\le \rho \le \max\{a,b\}$. By the
mean value theorem, using \eqref{eq:elem_mprime}
\[
|m(b)-m(a)|\,a\le a\int_{\min\{a,b\}}^{\max\{a,b\}} |m'(\rho)|\,d\rho
\simeq \int_{\min\{a,b\}}^{\max\{a,b\}} a\frac{\phi'(\rho)}{\rho^2}\,d\rho
\simeq z\frac{\phi'(a+z)}{a+z}
=\phi_a'(z).
\]
When $b>2a$, then $z = b-a$ and $z\simeq b$. In the case of convex $\Phi$, $m$ is
increasing, and thus, using the definition of $\phi_a'(z)$ in \eqref{eq:shifted_nfunc}
\[
|m(b)-m(a)|\,a = (m(b)-m(a))\,a\le m(b)a \le m(b)b =\phi'(b) = \frac{b}{b-a}\phi_a'(z) \le 2 \phi_a'(z) \lesssim \phi_a'(z).
\]
When $\Phi$ is concave, $m$ is decreasing, and
\[
|m(b)-m(a)|\,a = (m(a)-m(b))\,a\le m(a)a = \phi'(a) \le \phi'(b) \simeq \phi_a'(z).
\]
Thus, in both cases, \eqref{eq:temp_m_bound_1} holds and, since $\phi_a^*$ is increasing, the uniform $\Delta_2$-condition for
$\phi_a^*$ gives
\[
\phi_a^*(|m(b)-m(a)|a)\lesssim \phi_a^*(\phi_a'(z)) \simeq \phi_a(z),
\]
where we have used the standard Orlicz identity $\phi_a^*(\phi_a'(z))\simeq \phi_a(z)$, see e.g. \cite[Eq. (2.6)]{DieningKreuzer2008}. By (ii) in Lemma \ref{lem:DK-shifted},
\[
\phi_a(z)\simeq z^2\frac{\phi'(a+z)}{a+z}.
\]
Therefore, since $b\ge a/2$, we have $a+z\simeq b$, $a+b\simeq b$; using \eqref{eq:elem_muprime} yields
\[
\phi_a(z) \simeq z^2\frac{\phi'(b)}{b} \simeq|\mu'(b^2)|z^2b^2\simeq|\mu'(b^2)|\,(b^2-a^2)^2,
\]
which proves the desired estimate for $b\ge a/2$.

Now assume $b<a/2$. This is where we use the additional restrictions $\smphi\ge1/3$ in the concave case and $\spphi\le 3$ in the convex case. First consider the convex case; then $m(t)$ is increasing and
\[
|m(b)-m(a)|\,a = (m(a)-m(b))\,a\le m(a)a =\phi'(a).
\]
Moreover, since by definition $\phi'_a(a) = \phi'(2a)/2$, the growth estimate in Lemma \ref{lem:phiprime_bounds} gives
\[
\phi'(a) \leq 2^{-\smphi}\phi'(2a) = 2^{1-\smphi}\phi_a'(a).
\]
Using the monotonicity of $\phi^*_a$,
\[
\phi^*_a(|m(b)-m(a)|\,a) 
\lesssim
\phi_a^*(\phi'_a(a))
\lesssim
\phi_a(a)
\simeq a\phi'(a),
\]
where the last step follows by combining (i) and (ii) in Lemma \ref{lem:DK-shifted}.
Using the growth estimate from Lemma \ref{lem:phiprime_bounds} again and $\spphi\le3$ yields
\[
a\phi'(a)\le a\phi'(b)\left(\frac ab\right)^\spphi\leq \frac{\phi'(b)}{b^3}a^4.
\]
Since $b<a/2$ we have $
(b^2-a^2)^2\simeq a^4$,
and thus, using \eqref{eq:elem_muprime}
\[
\phi_a^*(|m(b)-m(a)|\,a)
\lesssim
|\mu'(b^2)|\,(b^2-a^2)^2.
\]
Finally, consider the concave case; since $m(t)$ is decreasing we have
\[
 |m(b)-m(a)|a \le m(b)a
 = \frac{\phi'(b)}{b}a = \lambda\phi'(b), \quad \lambda:=\frac{a}{b}>2.
\]
We then prove that
\[
(\phi_a')^{-1}(\lambda\phi'(b))
\lesssim b\lambda^{1/\smphi}.
\]
To this end, we set
\[
y:=C b\lambda^{1/\smphi}, \quad C:=2^{1/\smphi}.
\]
Since $\smphi\le1$ and $\lambda>1$, we have $y\ge C b\lambda \ge a$ and $a+y\le 2y$, yielding
\[
\phi_a'(y)
=
y\frac{\phi'(a+y)}{a+y}
\ge
\frac12 \phi'(y).
\]
Since $y\ge 2 b$, the two-point growth estimate for $\phi'$ in Lemma \ref{lem:phiprime_bounds}
\[
\phi'(y)
\ge
\left(\frac{y}{b}\right)^\smphi \phi'(b)
=
C^\smphi \lambda \phi'(b)
=
2\lambda\phi'(b).
\]
Thus
\[
\phi_a'(y)\ge \lambda\phi'(b),
\]
and since $\phi_a'$ is increasing,
\[
(\phi_a')^{-1}(\lambda\phi'(b))
\le y
=
C b\lambda^{1/\smphi}.
\]
Using the basic inequality from \eqref{eq:phi_complement}
\[
\phi_a^*(t)\le (\phi_a')^{-1}(t)t,\quad t\ge0,
\]
with $t=\lambda\phi'(b)$, we get, using the extra assumption $\smphi\ge1/3$
\[
\phi_a^*(\lambda\phi'(b))
\lesssim
b\phi'(b)\lambda^{1+1/\smphi}
\le
b\phi'(b)\lambda^4
=
\frac{\phi'(b)}{b^3}a^4.
\]
Since $b<a/2$, we have
\[
(b^2-a^2)^2\simeq a^4.
\]
Using \eqref{eq:elem_muprime} we obtain
\[
\phi_a^*(|m(b)-m(a)|a)
\lesssim
|\mu'(b^2)|\,(b^2-a^2)^2.
\]
This proves the estimate in the concave regime.
\myendproof

\begin{theorem}\label{thm:ode_fem_convergence_proof}
Let Assumption~\ref{ass:sindices_strict} and the assumptions of Theorems~\ref{thm:ode_bound} and~\ref{thm:mielke_lambda_discrete} be satisfied. Then there exists a constant $C>0$, depending only on $\smphi$, $\spphi$, and the $\Delta_2$ constants of $\phi$ and $\phi^*$, such
that
\[
\|V(Du_h[c_h(t)])-V(D\bar u_h)\|_\Ltwo \le C e^{-\lambda t}\sqrt{d_h(0)},
\]
where $c_h(t)$ is the solution of the ODE \eqref{eq:ode_fem}, $\bar u_h$ is the solution of \eqref{eq:galerkin_fem} and $\lambda$ is given in Theorem \ref{thm:mielke_lambda_discrete}. Therefore, if $\lambda>0$, we have
\[
\lim_{t\to+\infty}\|V(Du_h[c_h(t)])-V(D\bar u_h)\|_\Ltwo=0.
\] 
\end{theorem}
\mybegproof
Using the notation 
\[
c_t:=c_h(t),\quad u_t:=u_h[c_h(t)],\quad\mathcal{R}_t:= c_t - |D u_t|^2.
\]
we have
\[
\dot c_h(t)=-\mathcal{R}_t,\quad d_h(t)=\langle \mathcal{G}(c_t)\mathcal{R}_t,\mathcal{R}_t\rangle,
\]
where $d_h(t)$ is defined in \eqref{eq:Dht_def}.
Differentiating $d_h(t)$ yields
\[
\begin{aligned}
d_h'(t)&=\langle \mathcal{DG}(c_t)[\dot c_h(t)]\mathcal{R}_t,\mathcal{R}_t\rangle
+
 2\langle \mathcal{G}(c_t)\mathcal{R}_t,\mathcal{DR}_t(c_t)\dot c_h(t)\rangle \\
 &=
 -\langle \mathcal{DG}(c_t)[\mathcal{R}_t]\mathcal{R}_t,\mathcal{R}_t\rangle
 -
 2\langle \mathcal{G}(c_t)\mathcal{R}_t,\mathcal{DR}_t(c_t)\mathcal{R}_t\rangle
 \\
 &=
 -2
 \left(
 \frac12
 \langle \mathcal{DG}(c_t)[\mathcal{R}_t]\mathcal{R}_t,\mathcal{R}_t\rangle
 +
 \langle \mathcal{G}(c_t)\mathcal{R}_t,\mathcal{DR}_t(c_t)\mathcal{R}_t\rangle
 \right)
 \\
 &=
 -2\mathcal{B}(c_t,\mathcal{R}_t),
\end{aligned}
\]
where $\mathcal{B}$ is given in Equation \eqref{eq:Bdef}. Under the assumptions of Theorem \ref{thm:mielke_lambda_discrete}, 
we have the metric Hessian lower bound
\[
\mathcal{B}(c_t,\mathcal{R}_t)\ge \lambda g_{c_t}(\mathcal{R}_t,\mathcal{R}_t)=\lambda d_h(t),
\]
and thus 
\[
d_h'(t)\le -2\lambda d_h(t).
\]
Therefore, by Gronwall's inequality,
\begin{equation}\label{eq:temp_D_h_dec}
d_h(t)\le e^{-2\lambda t}d_h(0).
\end{equation}

Recall that $u_t$ is the solution of the weighted linear problem
\[
\int_\Omega \mu(c_t)Du_t:Dv_h\dx=F(v_h),\quad \forall v_h\in W_h,
\]
and $\bar{u}_h$ the solution of the nonlinear problem \eqref{eq:galerkin_fem}
\[
\int_\Omega A(D\bar{u}_h):Dv_h\dx=F(v_h),\quad \forall v_h\in W_h.
\]
Choosing $v_h=u_t-\bar u_h$ yields
\[
\int_\Omega
A( D\bar{u}_h):(Du_t-D\bar{u}_h)\dx
=
\int_\Omega
 \mu(c_t)Du_t:(Du_t- D\bar{u}_h)\dx.
\]
Multiplying by $-1$ and adding
\[
\int_\Omega A(Du_t):(Du_t-D\bar{u}_h)\dx
\]
to both sides yields
\[
\begin{aligned}
\int_\Omega
(A(Du_t)-A( D\bar{u}_h)):(Du_t-D\bar{u}_h)\dx
 &=
 \int_\Omega
 \bigl(A(Du_t)-\mu(c_t)Du_t\bigr):(Du_t- D\bar{u}_h)\dx
 \\
 &=
 \int_\Omega
 \bigl(\mu(|Du_t|^2)-\mu(c_t)\bigr)
 Du_t:(Du_t-D\bar{u}_h)\dx .
\end{aligned}
\]
From statement (iii) in Lemma \ref{lem:DK-shifted}
\[
 (A(Du_t)-A(D\bar{u}_h)):(Du_t-D\bar{u}_h)
 \simeq
 |V(Du_t)-V(D\bar{u}_h)|^2,
\]
and therefore there exists $c_0>0$ such that
\[
 c_0
 \|V(Du_t)-V(D\bar{u}_h)\|_\Ltwo^2
 \le
 \int_\Omega
 |\mu(|Du_t|^2)-\mu(c_t)|
 |Du_t|
 |Du_t-D\bar{u}_h|\dx .
\]
From the shifted Young inequality in statement (iv) of Lemma \ref{lem:DK-shifted} with shift $a:=|Du_t|$, for every
$\delta>0$,
\[
|\mu(|Du_t|^2)-\mu(c_t)|
 |Du_t|
 |Du_t-D\bar{u}_h|
\le
 \delta
 \phi_{|Du_t|}(|Du_t-D\bar{u}_h|)
 +
 C_\delta
 \phi^*_{|Du_t|}
 \left(
 |\mu(|Du_t|^2)-\mu(c_t)|
 |Du_t|
 \right),
\]
and thus
\[
\begin{aligned}
c_0 \|V(Du_t)-V(D\bar{u}_h)\|_\Ltwo^2
&\le
\delta
\int_\Omega
\phi_{|Du_t|}(|Du_t-D\bar{u}_h|)\dx
+
C_\delta \mathcal{C}_h(t)\\
&\le
C_0\,\delta
\|V(Du_t)-V(D\bar{u}_h)\|_\Ltwo^2
+
C_\delta \mathcal{C}_h(t),
\end{aligned}
\]
where
\[
\mathcal C_h(t)
:=
\int_\Omega
\phi^*_{|D u_t|}
\left(
 |\mu(c_t)-\mu(|Du_t|^2)|
 |D u_t|
 \right)\dx,
\]
and the last step follows from the shifted natural-distance equivalence in statement (iii) of Lemma \ref{lem:DK-shifted}. Choosing $\delta$ sufficiently small, for example $C_0 \delta \le c_0/2$, we obtain
\begin{equation}\label{eq:temp_VDU}
\|V(Du_t)-V(D\bar{u}_h)\|_\Ltwo^2 \le \frac{2C_\delta}{c_0} \mathcal C_h(t).
\end{equation}
By Lemma \ref{lem:coefficient-residual} we have the pointwise estimate
\[
\phi^*_{|Du_t|}
\left(
|\mu(c_t)-\mu(|Du_t|^2)|\,|Du_t|
\right)
\lesssim
|\mu'(c_t)|\,|c_t-|Du_t|^2|^2.
\]
Hence, using \eqref{eq:temp_D_h_dec} 
\[
\mathcal C_h(t) \lesssim
\int_\Omega
|\mu'(c_t)|
(c_t-|Du_t|^2)^2\dx
= 2d_h(t)\lesssim e^{-2\lambda t}d_h(0).
\]
Combining this with \eqref{eq:temp_VDU}  and taking the square root, we conclude the proof.
\myendproof

\begin{lemma}\label{lem:schur}
Under Assumptions \ref{ass:sindices} and \ref{ass:fem_space}, and if $c_{h,n} \in \Ch^+$,  the Schur complement matrix given in Equation \eqref{eq:schur} is positive definite.
\end{lemma}
\mybegproof
For compactness, with the notation from Section \ref{sec:impl_time}, set $\widehat c_h:=c_{h,n+1/2}$, $\widehat u_h:=u_{h,n}$, and $\alpha_n:=\frac{\tau_n}{1+\tau_n}$.
We next identify the bilinear form represented by $S^n$ given in \eqref{eq:schur}. Let
$P_h:\Ltwo\to C_h$ denote the $L^2$-projection onto $C_h$,
defined by
\[
\int_\Omega P_h f\,w_h\dx=\int_\Omega f w_h\dx,\quad \forall w_h\in C_h.
\]
For $z_h\in U_h$, the definition of $B^n$ implies
\[
M^{-1}B^n z_h= P_h(\grad \widehat u_h\cdot \grad z_h).
\]
Therefore, for $z_h,v_h\in U_h$,
\[
\langle S^n z_h,v_h\rangle=
\int_\Omega
\mu(\widehat c_h)\grad z_h\cdot\grad v_h\dx
+
2\alpha_n\int_\Omega\mu'(\widehat c_h)P_h(\grad \widehat u_h\cdot\grad z_h)(\grad \widehat u_h\cdot\grad v_h)\dx.
\]
By Assumption \ref{ass:fem_space}
\[
\grad z_h\cdot\grad v_h=
\frac{1}{2}
\left(
|\grad(z_h+v_h)|^2-|\grad z_h|^2-|\grad v_h|^2\right)\in C_h
\]
for all $z_h,v_h\in U_h$. Therefore  the $L^2$-projection is exact
$P_h(\grad \widehat u_h\cdot\grad z_h)=\grad \widehat u_h\cdot\grad z_h$,
and the Schur complement is represented by the symmetric bilinear form
\[
s^n(z_h,v_h)=\int_\Omega\mu(\widehat c_h)\grad z_h\cdot\grad v_h\dx
+
2\alpha_n\int_\Omega\mu'(\widehat c_h)(\grad \widehat u_h\cdot\grad z_h)
(\grad \widehat u_h\cdot\grad v_h)\dx.
\]

We now prove the positive definiteness and take $z_h=v_h$
\[
s^n(v_h,v_h)=\int_\Omega\mu(\widehat c_h)|\grad v_h|^2\dx
+
2\alpha_n
\int_\Omega\mu'(\widehat c_h)(\grad \widehat u_h\cdot\grad v_h)^2\dx.
\]
When $\Phi$ is convex, $\mu'(\widehat c_h)\geq 0$, and the second term is nonnegative. Therefore, the statement is proved since $\mu(\widehat c_h)> 0$. When $\Phi$ is concave
$\mu'(\widehat c_h)<0$; from \eqref{eq:be_interp} we have the pointwise inequality
$\alpha_n|\grad \widehat u_h|^2 %
\leq
\widehat c_h.$
Therefore, since by Assumption \ref{ass:sindices} 
\[
\mu(r)+2r\mu'(r)=
\frac{\sqrt r\,\phi''(\sqrt r)}{\phi'(\sqrt r)}
\mu(r)\ge \smphi\mu(r)>0,\quad\forall r > 0
\]
we obtain
\[
\begin{aligned}
\mu(\widehat c_h)|\grad v_h|^2+2\alpha_n\mu'(\widehat c_h)(\grad \widehat u_h\cdot\grad v_h)^2 &\ge
\mu(\widehat c_h)|\grad v_h|^2+2\alpha_n\mu'(\widehat c_h)|\grad \widehat u_h|^2|\grad v_h|^2\\
&\ge
\left(
\mu(\widehat c_h)+2\widehat c_h\mu'(\widehat c_h)\right)|\grad v_h|^2\\
&\ge
\smphi\mu(\widehat c_h)|\grad v_h|^2,
\end{aligned}
\]
and thus
\[
s^n(v_h,v_h)\geq \smphi \int_\Omega\mu(\widehat c_h)|\grad v_h|^2\dx >0,
\]
and the statement is proved.
\myendproof

\end{document}